\documentclass[a4paper, 11pt]{article}
\usepackage{amsfonts}
\usepackage{amsmath}
\usepackage{amssymb}
\usepackage{amsthm}

\usepackage{logreq}
\usepackage{tikz}
\usetikzlibrary{shapes,snakes}
\usepackage{paralist}	 
\usepackage{comment}
\usepackage{cite}
\usepackage{placeins}
\usepackage{multicol}
\usepackage{multirow}
\usepackage{caption}
\usepackage{subfigure}
\usepackage[latin1]{inputenc}
\usepackage{graphicx}
\usepackage{epsf}
\usepackage{fancyhdr}
\usepackage{pstricks}
\usepackage{algorithm}
\usepackage{algorithmic}
\usepackage{lscape}
\newpsobject{showgrid}{psgrid}{subgriddiv=1, griddots=10,gridlabels=6pt}
\usepackage{authblk}
\usepackage[capitalise]{cleveref}

\usepackage{environ}
\makeatletter
\newsavebox{\measure@tikzpicture}
\NewEnviron{scaletikzpicturetowidth}[1]{%
  \def\tikz@width{#1}%
  \def\tikzscale{1}\begin{lrbox}{\measure@tikzpicture}%
  \BODY
  \end{lrbox}%
  \pgfmathparse{#1/\wd\measure@tikzpicture}%
  \edef\tikzscale{\pgfmathresult}%
  \BODY
}

\begin{document}

\begin{center} 
{A discrete optimisation approach for target path planning whilst evading sensors} ~\\
~\\
{J.E. Beasley} ~\\ 
~\\
~\\
Mathematics, Brunel University, Uxbridge UB8 3PH, UK 

 ~\\
 ~\\
john.beasley@brunel.ac.uk ~\\
{http://people.brunel.ac.uk/$\sim$mastjjb/jeb/jeb.html}

~\\
~\\
June 2021,~Revised September 2021\\
\end{center}

\begin{abstract}
In this paper we deal with a practical problem that arises in military mission planning. The problem is to plan a path for one, or more, agents to reach a target without being detected by enemy sensors. 

Agents are not passive, rather they can 
initiate actions which aid evasion.  They can knockout sensors. Here to knockout a sensor means to completely disable the sensor. They can also confuse sensors. Here to confuse a sensor  means to  reduce the probability that the sensor can detect an agent.

Agent actions are path dependent and time limited.  By path dependent we mean that an agent needs to be sufficiently close to a sensor to knock it  out. By time limited we mean that a limit is imposed on how long a sensor is knocked out or confused before it reverts back to its original operating state.

The approach adopted breaks the continuous space in which agents move into a discrete space. This enables the problem to be formulated  as a zero-one integer program with linear constraints. The advantage of representing the problem in this manner is that powerful commercial software optimisation packages exist to solve the problem to proven global optimality.  A heuristic for the problem based on successive shortest paths is also presented.

Computational results are presented for a number of randomly generated test problems that are made publicly available.

\end{abstract}


\section{Introduction}

\subsection{Problem outline}
To outline the problem dealt with in this paper consider the scenario of a single agent who has to navigate their way through enemy territory to reach a target location. Clearly the agent is concerned with reaching the target location as quickly as possible, whilst evading detection by the enemy. 

Detection of the agent by the enemy is typically achieved by means of electronic sensors. Naturally any individual sensor will have a limited area within which it can detect an agent. However a number of such   sensors, especially if linked together into a connected sensor network, have the potential to provide the enemy with wide geographic coverage to aid in detecting any intruding agent and hence defend the  target location. 

A natural question is: why not simply take out the target location with a missile? This may  be inappropriate for a number of reasons. For example to avoid collateral damage, or due to the nature of the  mission. Missions such as hostage rescue,  or infiltration to gather intelligence from captured personnel and/or equipment, are clear examples where  a  missile strike would be inappropriate. 

In the problem considered in this paper the agent has two actions that they can take which will help them evade detection. 
The first such action is that the agent has the   ability to knockout sensors. Here to knockout a sensor means to completely disable the sensor.
The second such action is that the agent has the ability to confuse sensors. Here to confuse a sensor  means 
that the probability of a sensor detecting an agent, when the agent is within sensor detection range, is reduced.

The decision problem dealt with in this paper therefore is to decide the path the agent should adopt, and the appropriate knockout/confusion actions to take, so as to reach the target location whilst best avoiding detection.

\subsection{Approach taken}

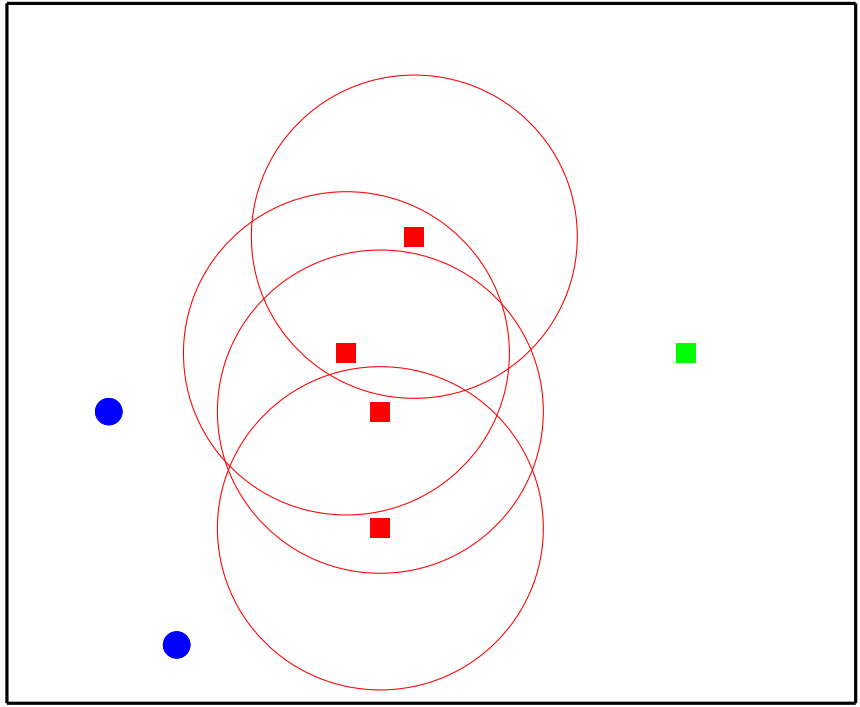
\begin{figure}[!htb]
\centering	  
\begin{scaletikzpicturetowidth}{\textwidth*0.70} 
\begin{tikzpicture}[scale=\tikzscale]	
\draw [red] (   5.50000 ,    2.59808) circle [radius=   2.40000];
\draw [red] (   5.00000 ,    5.19615) circle [radius=   2.40000];
\draw [red] (   6.00000 ,    6.92820) circle [radius=   2.40000];
\draw [red] (   5.50000 ,    4.33013) circle [radius=   2.40000];

\node at (  5.50000 ,    2.59808) [rectangle    , fill=red    ] {};
\node at (  5.00000 ,    5.19615) [rectangle    , fill=red    ] {};
\node at (  6.00000 ,    6.92820) [rectangle    , fill=red    ] {};
\node at (  5.50000 ,    4.33013) [rectangle    , fill=red    ] {};

\draw[very thick] (    .00000 ,     .00000) -- (  12.50000 ,     .00000);
\draw[very thick] (  12.50000 ,     .00000) -- (  12.50000 ,   10.39230);
\draw[very thick] (  12.50000 ,   10.39230) -- (    .00000 ,   10.39230);
\draw[very thick] (    .00000 ,   10.39230) -- (    .00000 ,     .00000);
\node at (   2.50000 ,     .86603) [circle    , fill=blue     ] {};
\node at (   2.50000 ,     .86603) [circle    , fill=blue     ] {};
\node at (   1.50000 ,    4.33013) [circle    , fill=blue     ] {};
\node at (   1.50000 ,    4.33013) [circle    , fill=blue     ] {};
\node at (  10.00000 ,    5.19615) [rectangle    , fill=green    ] {};
\end{tikzpicture}
\end{scaletikzpicturetowidth}
\caption{Example problem}
\label{fig1}
\end{figure}

To illustrate the approach taken to the problem dealt with in this paper  consider Figure~\ref{fig1}. 
In that figure we have two agents, shown as solid blue circles. We want one or other of these agents to reach the target, shown as a solid  green square. The complication is that there are four sensors, shown as solid red squares. Each sensor has  an associated circular area, shown as red circles, such that within a circle an agent can be detected by the associated sensor. 
The problem is to reach the target in a desired target time, but avoiding detection by the sensors.

In the context of the problem considered in this paper,  detection of an agent by a single sensor is not regarded as relevant. Rather detection of the agent is  \textbf{\emph{only}} relevant if the agent is detected by two or more sensors. This relates to the fact that typically detection by two or more sensors is needed to locate the position of the agent.

Agents are not passive, rather they can initiate actions which aid evasion. These are knockout and confusion.  Sensor knockout can be achieved, for example, by means of a kinetic or non-kinetic   attack. Sensor confusion, also referred to as sensor degradation, can be achieved, for example, by
electronic countermeasures.

Agent actions are path dependent and time limited.  By path dependent we mean that an agent needs to be sufficiently close to a sensor to knock it  out. By time limited we mean that a limit is imposed on how long a sensor is knocked out or confused before it reverts back to its original operating state.

\begin{figure}[!htb]
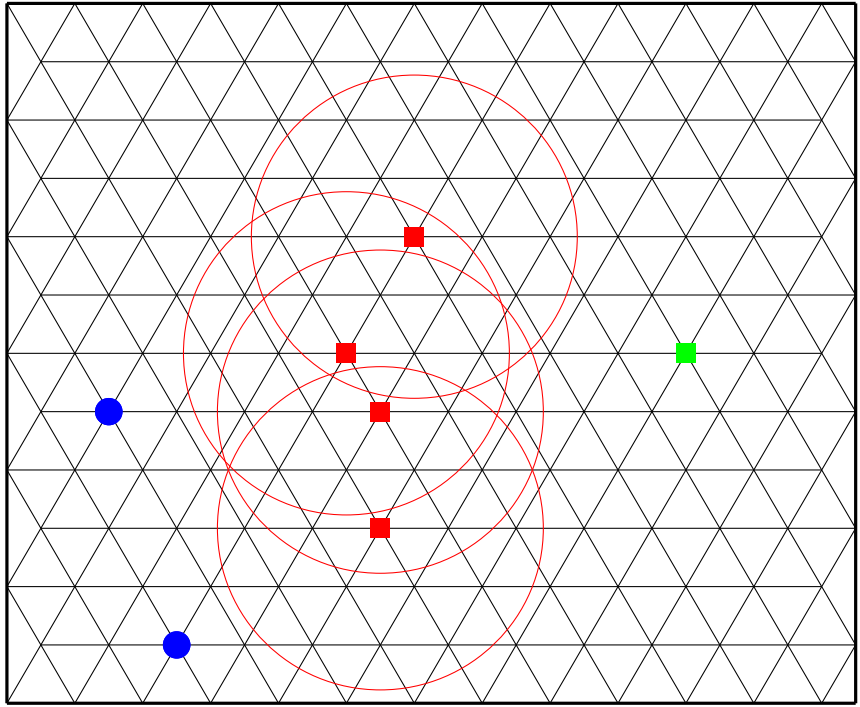

\centering	  
\begin{scaletikzpicturetowidth}{\textwidth*0.70} 

\end{scaletikzpicturetowidth}
\caption{Example problem: triangular mesh}
\label{fig2}
\end{figure}

We will discretise the problem both in space and time. Consider the  mesh graph shown in 
Figure~\ref{fig2} 
placed over the two-dimensional space shown in 
Figure~\ref{fig1}. Agents move from vertex to vertex on this graph, where a movement from one vertex to an adjacent vertex takes one time step. The aim is to have one, or more, agents reach the target vertex in the desired target time without detection by two or more sensors. 

Whilst conceptually any mesh is possible, the requirement that agents move from vertex to vertex necessities using a mesh such that the maximum time taken for an agent to move between any two adjacent vertices is 
the same. Hence, we adopt here a triangular mesh such as that  shown in 
Figure~\ref{fig2}. In that figure, we have a mesh composed of equilateral triangles.  In any known practical application  the mesh would be set having regard to geographic constraints, e.g.~rivers, buildings, etc. This is easily done - simply impose an arbitrary triangular mesh over the area of interest and then delete from the mesh any vertices that are not reachable due to geographic constraints. Note here  with reference to Figure~\ref{fig2} that, although the agents move between vertices, and the target is also at a vertex, there is no requirement for the sensors to be positioned at vertices.

With regard to the practical military problem which underlies the problem considered in this paper, an archetypal example would be that each agent represents a small force of ground troops who have to move through enemy territory to reach a known target location whilst evading detection by enemy sensors.

It is important to stress here that we do not envisage the use of the formulation given in this paper to be such that it gives in precise detail the path to be followed by agents. Rather we envisage the work reported here to be a decision aid for mission planners that enables them to generate a small set of options for more detailed mission planning.

The structure of this paper is as follows. 
In Section~\ref{litsur} we present our literature survey and state what we believe are the main contributions of this paper to the literature.
In Section~\ref{sec2} we present the constraints associated with agent movement. We indicate how making use of shortest path calculations can significantly reduce the number of decision variables required. We also indicate how to deal with the situations where we have more agents than we need and with agent exit strategy.

In Section~\ref{sec3a} we present the constraints associated with agent detection when agents can initiate sensor knockouts.  We discuss how a solution can be found when no knockouts are allowed and give computational results for an example problem.

In Section~\ref{sec4} we discuss introducing detection probabilities into the problem. We also discuss how we might combine detection probabilities with knockout.
In Section~\ref{sec5} we consider sensor confusion actions, as well as the relationship between the problem considered in this paper and the biobjective shortest path problem. A heuristic for the problem considered here based on successive shortest paths is also presented.  In Section~\ref{sec6} we discuss some further relevant restrictions on agent actions. 
Section~\ref{sec7} presents computational results 
for randomly generated test problems that are made publicly available and Section~\ref{sec8} presents our conclusions.

\section{Literature survey} \label{litsur}

In this section we present a  survey discussing papers relevant to the problem considered here. It is important for the reader to  be aware here that because the sensor-related papers discussed below are typically directly motivated by an underlying practical problem there is no single underlying problem for which a number of different solution approaches have been developed. Rather there are a set of papers considering a  specific sensor-related problem. This contrasts to more classical Operational Research areas where there is a single underlying problem 
and academic work focuses on developing more effective solution algorithms for that problem.

One problem considered in the literature relates to 
pursuit-evasion, where the problem is to determine a strategy for a team of pursuers, equivalently mobile sensors, to capture an evader. 
Croft~\cite{croft64} provides an accessible early reference to a pursuit-evasion problem entitled: Lion and Man. In that problem the lion and man are constrained to be in a circular arena, moving with equal maximum speeds, and the question is: can the lion catch the man in finite time? In the problem considered there both time and space were continuous. Sgall~\cite{sgall01} dealt with the problem when time is discrete and space is continuous.

DelBalzo and Hemsteter~\cite{del02}  considered a problem motivated by an underwater vessel evading searchers equipped with sonar. In their problem the underwater vessel  can deploy evasion tactics, e.g.~move in a different direction, on detection of searcher sonar. Their focus is on the appropriate search strategy for the sensors to maximise their cumulative detection probability. They proposed a genetic algorithm to solve this problem.
Yan et al~\cite{yan2013} considered a pursuit-evasion game with two stages. In the first stage the evader is always static, in the second stage  the evader moves once it detects pursuers. They proposed a control theory based solution procedure.

Krishnamoorthy et al~\cite{krish13} considered the problem of the optimal strategy for a pursuer searching for an evader on a Manhattan grid network. Their work was motivated in the context of a road network with sensors at selected junctions which detect the evader as it passes. The pursuer accesses this sensor information only as they reach each sensor junction. They considered the worst case time to capture the evader and the corresponding pursuit policy. 
Adams and Carlsson~\cite{adams2015} considered ball-shaped sensors moving randomly and continuously in a bounded domain. Sensors know when they overlap with another sensor. A moving intruder has to avoid detection by the sensors. Their work focused on the necessary and sufficient conditions required for the intruder to evade the sensors.

Karabulut et al~\cite{karabulut17} presented a bilevel formulation of the problem of determining the best sensor locations so as to maximize coverage intensity, and hence increase the probability of detecting an intruder. In their approach an intruder can destroy sensors so as to generate a path evading sensors, so a minimal exposure path. A number of heuristics to solve the problem were presented and computational results given. Note here however that in Karabulut et al~\cite{karabulut17} sensor destruction is represented by a continuous variable (known as the partial interdiction case), rather than a binary variable (known as the complete interdiction case), as in the work presented in this paper. Another difference between the work presented here and Karabulut et al~\cite{karabulut17} is  that in our work sensor destruction is agent path dependent, i.e.~an agent needs to be sufficiently close to a sensor to destroy it.

Lessin et al~\cite{lessin18} presented a bilevel formulation of the problem of determining the best sensor locations to detect an intruder following a path through the sensed area. They transformed the bilevel problem into a mixed-integer linear program solved using a standard solver.  Computational results were presented.
Lessin et al~\cite{lessin19} presented a multi-objective bilevel formulation of the problem of how to relocate sensors from their original locations in the event that some sensors have been incapacitated and others degraded. They presented a reformulation of the problem to identify the Pareto frontier representing objective tradeoffs. Computational results were presented.

Craparo et al~\cite{craparo2019} considered the problem of where to place sources and receivers in a sonar senor network so as  monitor a known set of locations. They presented two integer linear programs as well as a two-step heuristic and an iterative approach. Computational results were given. 
Delavernhe et al~\cite{Delavernhe2020} considered the problem of activating sensors to record the passage of moving targets, but where some uncertainty as to target movement exists. A formulation of the problem with a number of different objectives was given, together with a solution approach and computational results.

Fauske et al~\cite{fauske20} considered the problem of finding routes for force elements involved in maritime surveillance operations. In their problem different force elements have different sensors with different properties, e.g.~different detection ranges. They formulated the problem as a variant of the periodic vehicle routing problem with an objective related to maximising the number of small areas/cells sensed. They used column generation with both exact and heuristic column generation pricing. Computational results were presented.
Delavernhe et al~\cite{Delavernhe2021} considered the problem of moveable sensors searching for a moveable  target. They considered the cost of sensor travel between different areas and presented a formulation of the problem with linear constraints but a nonlinear objective.  A heuristic algorithm for the solution of the problem was presented and computational results given.

In  light of the literature survey given above we believe that the main contributions of this paper to the literature are as follows. Firstly, 
to formulate, as a zero-one integer program with linear constraints, the problem of finding paths for agents to reach a desired target location in minimal time.
Secondly, to include in the formulation agent actions which are path dependent and time limited, namely sensor knockout 
and sensor confusion. 
By path dependent we mean that an agent needs to be sufficiently close to a sensor to knock it  out. By time limited we mean that a limit is imposed on how long a sensor is knocked out or confused before it reverts back to its original operating state.
Thirdly, to present shortest path based constraints which considerably reduce the size of the problem needed to be solved.

\section{Agent movement} \label{sec2}								 
\subsection{Notation}
Our notation is as follows.
We have a graph with $N$ vertices, 
where $\Gamma(v)$ is the set of vertices adjacent to vertex $v$ in this graph ($v \notin \Gamma(v)$).  
We have $A$ agents, where agent $a$ starts at time 0 at vertex $\gamma(a)$.
The target is vertex $N$ and
the  number of time steps (time horizon) allowed is $T$, so we aim to  reach the target at time $T$.
We introduce an artificial vertex, vertex $0$.
This vertex acts as an absorbing vertex, so $\Gamma(0)=\emptyset$, and note that all agents can make a transition to vertex $0$ from any other vertex $v \neq 0$ at any time $t \geq 1$. For ease of explanation below however, note that $0 \notin \Gamma(v)$ and we introduce explicit terms for a transition to vertex $0$.  
The variables are
$x_{avt}=1$ if agent $a$ is at vertex $v$ at time $t$, $= 0$ otherwise.

\subsection{Movement constraints}
We have the following constraints which govern the movement of agents at each time step:
\begin{equation}
\sum_{v=0}^N x_{avt} = 1 ~~~~ a=1,\ldots,A; ~t=1,\ldots,T
\label{eq1}
\end{equation}

\begin{equation}
x_{av,t+1}  \leq x_{avt} + \sum_{k: ~v \in \Gamma(k)} x_{akt}  ~~~~a=1,\ldots,A; ~v=1,\ldots,N; ~t=0,\ldots,T-1
\label{eq2}
\end{equation}

\begin{equation}
x_{av,t+1} + x_{a0,t+1}+ \sum_{k \in \Gamma(v)} x_{ak,t+1}  \geq x_{avt} ~~~~a=1,\ldots,A; ~v=1,\ldots,N; ~t=0,\ldots,T-1
\label{eq2a}
\end{equation}

\begin{equation}
 x_{a\gamma(a)0}=1~~~~a=1,\ldots,A
\label{jebi1}
\end{equation}
\begin{equation}
 x_{av0}=0~~~~a=1,\ldots,A; v=0,\ldots,N; v \neq \gamma(a)
\label{jebi2}
\end{equation}

\begin{equation}
x_{a0,t+1} \geq x_{a0t}
~~~~ a=1,\ldots,A; ~t=1,\ldots,T-1
\label{add2}
\end{equation}

\begin{equation}
\sum_{a=1}^A x_{aNT} \geq 1
\label{eq4aaa}
\end{equation}

Equation~(\ref{eq1}) ensures that each agent is at some vertex at each time step.
Equation~(\ref{eq2}) ensures that at time $t+1$ an agent can only be at vertex $v$ if at the previous time step $t$ it was either already at vertex $v$, so has not moved, or was at some vertex $k$ from which it could reach vertex $v$ in one time step. 

Equation~(\ref{eq2a}) ensures that if at time $t$ an agent was at vertex $v$ then at time $t+1$ the agent is either at vertex $v$, so has not moved, or has moved to vertex $0$ or is at some vertex $k$ which could have been reached from  vertex $v$ in one time step. Note here that Equation~(\ref{eq2}) and Equation~(\ref{eq2a}) allow an agent to remain at the same vertex, so be stationary for one or more time steps.

Equations~(\ref{jebi1}) and (\ref{jebi2}) are the initial conditions.
Equation~(\ref{add2}) ensures that once an agent reaches vertex $0$ they remain there.
Equation~(\ref{eq4aaa})  ensures that at least one agent reaches vertex $N$ at time $T$.

Collectively the constraints above ensure that at each time step an agent is only at one vertex, and moreover that vertex must be one which is reachable from the vertex where the agent was at the previous time step. Hence the path for an agent must consist of a succession of adjacent vertices. 

Note  that we have not yet given an objective function above. This is because above we are treating the problem as a feasibility problem, namely find a path such that one, or more, agents reach the target $N$ at time $T$. There are essentially two reasons as to why treating the problem as a feasibility problem may be 
relevant.
Firstly, the time at which agents reach the target may be specified to   align with other military operations that are also taking place.  
Secondly, specifying the time $T$ at which agents reach the target and varying $T$ will enable the decision-maker to explore how the path changes as the time at which agents reach the target changes, and enable them to consider military considerations that lie outside the model formulated

However if,
in the absence of military considerations,
 we wish to directly minimise the time to target within a time horizon of $T$ time steps 
then we replace Equation~(\ref{eq4aaa})  by
a constraint which ensures that at least one agent reaches vertex $N$ within the time horizon $T$:
\begin{equation}
\sum_{a=1}^A \sum_{t=1}^T x_{aNt} \geq 1 
\label{eq4}
\end{equation}

In addition we introduce  a constraint that forces a transition to vertex $0$ from vertex $N$: 
\begin{equation}
x_{a0,t+1} \geq x_{aNt} ~~~~a=1,\ldots,A; ~t=1,\ldots,T-1
\label{eq3}
\end{equation}
Equation~(\ref{eq3}) ensures that if agent $a$ is at vertex $N$ at time $t$ then the agent is at vertex $0$ at time $t+1$. It will remain at vertex $0$ thereafter as it cannot transition away, see Equation~(\ref{add2}).

Our objective function then is:
\begin{equation}
\mbox{Minimise} ~\sum_{a=1}^A \sum_{t=1}^T tx_{aNt}
\label{eq5}
\end{equation}
In Equation~(\ref{eq5}) we minimise the sum of the times at which each agent reaches the target. Because Equation~(\ref{eq4}) means that at least one agent reaches the target the minimisation here will mean that once one agent reaches the target all agents will make a transition to vertex $0$ signifying that the target has been reached.
Equation~(\ref{eq5}) will therefore  in the optimal solution have only one nonzero term corresponding to the time at which some agent reaches the target.

\subsection{Reducing the number of variables}
A key insight which enables us to considerably improve the computational performance of the formulation given here relates to reducing the number of $x_{avt}$ variables that we need to consider. This reduction is important because there are potentially a large number of such variables, more precisely  $O(ANT)$ variables. 
In addition, reducing the number of such variables has reduction implications for the expansion of the formulation which considers knockout and detection probabilities as given later below.

Suppose we use a standard shortest path algorithm, e.g.~Dijkstra~\cite{dijkstra59}, to compute $D_{ij}$, the minimal number of time steps, ignoring sensors, between vertices $i$ and $j$ in the underlying graph, c.f.~Figure~\ref{fig2}.
Note here that $D_{ij}$ is symmetric, so $D_{ij}=D_{ji}$.
Then we have:
 \begin{equation}
x_{avt} = 0 ~~~~D_{\gamma(a)v} > t; ~a=1,\ldots,A;~v=1,\ldots,N;~t=1,\ldots,T
\label{reduction1}
\end{equation}
\begin{equation}
x_{avt} = 0 ~~~~D_{\gamma(a)v} \leq t;~t + D_{vN} >  T; ~a=1,\ldots,A;~v=1\ldots,N;~t=1,\ldots,T
\label{reduction2}
\end{equation}
Equation~(\ref{reduction1}) ensures that $x_{avt}$ is zero if an agent starting at vertex $\gamma(a)$ at time zero cannot reach vertex $v$ within $t$ time steps. Equation~(\ref{reduction2}) ensures that $x_{avt}$ is zero if an agent starting at vertex $\gamma(a)$ at time zero and going via vertex $v$, so can reach vertex $v$ by time $t$, cannot reach the target vertex $N$ within the time horizon $T$ under consideration.

For ease of presentation we have used explicit constraints, 
Equations~(\ref{reduction1}) and (\ref{reduction2}), to indicate that a particular variable $x_{avt}$ plays no active part in the solution to the problem.  However, in a computational implementation  these would not be included as explicit constraints, but rather any variable $x_{avt}$ satisfying these constraints would be excluded from consideration.

\subsection{Excess agents and agent exit strategy}

\subsubsection{Excess agents}
In our approach we have $A$ agents. The formulation given above has a high degree of degeneracy in that we may have many  solutions, but with a different structure.  So for example if one agent $a$ is used to reach the target at time $T$, then any solution with all other agents making random walks beyond sensor range, before making a transition to vertex zero as agent $a$ reaches the target, will be equivalent solutions. 

In order to eliminate such solutions we force any agents that never reach the target to make a transition to vertex $0$ at time $1$. In other words our formulation explicitly deals with the situation where if we have an excess of agents then we only make use of the number that we need. The constraint that we impose is:
\begin{equation}
x_{a01} \geq 1 - \sum_{t=1}^T x_{aNt} ~~~~ a=1,\ldots,A
\label{eqexcess0}
\end{equation}
Equation~(\ref{eqexcess0}) ensures that $ x_{a01}$ is forced to be one if agent $a$ never reaches the target vertex $N$. 

Note here that excess agents may arise because, in reality, we have only one agent but a number of different possible starting positions for that agent. In the formulation given above each such starting position corresponds to a different agent. In this situation we need to choose just one agent/starting position and to do this we add the constraint $\sum_{a=1}^A x_{a01} = (A-1)$ to ensure that all but one agent moves directly to vertex $0$ at time $1$.

\subsubsection{Agent exit strategy}

As mentioned previously above 
an archetypal example of the problem considered in this paper would be that each agent represents a small force of ground troops who have to move through enemy territory to reach  a target.

As formulated above the goal is for one or more agents to reach the target at time $T$. However in some situations the agents, having reached the target, need an exit strategy. An example exit strategy would be for the agents to make their escape by moving to an extraction   location. This situation can be easily dealt in our formulation with minor modifications. The principal modifications needed are as follows.
Firstly, to regard the vertex $N$ as the extraction location, and add a constraint ensuring that the vertex $\xi$, at which the target 
is located, is visited, i.e.~add $\sum_{a=1}^A \sum_{t=1}^T x_{a \xi t} \geq 1$.
Secondly, to alter the shortest path reduction constraints, Equations~(\ref{reduction1}) and (\ref{reduction2}), to account for the fact that we visit vertex $\xi$ before vertex $N$.

\section{Agent detection} 
\label{sec3a}

We assume that detection of an agent only occurs at vertices, so not during movement between vertices, and hence at time step $t$ detection, if any, occurs at the vertex to which the agent has just moved. In other words detection is based upon the values of $x_{avt}~t \geq 1$. For ease of explanation we also assume here that at time 0 the agents are positioned so as to be undetectable. 

With reference to terminology here note that below we refer to detection of an agent. This, depending upon the context, may refer to detection of an agent by a  single sensor, or detection of an agent by two or more sensors.

Detection of any agent by two or more sensors is of importance in this paper since the assumption here is that combining detection information from two or more sensors enables the operators of the sensor network to achieve much greater accuracy as to the position of the detected agent. An analogy here would be that if each sensor gives a precise bearing to a detected agent then detection by two sensors would be required to precisely and accurately position the agent, i.e.~in two-dimensional Euclidean space to deduce from  two sensor bearings the $(x,y)$ position of the agent. Torrieri~\cite{torr84} gives technical details as to how to estimate the location of an agent given imprecise sensor information.

We assume that
detection is binary, so either a sensor detects an agent or not. 
We also assume that each agent  can  knockout sensors to render them inactive, so incapable of detecting any agents.
The case where sensors have a probability of agent detection is considered later below.

\subsection{Notation and assumptions}
The framework for our formulation in terms of notation and assumptions is as follows.
We have $S$ sensors with fixed positions at known locations.
We set $d_{sv}=1$ if sensor $s$ can detect an agent at vertex $v$ under normal operating conditions, $=0$ otherwise. Without significant loss of generality we assume that $d_{sN}=0,~ s=1,\ldots,S$ so no sensor can detect an agent at the target vertex $N$. Similarly $d_{s0}=0,~ s=1,\ldots,S$. Here $d_{sv}$ captures the range of sensor $s$ in terms of which vertices it covers, i.e.~the vertices at which it can detect an agent under normal operating conditions. Note in passing here that $d_{sv}$ is not dependent on each sensor having a circular detection area, such as shown in Figure~\ref{fig2}. As such the detection area for each sensor can be of arbitrary shape, and different for different sensors.

Let $K_{sv}=1$ if an agent at vertex $v$ can knockout sensor $s$, i.e.~it is in knockout range, $=0$ otherwise. Here $K_{sv}$ captures for agents their knockout range in terms of which sensors an agent can knock out when positioned at vertex $v$. Here we assume that each agent has the same knockout range and can knockout out all sensors within range. With regard to this assumption then although knockout can be achieved by kinetic measures note that if we are using electromagnetic measures to achieve knockout then we have simplified the problem here.  Essentially we are assuming a range based propagation model of undirected electromagnetic energy. To achieve the desired knockout effect this electromagnetic energy has to arrive at a sensor with sufficient power. So our simplification is to assume omnidirectional electromagnetic propagation through a uniform environment, where the delivered power decreases with respect to distance from the agent.

We let $C_a$ be the cost of a single knockout action by agent $a$, with the total knockout action cost being limited by 
$B$. We define $\Delta_a$ as the number of time steps for which a sensor is disabled once it has been first knocked out by agent $a$. Note here that an agent knocking out a sensor renders the sensor completely inoperative, so it has no detection ability at all, for $\Delta_a$  time steps.

No agent can be detected by $\Omega$ or more sensors at any time, where typically we are interested in $\Omega=2$, as discussed above. 
Since detection of an agent relies on at least $\Omega$  sensors we automatically know that any vertex $v$ which not in range of $\Omega$  or more sensors is of no detection interest. Hence in terms of detection we need only concern ourselves with the set of vertices $V$ where $V \subseteq [v~|~ v=1,\ldots,N-1]$ is the set of vertices which are in range of $\Omega$, or more, sensors, i.e.~$V=[v~|~ \sum_{s=1}^S d_{sv} \geq \Omega~ v=1,\ldots,N-1]$.

\subsection{Variables}
The variables are as follows.
Let $\alpha_{avt}=1$ if agent $a$ at vertex $v$ at time $t$ chooses to knockout all sensors that are in range, $= 0$ otherwise.
Let $\lambda_{st}=1$ if sensor $s$ is no longer active at time $t$ because it is in a knocked out state, $= 0$ otherwise.
Let $y_{sat}=1$ if sensor $s$ detects agent $a$ at time $t$, $= 0$ otherwise.

Note there that although we need only restrict potential detections  to vertices $v \in V$ it is possible that performing a knockout at some vertex $v \notin V$ would be worthwhile due to the position of that vertex and the sensors within knockout range. For this reason the variables $\alpha_{avt}$ need to be defined over the entire set of vertices.

In our formulation we assume that knockout has  priority over detection. Recall here that detection only occurs at vertices. If at time $t$ an agent $a$ has just moved to vertex $v$, so $x_{avt}=1$, and that agent then initiates a knockout action, so $\alpha_{avt}=1$, then there can be no detection of that agent at that vertex by any sensor affected by the knockout action.

\subsection{Detection constraints}

The constraints associated with agent detection are as follows:

\begin{equation}
\alpha_{avt} \leq x_{avt} 
~~~~ a=1,\ldots,A; ~v=1,\ldots,N;~t=1,\ldots,T
\label{for1}
\end{equation}

\begin{equation}
\alpha_{avt} = 0 
~~~~ \sum_{s=1}^S K_{sv}=0;~a=1,\ldots,A; ~v=1,\ldots,N;~t=1,\ldots,T
\label{alphared}
\end{equation}

\begin{equation}
\sum_{a=1}^A \sum_{v=1}^N \sum_{t=1}^T C_a \alpha_{avt} \leq B 
\label{deq1}
\end{equation}

\begin{equation}
\lambda _{st} \leq \sum_{j=1}^N K_{sj} \sum_{a=1}^A ~\sum_{\tau=max(1,t-\Delta_a)}^t \alpha_{aj\tau}
~~~~s=1,\ldots,S;~t=1,\ldots,T
\label{eqlam1}
\end{equation}

\begin{equation}
\lambda _{st} \geq [\sum_{j=1}^N K_{sj} \sum_{a=1}^A ~\sum_{\tau=max(1,t-\Delta_a)}^t \alpha_{aj\tau}]/M
~~~~s=1,\ldots,S;~t=1,\ldots,T
\label{eqlam2}
\end{equation}

\begin{equation}
y_{sat} \geq x_{avt} - \lambda_{st}
~~~~d_{sv}=1;~s=1,\ldots,S;~a=1,\ldots,A;~\forall v \in V;~t=1,\ldots,T
\label{deq2}
\end{equation}

\begin{equation}
\sum_{s=1}^S  y_{sat} \leq \Omega - 1
~~~~a=1,\ldots,A;~t=1,\ldots,T 
\label{deq3}
\end{equation}

Equation~(\ref{for1}) ensures that we cannot perform a knockout unless an agent is at a vertex at the appropriate time.  
Any agent performing a knockout at a vertex $v$ that cannot affect any sensor at all is clearly irrelevant and this leads to Equation~(\ref{alphared}).
Equation~(\ref{deq1}) limits the total knockout cost.

If sensor $s$ is no longer active at time $t$ because it is in a knocked out state then this must because it has been knocked out by some agent, either at time step $t$ or in one of the preceding 
time steps. This is dealt with using Equation~(\ref{eqlam1})  and Equation~(\ref{eqlam2}). 

In Equation~(\ref{eqlam1}) we sum over all vertices $j$ that are capable of knocking out sensor $s$, i.e.~with $K_{sj}=1$. For each such vertex some agent $a$ may have initiated a knockout action at some time $\tau$ within the appropriate time interval ending at time $t$. 
Equation~(\ref{eqlam1}) ensures that $\lambda _{st}$ is zero if no such knockout actions have taken place. 

In Equation~(\ref{eqlam2}) $M$ is a positive constant and the equation ensures that  $\lambda _{st}$ is forced to be one if any term in the numerator on the right-hand side of the inequality is non-zero, i.e.~a knockout action for sensor $s$ has occurred. Considering Equation~(\ref{deq1}) a suitable value for $M$ here is $M=B/\mbox{min}[C_a~|~C_a >0~ a=1,\ldots,A]$. Equations~(\ref{eqlam1}) and (\ref{eqlam2}) together ensure that  $\lambda _{st}$  is forced to be zero or one as appropriate, with both constraints being required to achieve this, i.e.~neither constraint is redundant.

Note that, for ease of formulation, we have assumed here in
Equations~(\ref{eqlam1}) and (\ref{eqlam2})
that if knockout is initiated at some time $\tau < t$ it continues to be active at time $t$. If knockout is achieved by non-kinetic means such as electronic jamming, for example, then this assumes that jamming continues for sensor $s$ after first being initiated at time $\tau$ by some agent $a$ at some vertex $j$ in knockout range of sensor $s$. Since knockout will typically be initiated by an agent as they move through the detection area of a sensor this seems a reasonable assumption.

Now at vertex $v \in V$ an agent $a$ can only be detected by sensor $s$ at that vertex at time step $t$ if $x_{avt}=1$. However that detection is negated if sensor $s$ has been knocked out. This therefore leads to Equation~(\ref{deq2}). 
To clarify Equation~(\ref{deq2}), first note that this constraint only applies if $d_{sv}=1$, i.e.~if sensor $s$ can under normal operating conditions detect an agent at the vertex $v$ under consideration on the right-hand side of Equation~(\ref{deq2}). 

Equation~(\ref{deq2}) is inactive if
$x_{avt}=0$. Hence
Equation~(\ref{deq2})
can only force $y_{sat}$ to be one for agent $a$ at time $t$ if we have a vertex $v$ for which
$d_{sv}=1$ and $x_{avt}=1$ (recall here that each agent can only be at one vertex at each time step).
In the right-hand side of Equation~(\ref{deq2}) the second term will be non-zero, so rendering the constraint inactive, if sensor $s$ has been knocked out at time $t$.

Each agent must remain undetected at all times and this means that no agent can be detected by $\Omega$ or more sensors at any time. The constraint which ensures this is Equation~(\ref{deq3}). 

Note here that technically Equation~(\ref{deq2}) forces $y_{sat}$ to be one if the right-hand side is one, but places no restriction on $y_{sat}$ if the right-hand side is zero. However since Equation~(\ref{deq3}) restricts the sum of the $y_{sat}$ then $y_{sat}$ will in many cases be automatically assigned a value of zero to satisfy this equation. 

However, we can regard dealing with values assigned to $y_{sat}$ of one that are not required to be one when considering the values attained by the right-hand side of Equation~(\ref{deq2}) as simply a post-processing exercise. We simply change such values to zero without affecting the feasibility, or optimality, of the solution. Hence the formulation remains valid.

Note here that with agent knockout actions added to the problem the previous constraint, Equation~(\ref{eqexcess0}), dealing with excess agents changes to:
\begin{equation}
x_{a01} \geq 1 - \sum_{t=1}^T x_{aNt} - \sum_{v=1}^N \sum_{t=1}^T
\alpha_{avt}~~~~ a=1,\ldots,A
\label{eqexcess1}
\end{equation}
Equation~(\ref{eqexcess1}) ensures that $ x_{a01}$ is forced to be one if agent $a$ never reaches the target vertex $N$, or engages in a knockout action.

In the formulation given above we allow for multiple agents, so $A \geq 2$ agents. The reason for this is that it easy to conceive of practical situations where we need more than one agent. For example suppose the area the agents are traversing has a geographic barrier, such as a river, which impedes travel. In such a situation there may be a sensor on one side of the river that one agent needs to be assigned to knockout, to enable a second agent to travel without detection on the other side of the river.

\subsection{Solution, $B=0$}
Above we discussed how calculating the shortest path from each agent to the target could be used in problem reduction. In fact we can use a shortest path calculation to directly calculate a solution for $B=0$, i.e.~no knockout actions.

Recall that, as above, no agent can be detected by $\Omega$ or more sensors at any time. Find the shortest path for agent $a$ from their starting position to the target subject to the condition that the path cannot go through any vertex that can be detected by $\Omega$ or more sensors. In other words we avoid all vertices $v \in V$, which by definition satisfy
$\sum_{s=1}^S d_{sv} \geq \Omega$.
Then the maximum number of sensors which can detect an agent at any vertex on this path will be $\Omega-1$. 

Performing this shortest path calculation for each agent and taking the minimal path over all agents provides the minimal time to target solution for the case $B=0$. In this specific case there is no need to resort to the formulation given above for knockout.

\subsection{Example problem}

To provide insight into the formulation given above we give in this section some computational results for this formulation on the example problem shown in Figure~\ref{fig2}. These results were produced using the optimisation package Cplex~\cite{cplex1210}. 

To produce the results seen in this section we set the cost of a knockout action $C_a$ to one and so $B$, Equation~(\ref{deq1}),
corresponds to the total number of knockout actions allowed. Each agent must never be detected by two or more sensors, so detection by just a single sensor is allowed. This corresponds to $\Omega=2$, Equation~(\ref{deq3}).

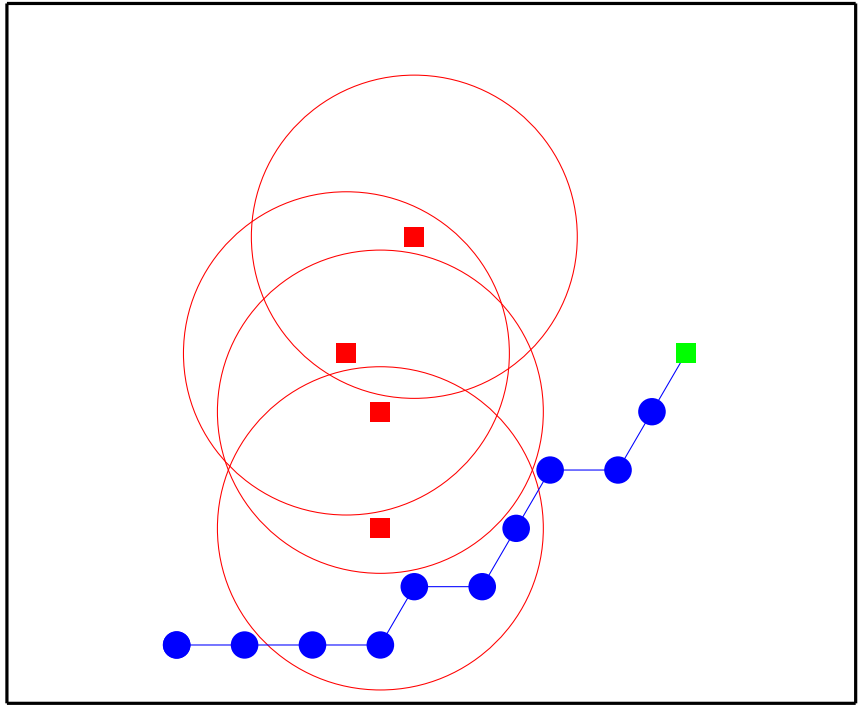
\begin{figure}[!htb]
\centering	  
\begin{scaletikzpicturetowidth}{\textwidth*0.70} 
\begin{tikzpicture}[scale=\tikzscale]	
\draw [red] (   5.50000 ,    2.59808) circle [radius=   2.40000];
\draw [red] (   5.00000 ,    5.19615) circle [radius=   2.40000];
\draw [red] (   6.00000 ,    6.92820) circle [radius=   2.40000];
\draw [red] (   5.50000 ,    4.33013) circle [radius=   2.40000];

\node at (  5.50000 ,    2.59808) [rectangle    , fill=red    ] {};
\node at (  5.00000 ,    5.19615) [rectangle    , fill=red    ] {};
\node at (  6.00000 ,    6.92820) [rectangle    , fill=red    ] {};
\node at (  5.50000 ,    4.33013) [rectangle    , fill=red    ] {};

\draw[very thick] (    .00000 ,     .00000) -- (  12.50000 ,     .00000);
\draw[very thick] (  12.50000 ,     .00000) -- (  12.50000 ,   10.39230);
\draw[very thick] (  12.50000 ,   10.39230) -- (    .00000 ,   10.39230);
\draw[very thick] (    .00000 ,   10.39230) -- (    .00000 ,     .00000);
\node at (   2.50000 ,     .86603) [circle    , fill=blue     ] {};
\node at (   2.50000 ,     .86603) [circle    , fill=blue     ] {};
\node at (   3.50000 ,     .86603) [circle    , fill=blue     ] {};
\draw [color=blue     ] (   2.50000 ,     .86603) -- (   3.50000 ,     .86603);
\node at (   4.50000 ,     .86603) [circle    , fill=blue     ] {};
\draw [color=blue     ] (   3.50000 ,     .86603) -- (   4.50000 ,     .86603);
\node at (   5.50000 ,     .86603) [circle    , fill=blue     ] {};
\draw [color=blue     ] (   4.50000 ,     .86603) -- (   5.50000 ,     .86603);
\node at (   6.00000 ,    1.73205) [circle    , fill=blue     ] {};
\draw [color=blue     ] (   5.50000 ,     .86603) -- (   6.00000 ,    1.73205);
\node at (   7.00000 ,    1.73205) [circle    , fill=blue     ] {};
\draw [color=blue     ] (   6.00000 ,    1.73205) -- (   7.00000 ,    1.73205);
\node at (   7.50000 ,    2.59808) [circle    , fill=blue     ] {};
\draw [color=blue     ] (   7.00000 ,    1.73205) -- (   7.50000 ,    2.59808);
\node at (   8.00000 ,    3.46410) [circle    , fill=blue     ] {};
\draw [color=blue     ] (   7.50000 ,    2.59808) -- (   8.00000 ,    3.46410);
\node at (   9.00000 ,    3.46410) [circle    , fill=blue     ] {};
\draw [color=blue     ] (   8.00000 ,    3.46410) -- (   9.00000 ,    3.46410);
\node at (   9.50000 ,    4.33013) [circle    , fill=blue     ] {};
\draw [color=blue     ] (   9.00000 ,    3.46410) -- (   9.50000 ,    4.33013);
\draw [color=blue     ] (   9.50000 ,    4.33013) -- (  10.00000 ,    5.19615);
\node at (  10.00000 ,    5.19615) [rectangle    , fill=green    ] {};
\end{tikzpicture}
\end{scaletikzpicturetowidth}
\caption{No knockouts allowed, $B=0, \Omega=2$, 10 time steps}
\label{fig3}
\end{figure}

Figure~\ref{fig3} shows the solution produced when no knockouts are allowed, so $B=0$. In this figure we can see that an agent reaches the target in ten time steps. The second agent shown in Figure~\ref{fig2} is not used at all, and this is automatically determined via
Equation~(\ref{eqexcess0}). Notice how in Figure~\ref{fig3}  the path followed by the agent never falls within a region covered by two or more sensors.

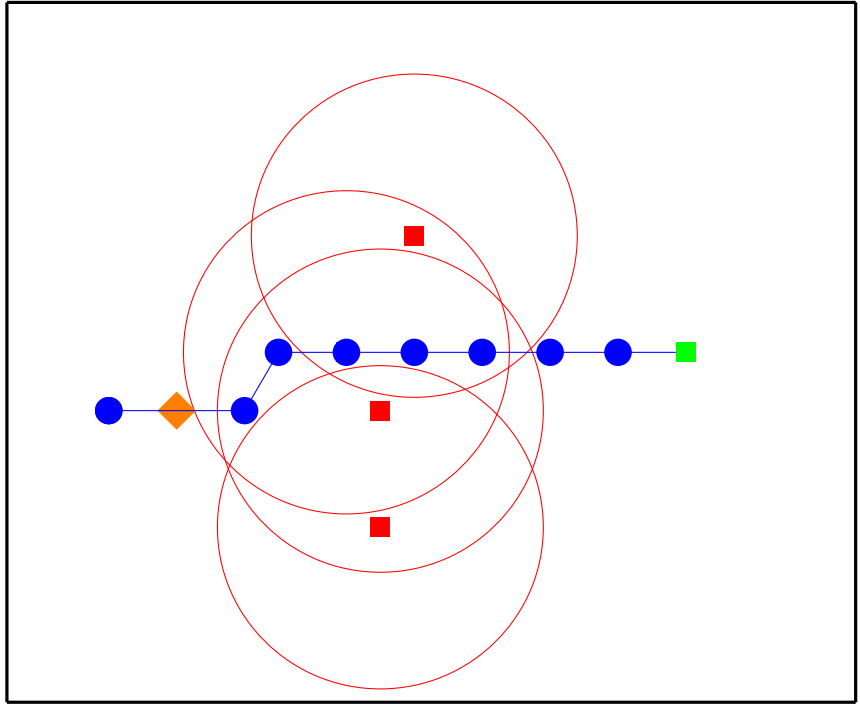
\begin{figure}[!htb]
\centering	  
\begin{scaletikzpicturetowidth}{\textwidth*0.70} 
\begin{tikzpicture}[scale=\tikzscale]	
\draw [red] (   5.50000 ,    2.59808) circle [radius=   2.40000];
\draw [red] (   5.00000 ,    5.19615) circle [radius=   2.40000];
\draw [red] (   6.00000 ,    6.92820) circle [radius=   2.40000];
\draw [red] (   5.50000 ,    4.33013) circle [radius=   2.40000];

\node at (  5.50000 ,    2.59808) [rectangle    , fill=red    ] {};
\node at (  5.00000 ,    5.19615) [rectangle    , fill=red    ] {};
\node at (  6.00000 ,    6.92820) [rectangle    , fill=red    ] {};
\node at (  5.50000 ,    4.33013) [rectangle    , fill=red    ] {};

\draw[very thick] (    .00000 ,     .00000) -- (  12.50000 ,     .00000);
\draw[very thick] (  12.50000 ,     .00000) -- (  12.50000 ,   10.39230);
\draw[very thick] (  12.50000 ,   10.39230) -- (    .00000 ,   10.39230);
\draw[very thick] (    .00000 ,   10.39230) -- (    .00000 ,     .00000);
\node at (   1.50000 ,    4.33013) [circle    , fill=blue     ] {};
\node at (   1.50000 ,    4.33013) [circle    , fill=blue     ] {};
\node at (   2.50000 ,    4.33013) [diamond    , fill=orange   ] {};
\draw [color=blue     ] (   1.50000 ,    4.33013) -- (   2.50000 ,    4.33013);
\node at (   3.50000 ,    4.33013) [circle    , fill=blue     ] {};
\draw [color=blue     ] (   2.50000 ,    4.33013) -- (   3.50000 ,    4.33013);
\node at (   4.00000 ,    5.19615) [circle    , fill=blue     ] {};
\draw [color=blue     ] (   3.50000 ,    4.33013) -- (   4.00000 ,    5.19615);
\node at (   5.00000 ,    5.19615) [circle    , fill=blue     ] {};
\draw [color=blue     ] (   4.00000 ,    5.19615) -- (   5.00000 ,    5.19615);
\node at (   6.00000 ,    5.19615) [circle    , fill=blue     ] {};
\draw [color=blue     ] (   5.00000 ,    5.19615) -- (   6.00000 ,    5.19615);
\node at (   7.00000 ,    5.19615) [circle    , fill=blue     ] {};
\draw [color=blue     ] (   6.00000 ,    5.19615) -- (   7.00000 ,    5.19615);
\node at (   8.00000 ,    5.19615) [circle    , fill=blue     ] {};
\draw [color=blue     ] (   7.00000 ,    5.19615) -- (   8.00000 ,    5.19615);
\node at (   9.00000 ,    5.19615) [circle    , fill=blue     ] {};
\draw [color=blue     ] (   8.00000 ,    5.19615) -- (   9.00000 ,    5.19615);
\draw [color=blue     ] (   9.00000 ,    5.19615) -- (  10.00000 ,    5.19615);
\node at (  10.00000 ,    5.19615) [rectangle    , fill=green    ] {};
\end{tikzpicture}
\end{scaletikzpicturetowidth}
\caption{One knockout allowed, $B=1, \Omega=2$, 9 time steps}
\label{fig4}
\end{figure}

Figure~\ref{fig4} shows the solution produced when one knockout is allowed, so $B=1$. In this figure we can see that an agent reaches the target in nine time steps. Note here that again only one agent is used, but a different agent from that shown in Figure~\ref{fig3}.

The solid orange diamond node shown after the first time step in 
Figure~\ref{fig4} indicates that at  the corresponding mesh vertex the agent initiates a knockout action and this action knockouts two of the sensors. This knockout leaves the agent free to proceed to the target without being detected by two, or more, sensors.

\begin{figure}[!htb]
\centering	  
\begin{scaletikzpicturetowidth}{\textwidth*0.70} 
\begin{tikzpicture}[scale=\tikzscale]	
\draw [red] (   5.50000 ,    2.59808) circle [radius=   2.40000];
\draw [red] (   6.00000 ,    6.92820) circle [radius=   2.40000];

\node at (  5.50000 ,    2.59808) [rectangle    , fill=red    ] {};
\node at (  6.00000 ,    6.92820) [rectangle    , fill=red    ] {};

\draw[very thick] (    .00000 ,     .00000) -- (  12.50000 ,     .00000);
\draw[very thick] (  12.50000 ,     .00000) -- (  12.50000 ,   10.39230);
\draw[very thick] (  12.50000 ,   10.39230) -- (    .00000 ,   10.39230);
\draw[very thick] (    .00000 ,   10.39230) -- (    .00000 ,     .00000);
\node at (   1.50000 ,    4.33013) [circle    , fill=blue     ] {};
\node at (   1.50000 ,    4.33013) [circle    , fill=blue     ] {};
\node at (   2.50000 ,    4.33013) [diamond    , fill=orange   ] {};
\draw [color=blue     ] (   1.50000 ,    4.33013) -- (   2.50000 ,    4.33013);
\node at (   3.50000 ,    4.33013) [circle    , fill=blue     ] {};
\draw [color=blue     ] (   2.50000 ,    4.33013) -- (   3.50000 ,    4.33013);
\node at (   4.00000 ,    5.19615) [circle    , fill=blue     ] {};
\draw [color=blue     ] (   3.50000 ,    4.33013) -- (   4.00000 ,    5.19615);
\node at (   5.00000 ,    5.19615) [circle    , fill=blue     ] {};
\draw [color=blue     ] (   4.00000 ,    5.19615) -- (   5.00000 ,    5.19615);
\node at (   6.00000 ,    5.19615) [circle    , fill=blue     ] {};
\draw [color=blue     ] (   5.00000 ,    5.19615) -- (   6.00000 ,    5.19615);
\node at (   7.00000 ,    5.19615) [circle    , fill=blue     ] {};
\draw [color=blue     ] (   6.00000 ,    5.19615) -- (   7.00000 ,    5.19615);
\node at (   8.00000 ,    5.19615) [circle    , fill=blue     ] {};
\draw [color=blue     ] (   7.00000 ,    5.19615) -- (   8.00000 ,    5.19615);
\node at (   9.00000 ,    5.19615) [circle    , fill=blue     ] {};
\draw [color=blue     ] (   8.00000 ,    5.19615) -- (   9.00000 ,    5.19615);
\draw [color=blue     ] (   9.00000 ,    5.19615) -- (  10.00000 ,    5.19615);
\node at (  10.00000 ,    5.19615) [rectangle   , fill=green    ] {};
\end{tikzpicture}
\end{scaletikzpicturetowidth}
\caption{Sensors remaining after knockout, $B=1, \Omega=2$}
\label{fig5}
\end{figure}
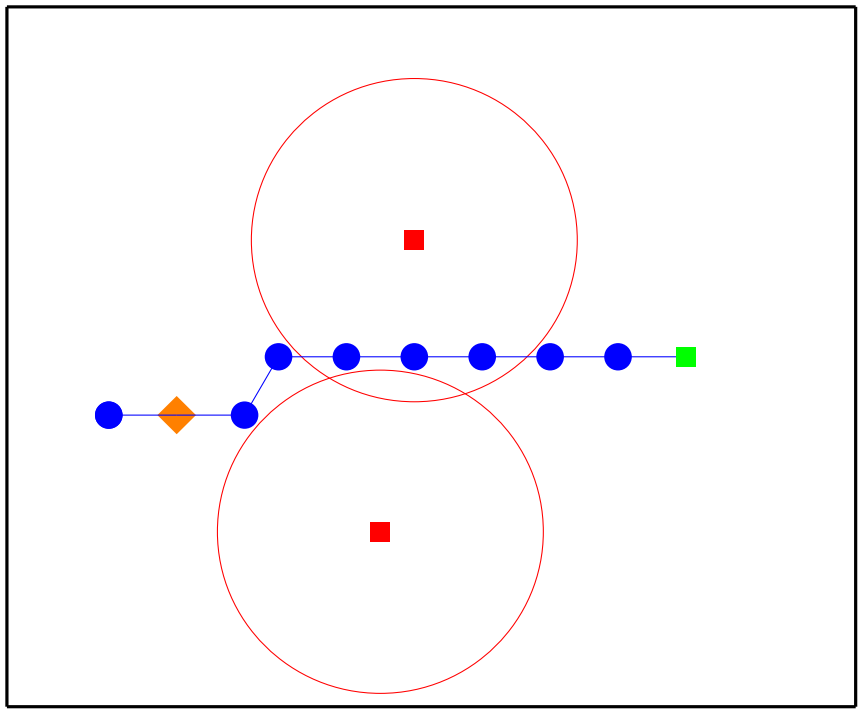

To make the situation clearly Figure~\ref{fig5} shows the same situation as Figure~\ref{fig4}, but with just the two sensors remaining after knockout shown.

\section{Detection probabilities} 
\label{sec4}	
The formulation above assumes detection by sensor $s$ at vertex $v$ is perfect if $d_{sv}=1$ and the sensor is not knocked out. Given the framework in terms of variables and constraints established above it is possible to include detection probabilities. 

Reinterpret $d_{sv}=1$ to mean that sensor $s$ can detect an agent at vertex $v$ under normal operating conditions, but with an associated probability $q_{sv}$ of \textbf{\emph{missing}} such a detection ($0 \leq  q_{sv} \leq 1$). 
This is equivalent to stating that sensor $s$ has a probability of $(1-q_{sv})$ of detecting an agent at vertex $v$. If $d_{sv}=0$, so vertex $v$ is outside the detection area of sensor $s$,  then $q_{sv}=1$.

Here we do not include the possibility of knockout, so all sensors are active, but with a detection probability that depends upon the sensor and each vertex in the mesh. Referring to Figure~\ref{fig2} any sensor will have zero probability of detecting an agent at any vertex $v$ outside its circular area, but a non-zero probability of detecting an agent within its circular area. As this probability is vertex dependent we can explicitly relate it to the distance between the sensor and the vertex.

Let $Q_v(\Omega)$ be the probability that an agent at vertex $v$ is not detected by the set of sensors $[s~|~s=1,\ldots,S]$. Note here that this probability depends upon $\Omega$, where an agent is detected if it is detected by $\Omega$ or more sensors,  it evades detection if it is detected by $(\Omega-1)$ or less sensors, 
c.f.~Equation~(\ref{deq3}).

For example suppose a vertex $v$ is covered by all sensors, so $d_{sv}=1~s=1,\ldots,S$. Suppose too that one constraint is that an agent must not be detected by two or more sensors, so $\Omega=2$. So in this case  $Q_v(2)$ represents the probability that an agent at vertex $v$ is only detected by at most one sensor and we have: 
\begin{equation}
Q_v(2) = \prod_{s=1}^S q_{sv} + \sum_{s=1}^S (1-q_{sv})\prod_{r=1,~ r \neq s}^S q_{rv}
\label{qdef}
\end{equation}
In Equation~(\ref{qdef}) the first term is the probability that all sensors detect nothing at vertex $v$, and the second term is the probability that just one sensor detects an agent at vertex $v$.
Note here that $Q_v(\Omega)$ is simply an input to our formulation and could be computed by any user-defined function if we so  wish. 

Then the  \textbf{\emph{probability of evading detection, henceforth PED,}} namely the probability that, over all time and over all agents,  there is no detection of any agent is given by:
\begin{equation}
PED = \prod_{t=1}^T  ~\prod_{a=1}^A~\prod_{v=1~x_{avt}=1}^{N-1}  Q_{v}(\Omega)  
\label{prob1n1}
\end{equation}

In Equation~(\ref{prob1n1}) we assume that all the probabilities are independent and the expression seen is the total probability of no detections over the entire time horizon $T$. 
Clearly Equation~(\ref{prob1n1}) is nonlinear but can be linearised to:
 \begin{equation}
\sum_{t=1}^T  ~\sum_{a=1}^A ~\sum_{v=1}^{N-1} ~ \ln (Q_v(\Omega)) x_{avt}
\label{prob3n}
\end{equation}
Notice how in Equation~(\ref{prob3n}) we have included the  $x_{avt}$ terms in the linear sum of  logarithmic terms. This equation represents PED in logarithmic form.

Note  here that when we are dealing with probabilities (and no knockouts) the variables $\alpha_{avt}$, $y_{sat}$ and $\lambda_{st}$ are irrelevant and can be deleted from consideration.

We have excluded  vertex $N$ from Equations~(\ref{prob1n1}) and (\ref{prob3n})   since we have assumed here for simplicity that reaching that vertex renders detection irrelevant.  For example this may be because the target is a command and control facility to which all sensors report and its destruction by an agent means the sensors  are effectively inactive, as they now lack a functioning facility  to which to report.    
Alternatively the sensors may remain active after the target has been reached and the agent  then needs to make their escape, whilst avoiding detection if possible. In this scenario the agent moves to the agent extraction (exfil)  location, which would have been chosen to be outside sensor detection. This can be dealt with as discussed above under agent exit strategy, with vertex $N$ now being the extraction location. However it is important to note here that if we wish to relax the assumption that   detection at vertex $N$ is irrelevant then this can be done with only minor changes to the formulation.

To incorporate detection probabilities into the problem a natural approach would be to maximise Equation~(\ref{prob3n}) for a given value of $T$, the time horizon. Varying $T$ will show how PED changes as we vary the time horizon. Note here that an alternative approach to varying $T$ and maximising PED would be to set  a lower limit on PED and then find the minimal time to target 
which will achieve the required PED.

As an illustration of what can be done in terms of detection probabilities we took the example problem considered above and set the detection probability for sensor $s$ and vertex $v$ (in sensor detection range) as equal to 1/[1+(distance from $s$ to $v$)$^2$/(sensor detection radius)$^2$].  
Detection of an agent at a vertex only occurs if two or more sensors detect the agent at the vertex, so $\Omega=2$. 

In this situation an agent can reach the target in ten time steps without being detected, so PED=1. This can be seen in 
Figure~\ref{fig7} and note how the path only involves sensor detection by a single sensor.

For a path involving one less time step, so nine time steps, the path with the maximum PED by two or more sensors is shown in Figure~\ref{fig8}. This path has a PED value much less than 1\%. Although clearly the precise probability will depend upon the situation under consideration, the point here is that this example exposes the decision-maker to an explicit trade-off. Do they go for a path involving nine time steps but a very low PED, or do they take a longer path involving ten time steps but with a high PED value?

\begin{figure}[!htb]
\centering	  
\begin{scaletikzpicturetowidth}{\textwidth*0.70} 
\begin{tikzpicture}[scale=\tikzscale]	
\draw [red] (   5.50000 ,    2.59808) circle [radius=   2.40000];
\draw [red] (   5.00000 ,    5.19615) circle [radius=   2.40000];
\draw [red] (   6.00000 ,    6.92820) circle [radius=   2.40000];
\draw [red] (   5.50000 ,    4.33013) circle [radius=   2.40000];
\node at (  5.50000 ,    2.59808) [rectangle    , fill=red    ] {};
\node at (  5.00000 ,    5.19615) [rectangle    , fill=red    ] {};
\node at (  6.00000 ,    6.92820) [rectangle    , fill=red    ] {};
\node at (  5.50000 ,    4.33013) [rectangle    , fill=red    ] {};

\draw[very thick] (    .00000 ,     .00000) -- (  12.50000 ,     .00000);
\draw[very thick] (  12.50000 ,     .00000) -- (  12.50000 ,   10.39230);
\draw[very thick] (  12.50000 ,   10.39230) -- (    .00000 ,   10.39230);
\draw[very thick] (    .00000 ,   10.39230) -- (    .00000 ,     .00000);
\node at (   2.50000 ,     .86603) [circle    , fill=blue     ] {};
\node at (   2.50000 ,     .86603) [circle    , fill=blue     ] {};
\node at (   3.50000 ,     .86603) [circle    , fill=blue     ] {};
\draw [color=blue     ] (   2.50000 ,     .86603) -- (   3.50000 ,     .86603);
\node at (   4.00000 ,    1.73205) [circle    , fill=blue     ] {};
\draw [color=blue     ] (   3.50000 ,     .86603) -- (   4.00000 ,    1.73205);
\node at (   5.00000 ,    1.73205) [circle    , fill=blue     ] {};
\draw [color=blue     ] (   4.00000 ,    1.73205) -- (   5.00000 ,    1.73205);
\node at (   6.00000 ,    1.73205) [circle    , fill=blue     ] {};
\draw [color=blue     ] (   5.00000 ,    1.73205) -- (   6.00000 ,    1.73205);
\node at (   7.00000 ,    1.73205) [circle    , fill=blue     ] {};
\draw [color=blue     ] (   6.00000 ,    1.73205) -- (   7.00000 ,    1.73205);
\node at (   7.50000 ,    2.59808) [circle    , fill=blue     ] {};
\draw [color=blue     ] (   7.00000 ,    1.73205) -- (   7.50000 ,    2.59808);
\node at (   8.50000 ,    2.59808) [circle    , fill=blue     ] {};
\draw [color=blue     ] (   7.50000 ,    2.59808) -- (   8.50000 ,    2.59808);
\node at (   9.00000 ,    3.46410) [circle    , fill=blue     ] {};
\draw [color=blue     ] (   8.50000 ,    2.59808) -- (   9.00000 ,    3.46410);
\node at (   9.50000 ,    4.33013) [circle    , fill=blue     ] {};
\draw [color=blue     ] (   9.00000 ,    3.46410) -- (   9.50000 ,    4.33013);
\draw [color=blue     ] (   9.50000 ,    4.33013) -- (  10.00000 ,    5.19615);
\node at (  10.00000 ,    5.19615) [rectangle    , fill=green    ] {};
\end{tikzpicture}
\end{scaletikzpicturetowidth}
\caption{Path with  10 time steps and  maximum probability of evading detection by two or more sensors, PED=1}
\label{fig7}
\end{figure}
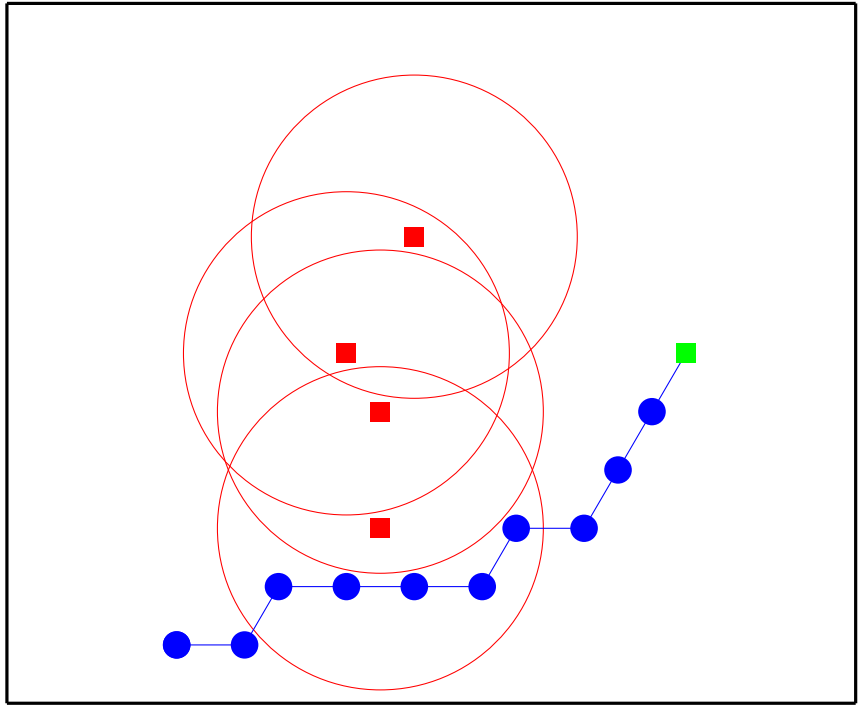

\begin{figure}[!htb]
\centering	  
\begin{scaletikzpicturetowidth}{\textwidth*0.70} 
\begin{tikzpicture}[scale=\tikzscale]	
\draw [red] (   5.50000 ,    2.59808) circle [radius=   2.40000];
\draw [red] (   5.00000 ,    5.19615) circle [radius=   2.40000];
\draw [red] (   6.00000 ,    6.92820) circle [radius=   2.40000];
\draw [red] (   5.50000 ,    4.33013) circle [radius=   2.40000];

\node at (  5.50000 ,    2.59808) [rectangle    , fill=red    ] {};
\node at (  5.00000 ,    5.19615) [rectangle    , fill=red    ] {};
\node at (  6.00000 ,    6.92820) [rectangle    , fill=red    ] {};
\node at (  5.50000 ,    4.33013) [rectangle    , fill=red    ] {};

\draw[very thick] (    .00000 ,     .00000) -- (  12.50000 ,     .00000);
\draw[very thick] (  12.50000 ,     .00000) -- (  12.50000 ,   10.39230);
\draw[very thick] (  12.50000 ,   10.39230) -- (    .00000 ,   10.39230);
\draw[very thick] (    .00000 ,   10.39230) -- (    .00000 ,     .00000);
\node at (   1.50000 ,    4.33013) [circle    , fill=blue     ] {};
\node at (   1.50000 ,    4.33013) [circle    , fill=blue     ] {};
\node at (   2.50000 ,    4.33013) [circle    , fill=blue     ] {};
\draw [color=blue     ] (   1.50000 ,    4.33013) -- (   2.50000 ,    4.33013);
\node at (   3.00000 ,    5.19615) [circle    , fill=blue     ] {};
\draw [color=blue     ] (   2.50000 ,    4.33013) -- (   3.00000 ,    5.19615);
\node at (   4.00000 ,    5.19615) [circle    , fill=blue     ] {};
\draw [color=blue     ] (   3.00000 ,    5.19615) -- (   4.00000 ,    5.19615);
\node at (   5.00000 ,    5.19615) [circle    , fill=blue     ] {};
\draw [color=blue     ] (   4.00000 ,    5.19615) -- (   5.00000 ,    5.19615);
\node at (   6.00000 ,    5.19615) [circle    , fill=blue     ] {};
\draw [color=blue     ] (   5.00000 ,    5.19615) -- (   6.00000 ,    5.19615);
\node at (   7.00000 ,    5.19615) [circle    , fill=blue     ] {};
\draw [color=blue     ] (   6.00000 ,    5.19615) -- (   7.00000 ,    5.19615);
\node at (   8.00000 ,    5.19615) [circle    , fill=blue     ] {};
\draw [color=blue     ] (   7.00000 ,    5.19615) -- (   8.00000 ,    5.19615);
\node at (   9.00000 ,    5.19615) [circle    , fill=blue     ] {};
\draw [color=blue     ] (   8.00000 ,    5.19615) -- (   9.00000 ,    5.19615);
\draw [color=blue     ] (   9.00000 ,    5.19615) -- (  10.00000 ,    5.19615);
\node at (  10.00000 ,    5.19615) [rectangle    , fill=green    ] {};
\end{tikzpicture}
\end{scaletikzpicturetowidth}
\caption{Path with 9 time steps and maximum probability of evading detection by two or more sensors, PED $<$ 0.01}
\label{fig8}
\end{figure}
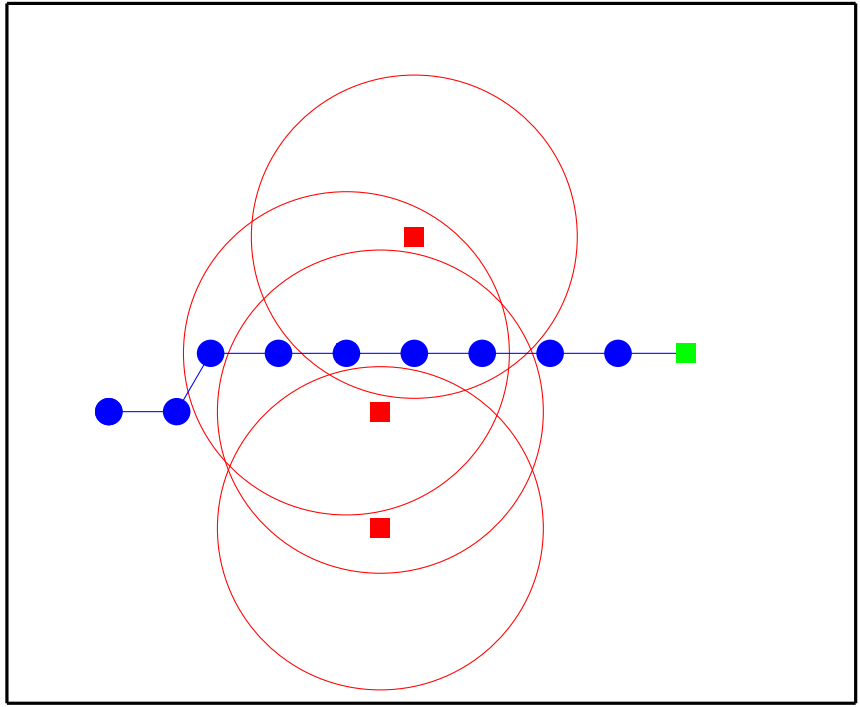

\subsection{Combining detection probabilities with knockout}

It does not (in general) appear possible to directly incorporate knockout into the mathematical formulation given above when we have detection probabilities present. The fundamental reason for this is that to calculate the value of the $Q_v(\Omega)$ term, which is involved in linearised Equation~(\ref{prob3n}) via a logarithmic expression, we need to know which sensors are active.

For example consider the calculation given above for $Q_v(2)$, namely
$Q_v(2) = \prod_{s=1}^S q_{sv} + \sum_{s=1}^S (1-q_{sv})\prod_{r=1,~ r \neq s}^S q_{rv}$. Incorporating into this expression variables indicating whether a sensor is active or knocked out, and then computing the logarithmic term, will in general make the formulation nonlinear.

However given that we have the formulation above that will find the path with the maximum PED for a given  set of active sensors then it is clear that an iterative scheme could be adopted to incorporate knockout with detection probabilities. Such a scheme might be as follows.
Firstly, given a desired time horizon to reach the target find a path with the maximum PED assuming all sensors are active. 
Secondly, examine the path and choose an active sensor which is influencing the path to target to knockout.
Thirdly, iteratively remove or retain sensors, until you have a path that you are satisfied with.
An iterative scheme such as this could incorporate human insight, or be purely algorithmic in nature. Developing such a scheme however is beyond the scope of the work reported here.

\section{Sensor confusion} 
\label{sec5}	
We interpret sensor confusion, also known as sensor degradation,
to mean that the previous probability $q_{sv}$ of a sensor missing a detection is increased. Let $q_{sv}^c$ be the probability that applies after confusion.  
Without significant loss of generality we assume for simplicity that $q_{sv}^c \geq q_{sv}$, so confusion increases the chances that a sensor misses detecting an agent at a vertex.

We assume here that confusion, if actioned by an agent, applies to all sensors equally. In other words confusion is not range restricted, unlike knockout. This assumption was made to simplify the  problem.  The underlying scenario here is that  it  may be possible to cause confusion in a connected sensor network by interacting with only one sensor in that connected network. This interaction may be by direct agent action  or by remote action.

Let
$\beta_{at}=1$ if a confusion action is initiated (started) by agent $a$ at time $t$, $= 0$ otherwise.
Let $z_t=1$ if all sensors are at confused at time $t$, $= 0$ otherwise.
Here the $\beta_{at}$ variables relate to when sensors are first confused, so initiation of a confusion action by an agent. Confusion initiated by agent $a$ applies for $\Delta_a^c$ time steps.
The $z_t$ variables relate to whether sensors are confused or not at time $t$, for example by a previous confusion initiation.

We have the following constraint:
\begin{equation}
z_t = \sum_{a=1}^A \sum_{\tau=\mbox{max}[1,t-\Delta_a^c]}^t \beta_{a\tau}
~~~~t=1,\ldots,T
\label{c1}
\end{equation}
Equation~(\ref{c1}) ensures that if at time $t$ a confusion action has been initiated by any agent $a$ in the appropriate past time period $[\mbox{max}[1,t-\Delta_a^c],t]$ then $z_t$ will be one, zero otherwise. 

We assume here that the nature of the problem considered is such that there is no additional benefit in confusing sensors that are already confused for a second or third, or more, time.  Consequently the right-hand side of Equation~(\ref{c1}) is always zero or one, as the constraint implies with a zero-one variable on the left-hand side of the constraint.

An agent cannot initiate a confusion action if they are no longer active, i.e.~if they have made a transition to vertex $0$, so we have:
\begin{equation}
\beta_{at} \leq 1- x_{a0t} ~~~~a=1,\ldots,A;~t=1,\ldots,T
\label{c1a}
\end{equation}

Obviously the optimal strategy is to always confuse sensors unless there is some limitation and so here we assume that agent $a$ initiating a confusion action involves a cost $C_a^c$ and so the constraint limiting the total cost, Equation~(\ref{deq1}), becomes:
\begin{equation}
\sum_{a=1}^A \sum_{v=1}^N \sum_{t=1}^T C_a \alpha_{avt} + 
\sum_{a=1}^A \sum_{t=1}^T C_a^c \beta_{at} \leq B 
\label{cost1}
\end{equation}

Note here that with agent confusion actions added to the problem the previous constraint, Equation~(\ref{eqexcess1}), dealing with excess agents changes to:
\begin{equation}
x_{a01} \geq 1 - \sum_{t=1}^T x_{aNt} - \sum_{v=1}^N \sum_{t=1}^T
\alpha_{avt} - \sum_{t=1}^T\beta_{at}~~~~ a=1,\ldots,A
\label{eqexcess2}
\end{equation}
Equation~(\ref{eqexcess2}) ensures that $ x_{a01}$ is forced to be one if agent $a$ never reaches the target vertex $N$, and never engages in a knockout or confusion action.

Let $Q_v^c(\Omega)$ be the probability that an agent at vertex $v$ is not detected by the set of sensors $[s~|~s=1,\ldots,S]$ when the sensors are confused. Here $Q_v^c(\Omega)$ is defined in exactly the same way as 
$Q_v(\Omega)$ above, but using $q_{sv}^c$ rather than $q_{sv}$.
Without significant loss of generality we assume $Q_v^c(\Omega) \geq  Q_v(\Omega)$, so confusion improves the chances of evading sensor detection. 

Equation~(\ref{prob1n1}), the probability that, over all time and over all agents,  there is no detection of any agent now becomes:
\begin{equation}
PED = \Bigl \{ \prod_{t=1}^T  ~\prod_{a=1}^A~\prod_{v=1~x_{avt}=1~z_t=0}^{N-1} Q_{v}(\Omega)\Bigr \} ~ 
\Bigl \{ \prod_{t=1}^T  ~\prod_{a=1}^A~\prod_{v=1~x_{avt}=1~z_t=1}^{N-1} Q_v^c(\Omega) \Bigr \}
\label{c2}
\end{equation}

Equation~(\ref{c2}) contains two product terms, one for time periods when no confusion is present (so $z_t=0$), and one for time periods when confusion is present (so $z_t=1$). Now Equation~(\ref{c2}) can be rewritten as:
\begin{equation}
\Bigl \{ \prod_{t=1}^T  ~\prod_{a=1}^A~\prod_{v=1~x_{avt}=1}^{N-1} Q_{v}(\Omega)\Bigr \} ~ 
\Bigl \{ \prod_{t=1}^T  ~\prod_{a=1}^A~\prod_{v=1~x_{avt}=1~z_t=1}^{N-1} (Q_v^c(\Omega)/Q_v(\Omega)) \Bigr \}
\label{c3}
\end{equation}
Here the first product term is now independent of the value of $z_t$, but we have amended the second product term to account for this.

Now in order to linearise Equation~(\ref{c3}) we need introduce variables $\delta_{avt}$ where $\delta_{avt}=1$ if $x_{avt}=1$ and $z_t=1$, zero otherwise. To represent this we need the constraints:
\begin{equation}
x_{avt} + z_t - 1 \leq \delta_{avt} \leq x_{avt}~~~~ a=1,\ldots,A; ~v=1,\ldots,N-1;~t=1,\ldots,T
\label{c4}
\end{equation}
\begin{equation}
\delta_{avt} \leq z_t~~~~ a=1,\ldots,A; ~v=1,\ldots,N-1;~t=1,\ldots,T
\label{c4a}
\end{equation}
Equations~(\ref{c4}) and (\ref{c4a}) ensure that $\delta_{avt}=1$ if and only if $x_{avt}=1$ and $z_t=1$. Equation~(\ref{c3}) can now be rewritten as:
\begin{equation}
\Bigl \{ \prod_{t=1}^T  ~\prod_{a=1}^A~\prod_{v=1~x_{avt}=1}^{N-1} Q_{v}(\Omega)\Bigr \} ~ 
\Bigl \{ \prod_{t=1}^T  ~\prod_{a=1}^A~\prod_{v=1~\delta_{avt}=1}^{N-1} (Q_v^c(\Omega)/Q_v(\Omega)) \Bigr \}
\label{c5}
\end{equation}
So we can now linearise Equation~(\ref{c5}) so that it becomes:
\begin{equation}
\Bigl \{ \sum_{t=1}^T  ~\sum_{a=1}^A~\sum_{v=1}^{N-1} \ln (Q_{v}(\Omega)) x_{avt} \Bigr \} ~+~ 
\Bigl \{ \sum_{t=1}^T  ~\sum_{a=1}^A~\sum_{v=1}^{N-1} \ln (Q_v^c(\Omega)/Q_v(\Omega))  \delta_{avt} \Bigr \}
\label{c6}
\end{equation}

Hence if sensor confusion is present we can simply proceed as before, maximising the probability of evading detection, PED, as represented, in logarithmic form, by Equation~(\ref{c6}).

\begin{figure}[!htb]
\centering	  
\begin{scaletikzpicturetowidth}{\textwidth*0.70} 
\begin{tikzpicture}[scale=\tikzscale]	
\draw [red] (   5.50000 ,    2.59808) circle [radius=   2.40000];
\draw [red] (   5.00000 ,    5.19615) circle [radius=   2.40000];
\draw [red] (   6.00000 ,    6.92820) circle [radius=   2.40000];
\draw [red] (   5.50000 ,    4.33013) circle [radius=   2.40000];

\node at (  5.50000 ,    2.59808) [rectangle    , fill=red    ] {};
\node at (  5.00000 ,    5.19615) [rectangle    , fill=red    ] {};
\node at (  6.00000 ,    6.92820) [rectangle    , fill=red    ] {};
\node at (  5.50000 ,    4.33013) [rectangle    , fill=red    ] {};

\draw[very thick] (    .00000 ,     .00000) -- (  12.50000 ,     .00000);
\draw[very thick] (  12.50000 ,     .00000) -- (  12.50000 ,   10.39230);
\draw[very thick] (  12.50000 ,   10.39230) -- (    .00000 ,   10.39230);
\draw[very thick] (    .00000 ,   10.39230) -- (    .00000 ,     .00000);
\node at (   1.50000 ,    4.33013) [circle    , fill=blue     ] {};
\node at (   1.50000 ,    4.33013) [circle    , fill=blue     ] {};
\node at (   2.50000 ,    4.33013) [diamond    , fill=orange   ] {};
\draw [color=blue     ] (   1.50000 ,    4.33013) -- (   2.50000 ,    4.33013);
\node at (   3.00000 ,    5.19615) [circle    , fill=blue     ] {};
\draw [color=blue     ] (   2.50000 ,    4.33013) -- (   3.00000 ,    5.19615);
\node at (   4.00000 ,    5.19615) [circle    , fill=blue     ] {};
\draw [color=blue     ] (   3.00000 ,    5.19615) -- (   4.00000 ,    5.19615);
\node at (   5.00000 ,    5.19615) [circle    , fill=blue     ] {};
\draw [color=blue     ] (   4.00000 ,    5.19615) -- (   5.00000 ,    5.19615);
\node at (   6.00000 ,    5.19615) [circle    , fill=blue     ] {};
\draw [color=blue     ] (   5.00000 ,    5.19615) -- (   6.00000 ,    5.19615);
\node at (   7.00000 ,    5.19615) [circle    , fill=blue     ] {};
\draw [color=blue     ] (   6.00000 ,    5.19615) -- (   7.00000 ,    5.19615);
\node at (   8.00000 ,    5.19615) [circle    , fill=blue     ] {};
\draw [color=blue     ] (   7.00000 ,    5.19615) -- (   8.00000 ,    5.19615);
\node at (   9.00000 ,    5.19615) [circle    , fill=blue     ] {};
\draw [color=blue     ] (   8.00000 ,    5.19615) -- (   9.00000 ,    5.19615);
\draw [color=blue     ] (   9.00000 ,    5.19615) -- (  10.00000 ,    5.19615);
\node at (  10.00000 ,    5.19615) [rectangle    , fill=green    ] {};
\end{tikzpicture}
\end{scaletikzpicturetowidth}
\caption{Confusion solution, 9 time steps, PED=0.95}
\label{fig10}
\end{figure}
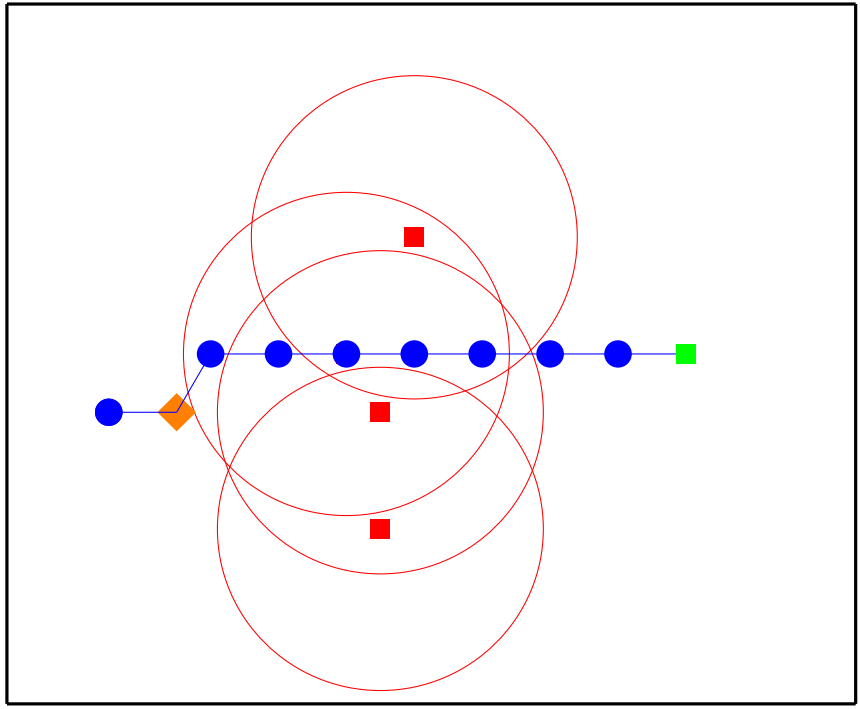

As an illustration of confusion we took the detection probability  example, with nine time steps, considered previously above in Figure~\ref{fig8} that had an associated PED value of less than 1\%. On the assumption that confusion reduces detection probability to 10\% of its previous value the solution for maximising PED is as shown in 
Figure~\ref{fig10}. In that figure a confusion action is initiated at time step 1. The probability of evading detection associated with Figure~\ref{fig10} is approximately 95\%.

\subsection{Required PED}

Note here that, as mentioned previously above, an alternative approach to varying $T$ and maximising detection probability would be to set  a lower limit on the detection probability and then find the minimal time to target 
which will achieve that detection probability. If we consider this approach together with confusion, then if the required probability of evading detection  is $Q^*$, it follows from  Equation~(\ref{c2}) that the constraint upon detection probability is:
\begin{equation}\Bigl \{ \prod_{t=1}^T  ~\prod_{a=1}^A~\prod_{v=1~x_{avt}=1~z_t=0}^{N-1} Q_{v}(\Omega)\Bigr \} ~ 
\Bigl \{ \prod_{t=1}^T  ~\prod_{a=1}^A~\prod_{v=1~x_{avt}=1~z_t=1}^{N-1} Q_v^c(\Omega) \Bigr \} \geq Q^*
\label{c2a}
\end{equation}
This can obviously be linearised in the same manner as done previously  above to become:
\begin{equation}
\Bigl \{ \sum_{t=1}^T  ~\sum_{a=1}^A~\sum_{v=1}^{N-1} \ln (Q_{v}(\Omega)) x_{avt} \Bigr \} ~+~ 
\Bigl \{ \sum_{t=1}^T  ~\sum_{a=1}^A~\sum_{v=1}^{N-1} \ln (Q_v^c(\Omega)/Q_v(\Omega))  \delta_{avt} \Bigr \} \geq \ln(Q^*)
\label{c6lim}
\end{equation}

Note here however that the left-hand side of  Equation~(\ref{c2a}) is the product of probabilities and we require this product to be at least  $Q^*$. It immediately follows that any term active in the product on the left-hand side of  Equation~(\ref{c2a}) must be at least $Q^*$. This leads to the constraints:
\begin{equation}
x_{avt} \leq z_t ~~~~ Q_v(\Omega) < Q^*; ~a=1,\ldots,A;~v=1\ldots,N-1;~t=1,\ldots,T 
\label{reqc2}
\end{equation}
\begin{equation}
x_{avt} = 0 ~~~~  Q_v^c(\Omega) < Q^*; ~a=1,\ldots,A;~v=1\ldots,N-1;~t=1,\ldots,T
\label{reqc1}
\end{equation}
Equation~(\ref{reqc2}) ensures that no agent $a$ can be at vertex $v$ at time $t$ if it would lead to violation of Equation~(\ref{c2a}) when confusion is not active at 
time $t$ (so $z_t=0$). Equation~(\ref{reqc1}) ensures that no agent $a$ can be at vertex $v$ at time $t$ if, even with confusion active at time $t$ (so $z_t=1$), it would lead to violation of Equation~(\ref{c2a}). 

\subsection{Biobjective shortest path problem}

In the literature  the  biobjective shortest path problem, e.g.~\cite{Duque2015, Ghoseiri2010,  Raith2009,  Sedeno2015,  Sedeno2019}, typically involves finding the shortest distance path between two vertices on a graph, but where the path has a secondary characteristic that has to be considered, typically  minimised.  This secondary characteristic may for example be cost, and we need to find a set of paths that enable us to trade off the distance travelled with the cost incurred. 

This path set ideally needs to be an efficient set, so that for any path in the set with the two characteristic values being $(c_1,c_2)$ it is not possible to improve upon one of these two values without deteriorating, increasing, the other value.

It is clear that the biobjective shortest path problem has similarities with the problem considered here of minimising the time to target, whilst also considering the associated PED value, in the special case of just one agent. 

There is a key difference between the standard  biobjective shortest path problem as considered in the literature however, and the PED problem considered here. 
This difference  is that in the literature once a path is fixed the values of the two characteristics associated with that path are automatically known, for example by summing the characteristic values for each arc on the path. 
By contrast in the problem considered here once a path is fixed we know the time to target, \textbf{\emph{but}} the PED value associated with that path is not known. Rather it is only known once a decision problem relating to when ($\beta_{at}$) to initiate confusion actions 
has been solved.

For this reason it is not clear how approaches presented in the literature associated with the  biobjective shortest path problem could be adapted to deal with the problem of simultaneously considering time to target and PED, even in the case of just a single agent. Hence this is left as a topic for future research.

\subsection{Successive shortest path  heuristic}
We  implemented a heuristic for the problem considered in this paper,  based upon generating successive shortest paths. This heuristic was as follows.
Consider each agent in turn. Find the shortest path from the agent to the target. Enumerate all possible combinations of knockout/confusion actions for the vertices on the path. If a combination is found for an agent that satisfies the requirements of the problem then terminate consideration of that agent. Otherwise delete from the mesh the vertex on the path that is  within detection range of the maximum number of sensors. Then repeat the procedure, finding the shortest path from the agent to the target, etc, as above. 
Once all agents have been considered choose the best path found over all agents. 

\section{Agent restrictions} 
\label{sec6}	

Currently agent actions such as knockout or confusion incur a cost, see Equation~(\ref{cost1}). It is possible to extend the formulation given above to include the situation where an agent action can only be carried out once within a specified number of time steps. It is also possible to extend the formulation given above to include the situation where
an agent action imposes a movement/time penalty, such that the agent has to remain at the vertex from which the action was initiated for a fixed number of time steps.
We deal with both of these  below.

\subsection{Time restrictions}
Let $\Phi_a$ ($1 \leq \Phi_a \leq T-1$) be for
for agent $a$ the specified number of time steps within which we can only carry out one knockout action. Then the appropriate constraint is:
\begin{equation}
\sum_{v=1}^N \sum_{\tau=t}^{t+\Phi_a}
\alpha_{av\tau}    \leq 1  
~~~~t=1,\ldots,(T-\Phi_a);~a=1,\ldots,A
\label{eqkt}
\end{equation}
Equation~(\ref{eqkt}) ensures that at most one $\alpha_{avt}$ can be one over all vertices within time periods of length $\Phi_a$.

Similarly if 
$\Phi_a^c$ ($1 \leq \Phi_a^c \leq T-1$) 
is for agent $a$ the specified number of time steps within which we can only carry out one confusion action then the appropriate constraint is:
\begin{equation}
\sum_{\tau=t}^{t+\Phi_a^c}
\beta_{a\tau}    \leq 1  
~~~~t=1,\ldots,(T-\Phi_a^c);~a=1,\ldots,A
\label{eqktc}
\end{equation}

\subsection{Movement restrictions}
Here we deal with the situation where an agent action imposes a movement/time  penalty,  
such that the agent has to remain at the vertex from which the action was initiated for a fixed number of time steps.

Let  $\Psi_a $ ($1 \leq \Psi_a \leq T-1$) be for agent $a$ the specified  number of time steps for which the agent must remain at a vertex from which it  initiated  
a knockout action. Then the appropriate constraint is:
\begin{equation}
x_{a0\tau} + x_{av\tau} \geq \alpha_{avt}~~~~\tau=t,\ldots,\mbox{min}[t+\Psi_a,T];
~t=1,\ldots,T;~v=1,\ldots,N;~a=1,\ldots,A 
\label{eqkstat}
\end{equation}
Equation~(\ref{eqkstat}) ensures that  if  $\alpha_{avt}$ is one then 
$x_{av\tau}$ is one for an appropriate  number of time periods $\tau \geq t$, or the agent moves to vertex $0$.

Similarly suppose
$\Psi_a^c $ ($1 \leq \Psi_a^c \leq T-1$)
is for agent $a$ the specified  number of time steps for which the agent must remain at a vertex from which it  initiated  
a confusion action.  Then the appropriate constraint is:
\begin{equation}
x_{a0\tau} + x_{av\tau} \geq x_{avt} - (1 - \beta_{at})~~~~\tau=t,\ldots,\mbox{min}[t+\Psi_a^c,T];
~t=1,\ldots,T;~v=1,\ldots,N;~a=1,\ldots,A 
\label{eqkstatc}
\end{equation}
In Equation~(\ref{eqkstatc}) if $\beta_{at}=0$ the constraint is inactive, since in that case the right-hand side is either zero or minus one. If $\beta_{at}=1$, 
i.e.~agent $a$ initiates a confusion action at time $t$ then there is only one vertex $v$ for which $x_{avt}$ is non-zero, recall here that each agent is only at one vertex at each time step, and so the constraint will ensure that the agent remains at that vertex for the appropriate number of time steps, or the agent moves to vertex $0$.

\section{Computational results} 
\label{sec7}	

Clearly whether, or not, the formulation given here is of practical value will depend on how it performs computationally. In other words can we use the formulation to solve problems of practical size in a reasonable computational time.

\begin{figure}[!htb]
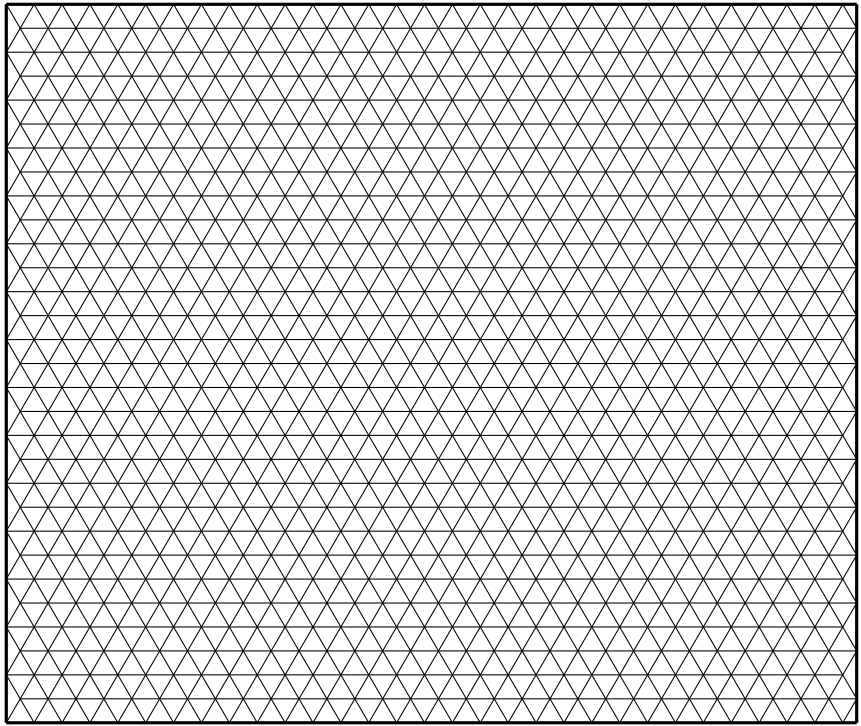

\centering	  
\begin{scaletikzpicturetowidth}{\textwidth*0.70} 

\end{scaletikzpicturetowidth}
\caption{Triangular mesh: $30 \times 30$}
\label{fig30}
\end{figure}

\begin{figure}[!htb]
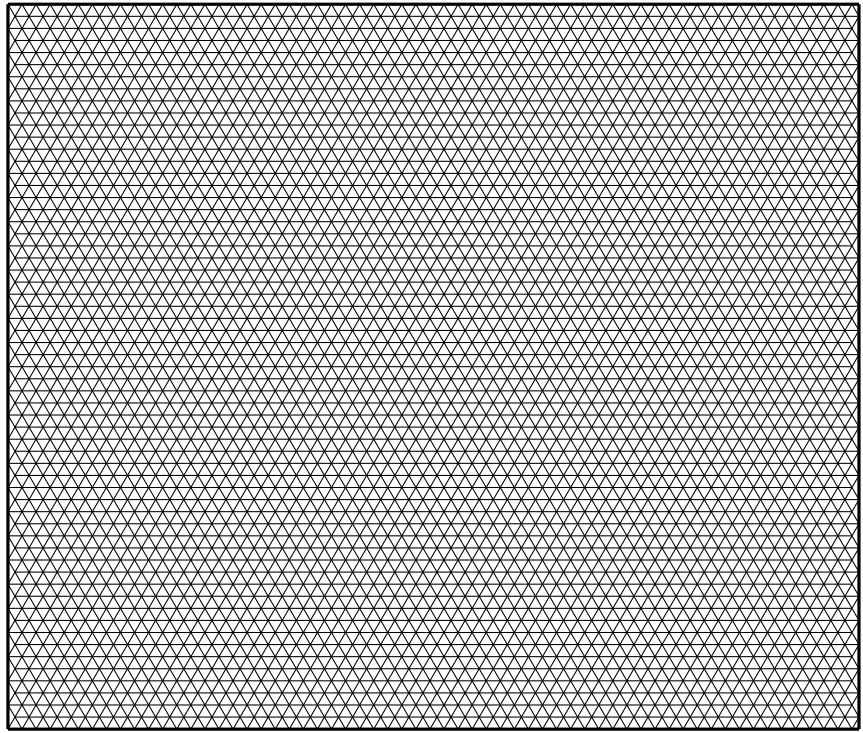

\centering	  
\begin{scaletikzpicturetowidth}{\textwidth*0.70} 
%
\end{scaletikzpicturetowidth}
\caption{Triangular mesh: $60 \times 60$}
\label{fig60}
\end{figure}

\begin{figure}[!htb]
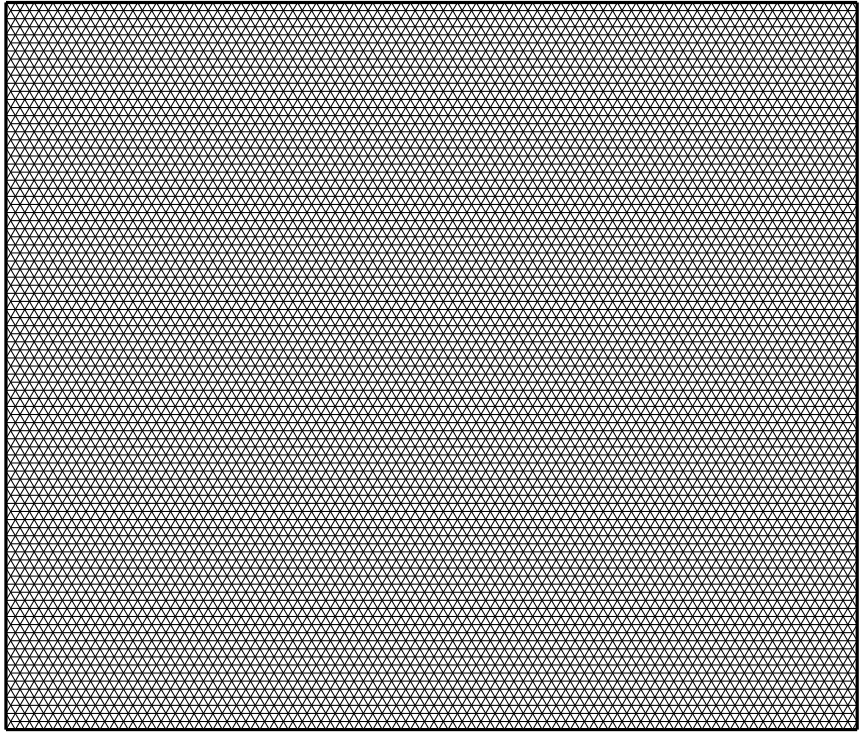

\centering	  
\begin{scaletikzpicturetowidth}{\textwidth*0.70} 
%
\end{scaletikzpicturetowidth}
\caption{Triangular mesh: $90 \times 90$}
\label{fig90}
\end{figure}

In order to investigate this issue we randomly generated a number of test problem instances using mesh sizes of $30 \times 30$, $60 \times 60$ and $90 \times 90$.
The triangular mesh for these test instances are shown in 
Figures~\ref{fig30}-\ref{fig90}. A simple interpretation of the 
$90 \times 90$ mesh
is that if it would take an agent 90 minutes to travel from one side of the area shown in Figure~\ref{fig90} to the other side then one time step corresponds to one minute. Similarly if one side of the area shown in 
Figure~\ref{fig90} is 5 kilometres then one time step corresponds to 55 metres.

Notice how in Figures~\ref{fig30}-\ref{fig90} the triangular discretisation of the area covers the area in greater detail as the mesh size increases. 

The maximum mesh size of 
$90 \times 90$
adopted means that there are approximately $90^2 = 8100$ vertices in the triangular mesh shown in Figure~\ref{fig90}. This in turn implies that for each agent we have of the order of $90^3 = 729000$ $x_{avt}$ variables. Although this seems very large note that as mentioned previously above the use of Equations~(\ref{reduction1}) and~(\ref{reduction2}), which set some $x_{avt}$ variables to zero, significantly reduces the number of variables that need explicit consideration. Also as mentioned above knowing that a specific  $x_{avt}$ is zero enables other reductions to be made, e.g.~from Equation~(\ref{for1}) knowing that $x_{avt}$ is zero means that we also know that $\alpha_{avt}$ is zero.

All of the problems considered in this paper had two agents ($A=2$) and 15 sensors ($S=15$). The reason why we only consider problems with a fixed number of agents and sensors, but varying mesh size, was based on input from the sponsor of this study. In essence their primary concern was mesh size, since the size of the mesh that can be dealt with computationally directly relates to how accurately the topology of the area to be traversed by the agents can be mapped. Mapping this area accurately is clearly of great importance in mission planning. Varying the number of agents and sensors considered, beyond the fixed values adopted, was considered to be much less important.

\subsection{Results}

In order to produce insight into the computational performance as mesh size  increases, we took  five random instances for each mesh size considered, so in total 15 test problem  instances.
We solved these instances for the cases shown in 
Table~\ref{table1}. 
We used a Windows pc with 6GB of memory and an Intel Core i5-2500S 2.5Ghz processor, a multi-core pc with four cores. 

The first two cases in Table~\ref{table1} deal with situations where knockout actions are possible. The third case deals with a situation where we have sensor detection probabilities, but no confusion. The last two cases deal with sensor detection probabilities where confusion actions are possible.
All of the test instances used in Table~\ref{table1} are available on request from the author.

In Table~\ref{table1} we have four columns associated with each of the triangular mesh sizes. Two of these columns relate to the optimal approach given in this paper, so taking the formulation presented and solving it to proven global optimality using Cplex~\cite{cplex1210}. The other two columns relate to the heuristic approach based upon successive shortest paths presented above.  

For the optimal approach the column headed \textbf{\emph{Time}} is the average total computation (CPU) time, in seconds, over the five test problem instances considered.
 The column  headed \textbf{\emph{\%}} is the percentage of this time that was required by the solver (Cplex), the remaining percentage being related to setting up the problem, 
e.g.~computing the information required for variable reduction based on Equations~(\ref{reduction1}) and (\ref{reduction2}). Here, as we used a multi-core pc with four cores, this calculation is based on the CPU time reported by Cplex, not the wall clock time.

For the heuristic approach the column headed \textbf{\emph{Time}} is the average total computation time, in seconds, over the five test problem instances considered.  The column  headed \textbf{\emph{Not opt}} is the number of instances, out of five for each case, where the heuristic terminated with a solution that was suboptimal (i.e.~so with a solution strictly greater than the solution from  the optimal approach).

So, for example, for case 5 and the mesh size of $90 \times 90$, the optimal approach required 11.69 seconds (on average over the five instances considered) with 40\% of that time being consumed by the solver. The average computation time for the heuristic was 42.60 seconds with the heuristic returning a suboptimal solution for one of the five instances.

 The final row in Table~\ref{table1} gives the average percentage of $x_{avt}$ variables excluded, over the five test problem instances considered,  based on applying  Equations~(\ref{reduction1}) and (\ref{reduction2}). These reductions, by their very nature, only apply to the optimal approach.

Considering Table~\ref{table1} it is clear that, as we would expect, as mesh size increases computation time increases. For case 1 the heuristic performs well, requiring on average less computation time than the optimal approach, and only reporting one non-optimal solution. However for the other four cases, and all mesh sizes, shown in 
Table~\ref{table1} it seems reasonable to conclude that the optimal approach dominates the heuristic approach. Over all mesh sizes for cases 2--5 the optimal approach has a lower average computation time. Obviously,   by its very nature, it always finds the optimal solution. By contrast for cases 2--5 there are 9 instances where the heuristic approach fails to find the optimal solution.

With respect to the percentage of computation time consumed by the solver we can see that this is high for case 1, over 70\%, whatever the mesh size. By contrast for cases 2--5 that percentage never exceeds 41\%, and is often much less.

We can also conclude, that even with the largest mesh considered of $90 \times 90$, the computation times seen for the optimal approach are very reasonable. 
Considering the final row in Table~\ref{table1} it is clear that the reduction in the number of variables achieved by making use of  Equations~(\ref{reduction1}) and (\ref{reduction2}) is very significant.

\begin{landscape}
\begin{table}[!htb]
\small
\centering
\renewcommand{\tabcolsep}{1mm} 
\renewcommand{\arraystretch}{1.3} 
\begin{tabular}{|c|l|cc|cc|cc|cc|cc|cc|}
\hline
\multirow{4}*{Case} & \multirow{4}*{Description} & \multicolumn{12}{c|}{Triangular mesh size} \\
\cline{3-14}
 & & \multicolumn{4}{c|}{$30 \times 30$} & \multicolumn{4}{c|}{$60 \times 60$}  & \multicolumn{4}{c|}{$90 \times 90$} \\ 
\cline{3-14}
 & & \multicolumn{2}{c|}{Optimal} & \multicolumn{2}{c|} {Heuristic}  &
\multicolumn{2}{c|}{Optimal} & \multicolumn{2}{c|} {Heuristic}  &
\multicolumn{2}{c|}{Optimal} & \multicolumn{2}{c|} {
Heuristic}  \\
 & & 
Time & \%  & Time & Not opt &
Time & \%  & Time & Not opt &
Time & \%  & Time & Not opt \\

\hline
1 & Minimise time to target &

0.53	& 	77\%	& 	0.27	&	1	&
4.94	& 	71\%	& 	2.81	&	0	&
40.11	& 	79\%	& 	13.59	&	0	

\\  
& ~~1 knockout allowed, $B=1$, $\Omega=2$ &

	& 		& 		&	&
	& 		& 		&		&
	& 		& 		&		
\\
& ~~Agent knockouts are restricted in range and time &  
	& 		& 		&	&
	& 		& 		&		&
	& 		& 		&	
 \\ 

\hline
2 & Minimise time to target &

0.15	&	19\%	& 	0.66	& 	2	&

1.50	&	4\%	& 	4.22	& 	1	&

7.10	&	1\%	& 	21.87	& 	1

\\ 
& ~~2 knockouts allowed, $B=2$, $\Omega=2$ &  
	& 		& 		&	&
	& 		& 		&		&
	& 		& 		&	
\\
& ~~Agent knockouts are restricted in range and time & 
	& 		& 		&	&
	& 		& 		&		&
	& 		& 		&	
\\ 

\hline
3 & Maximise PED, reach target in minimal time &
0.17	&	27\%	&	0.25	&	1	&
1.59	&	9\%	&	2.85	&	1	&

7.57	&	8\%	&	13.51	&	0

\\
& ~~~$T=\mbox{min}[D_{\gamma(a)N},~a=1,\ldots,A]$ &  
	& 		& 		&	&
	& 		& 		&		&
	& 		& 		&	
\\
 & ~~~No confusion, $\Omega=2$ & 
	& 		& 		&	&
	& 		& 		&		&
	& 		& 		&	
 \\

\hline
4 & Maximise PED, reach target in  minimal time &
0.21	&	41\%	&	7.62	&	1	&
2.06	&	30\%	&	148.58	&	1	&
8.43	&	17\%	&	821.25	&	0	

\\
& ~~~$T=\mbox{min}[D_{\gamma(a)N},~a=1,\ldots,A]$ &  
	& 		& 		&	&
	& 		& 		&		&
	& 		& 		&	
\\
& ~~~Confusion, 2 confusion actions allowed, $B=2$, $\Omega=2$ & 
	& 		& 		&	&
	& 		& 		&		&
	& 		& 		&	
\\
  & ~~~Confusion lasts for $\Delta_a^c=10$ time steps  & 
	& 		& 		&	&
	& 		& 		&		&
	& 		& 		&	
\\
& ~~~At most one confusion action every $\Phi_a^c=10$ time steps  & 
	& 		& 		&	&
	& 		& 		&		&
	& 		& 		&	
\\

\hline
5 & Minimise time to target, PED at least 95\% ($Q^* \geq 0.95$) &
0.15	&	19\%	&	1.81	&	0	&
1.76	&	18\%	&	10.08	&	0	&
11.69	&	40\%	&	42.60	&	1	

\\
 & ~~~Confusion, 1 confusion action allowed, $B=1$, $\Omega=2$ & 
	& 		& 		&	&
	& 		& 		&		&
	& 		& 		&	
\\
& ~~~Confusion lasts for $\Delta_a^c=10$ time steps  & 
	& 		& 		&	&
	& 		& 		&		&
	& 		& 		&	
\\

\hline
\multicolumn{2}{|l|}{Percentage of variables excluded by Equations~(\ref{reduction1}) and (\ref{reduction2}) } &
 \multicolumn{2}{c|}{97.98} &  &  & 
 \multicolumn{2}{c|}{99.19} & &  &
 \multicolumn{2}{c|}{99.24} & &  
\\
\hline
\end{tabular}
\caption{Comparative results for the optimal and heuristic approaches with varying mesh sizes}
\label{table1}
\end{table}
\end{landscape}

\normalsize
\subsection{Example solutions}
The figures below show example solutions for the $90 \times 90$ mesh for one of the instances considered in 
Table~\ref{table1}. The example instance shown here was
computationally the most 
time consuming instance of the 15 test instances in  
Table~\ref{table1}.

\begin{figure}[!htb]
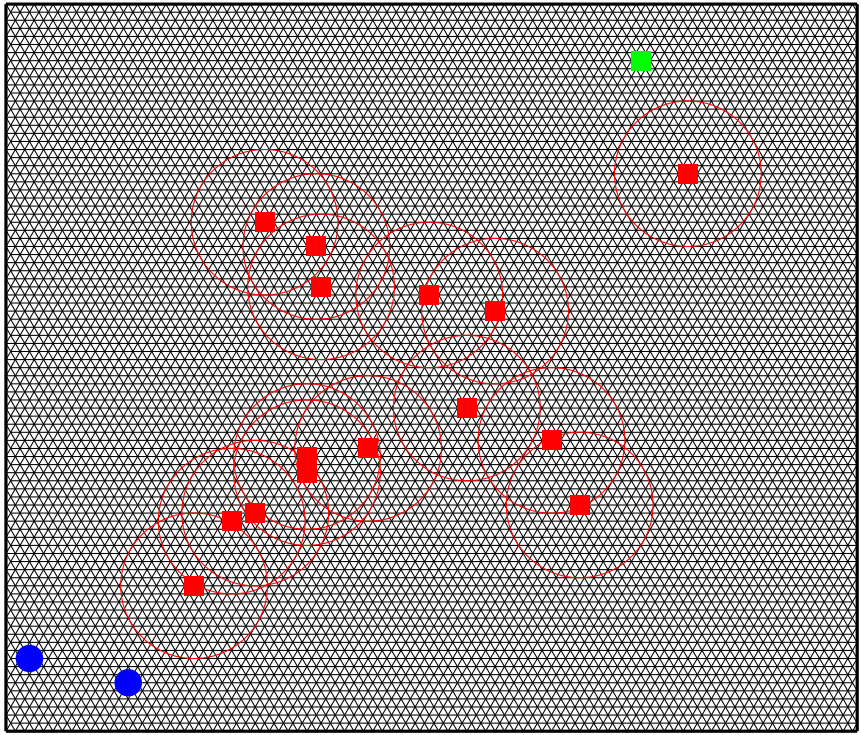

\centering	  
\begin{scaletikzpicturetowidth}{\textwidth*0.70} 

\end{scaletikzpicturetowidth}
\caption{Example problem: $90 \times 90$ mesh shown}
\label{fig12a}
\end{figure}

\begin{figure}[!htb]
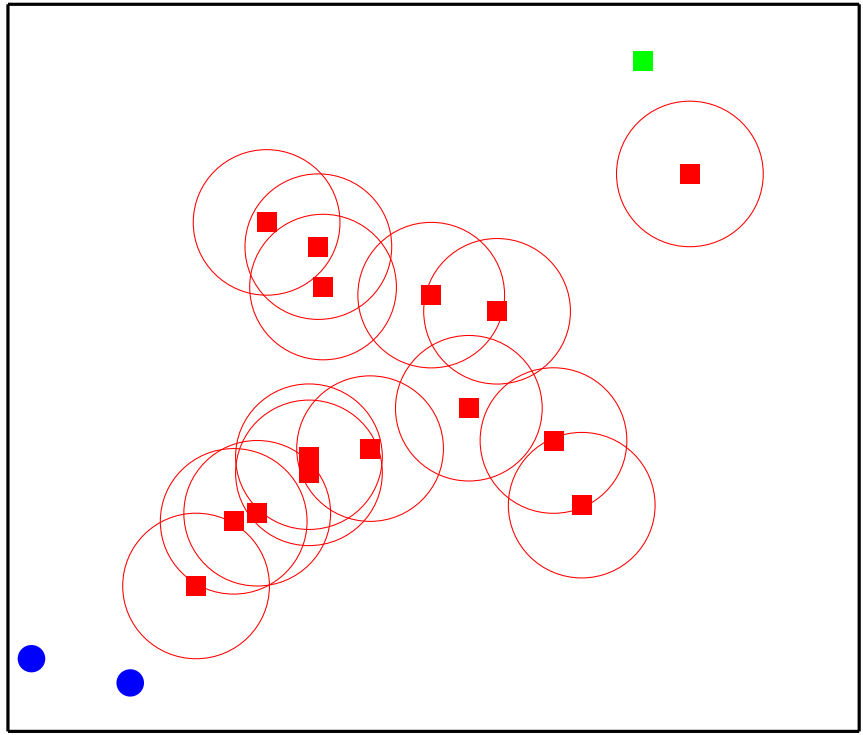

\centering	  
\begin{scaletikzpicturetowidth}{\textwidth*0.70} 
%
\end{scaletikzpicturetowidth}
\caption{Example problem: $90 \times 90$ mesh not shown}
\label{fig345a}
\end{figure}

Figure~\ref{fig12a} shows the example instance we consider here with a $90\times90$ mesh, as in Figure~\ref{fig90}.
Figure~\ref{fig345a} shows the same problem, but  the  mesh is not shown there for clarity. The two agents are shown as solid blue circles and the target as a solid green square. Here there are 15 sensors, each with an associated circular detection region shown in red.
Agent actions, knockout or confusion as appropriate, are shown as solid orange diamonds in the figures below.

\subsubsection{Cases 1,2}

\begin{figure}[!htb]
\centering	  
\begin{scaletikzpicturetowidth}{\textwidth*0.70} 
\begin{tikzpicture}[scale=\tikzscale]	
\draw [red] (  26.50000 ,   23.38269) circle [radius=   7.80000];
\draw [red] (  32.00000 ,   29.44486) circle [radius=   7.80000];
\draw [red] (  72.50000 ,   59.75575) circle [radius=   7.80000];
\draw [red] (  20.00000 ,   15.58846) circle [radius=   7.80000];
\draw [red] (  24.00000 ,   22.51666) circle [radius=   7.80000];
\draw [red] (  33.00000 ,   51.96152) circle [radius=   7.80000];
\draw [red] (  61.00000 ,   24.24871) circle [radius=   7.80000];
\draw [red] (  33.50000 ,   47.63140) circle [radius=   7.80000];
\draw [red] (  27.50000 ,   54.55960) circle [radius=   7.80000];
\draw [red] (  45.00000 ,   46.76537) circle [radius=   7.80000];
\draw [red] (  32.00000 ,   27.71281) circle [radius=   7.80000];
\draw [red] (  52.00000 ,   45.03332) circle [radius=   7.80000];
\draw [red] (  38.50000 ,   30.31089) circle [radius=   7.80000];
\draw [red] (  58.00000 ,   31.17691) circle [radius=   7.80000];
\draw [red] (  49.00000 ,   34.64102) circle [radius=   7.80000];
\node at (  26.50000 ,   23.38269) [rectangle    , fill=red    ] {};
\node at (  32.00000 ,   29.44486) [rectangle    , fill=red    ] {};
\node at ( 72.50000 ,   59.75575) [rectangle    , fill=red    ] {};
\node at (  20.00000 ,   15.58846) [rectangle    , fill=red    ] {};
\node at ( 24.00000 ,   22.51666) [rectangle    , fill=red    ] {};
\node at (  33.00000 ,   51.96152) [rectangle    , fill=red    ] {};
\node at (  61.00000 ,   24.24871) [rectangle    , fill=red    ] {};
\node at ( 33.50000 ,   47.63140) [rectangle    , fill=red    ] {};
\node at (  27.50000 ,   54.55960) [rectangle    , fill=red    ] {};
\node at (  45.00000 ,   46.76537) [rectangle    , fill=red    ] {};
\node at (  32.00000 ,   27.71281) [rectangle    , fill=red    ] {};
\node at (  52.00000 ,   45.03332) [rectangle    , fill=red    ] {};
\node at (  38.50000 ,   30.31089) [rectangle    , fill=red    ] {};
\node at (  58.00000 ,   31.17691) [rectangle    , fill=red    ] {};
\node at (  49.00000 ,   34.64102) [rectangle    , fill=red    ] {};

\draw[very thick] (    .00000 ,     .00000) -- (  90.50000 ,     .00000);
\draw[very thick] (  90.50000 ,     .00000) -- (  90.50000 ,   77.94229);
\draw[very thick] (  90.50000 ,   77.94229) -- (    .00000 ,   77.94229);
\draw[very thick] (    .00000 ,   77.94229) -- (    .00000 ,     .00000);
\node at (  13.00000 ,    5.19615) [circle    , fill=blue     ] {};
\draw [color=blue     ] (  13.00000 ,    5.19615) -- (  14.00000 ,    5.19615);
\draw [color=blue     ] (  14.00000 ,    5.19615) -- (  15.00000 ,    5.19615);
\draw [color=blue     ] (  15.00000 ,    5.19615) -- (  16.00000 ,    5.19615);
\draw [color=blue     ] (  16.00000 ,    5.19615) -- (  17.00000 ,    5.19615);
\draw [color=blue     ] (  17.00000 ,    5.19615) -- (  18.00000 ,    5.19615);
\draw [color=blue     ] (  18.00000 ,    5.19615) -- (  19.00000 ,    5.19615);
\draw [color=blue     ] (  19.00000 ,    5.19615) -- (  20.00000 ,    5.19615);
\draw [color=blue     ] (  20.00000 ,    5.19615) -- (  21.00000 ,    5.19615);
\draw [color=blue     ] (  21.00000 ,    5.19615) -- (  22.00000 ,    5.19615);
\draw [color=blue     ] (  22.00000 ,    5.19615) -- (  23.00000 ,    5.19615);
\draw [color=blue     ] (  23.00000 ,    5.19615) -- (  24.00000 ,    5.19615);
\draw [color=blue     ] (  24.00000 ,    5.19615) -- (  25.00000 ,    5.19615);
\draw [color=blue     ] (  25.00000 ,    5.19615) -- (  26.00000 ,    5.19615);
\draw [color=blue     ] (  26.00000 ,    5.19615) -- (  27.00000 ,    5.19615);
\draw [color=blue     ] (  27.00000 ,    5.19615) -- (  27.50000 ,    6.06218);
\draw [color=blue     ] (  27.50000 ,    6.06218) -- (  28.50000 ,    6.06218);
\draw [color=blue     ] (  28.50000 ,    6.06218) -- (  29.50000 ,    6.06218);
\draw [color=blue     ] (  29.50000 ,    6.06218) -- (  30.50000 ,    6.06218);
\draw [color=blue     ] (  30.50000 ,    6.06218) -- (  31.50000 ,    6.06218);
\draw [color=blue     ] (  31.50000 ,    6.06218) -- (  32.00000 ,    6.92820);
\draw [color=blue     ] (  32.00000 ,    6.92820) -- (  32.50000 ,    7.79423);
\draw [color=blue     ] (  32.50000 ,    7.79423) -- (  33.00000 ,    8.66025);
\draw [color=blue     ] (  33.00000 ,    8.66025) -- (  33.50000 ,    9.52628);
\draw [color=blue     ] (  33.50000 ,    9.52628) -- (  34.00000 ,   10.39230);
\draw [color=blue     ] (  34.00000 ,   10.39230) -- (  34.50000 ,   11.25833);
\draw [color=blue     ] (  34.50000 ,   11.25833) -- (  35.00000 ,   12.12436);
\draw [color=blue     ] (  35.00000 ,   12.12436) -- (  35.50000 ,   12.99038);
\draw [color=blue     ] (  35.50000 ,   12.99038) -- (  36.00000 ,   13.85641);
\draw [color=blue     ] (  36.00000 ,   13.85641) -- (  36.50000 ,   14.72243);
\draw [color=blue     ] (  36.50000 ,   14.72243) -- (  37.00000 ,   15.58846);
\draw [color=blue     ] (  37.00000 ,   15.58846) -- (  37.50000 ,   16.45448);
\draw [color=blue     ] (  37.50000 ,   16.45448) -- (  38.00000 ,   17.32051);
\draw [color=blue     ] (  38.00000 ,   17.32051) -- (  38.50000 ,   18.18653);
\draw [color=blue     ] (  38.50000 ,   18.18653) -- (  39.00000 ,   19.05256);
\draw [color=blue     ] (  39.00000 ,   19.05256) -- (  39.50000 ,   19.91858);
\draw [color=blue     ] (  39.50000 ,   19.91858) -- (  40.00000 ,   20.78461);
\draw [color=blue     ] (  40.00000 ,   20.78461) -- (  40.50000 ,   21.65064);
\draw [color=blue     ] (  40.50000 ,   21.65064) -- (  41.00000 ,   22.51666);
\draw [color=blue     ] (  41.00000 ,   22.51666) -- (  41.50000 ,   23.38269);
\draw [color=blue     ] (  41.50000 ,   23.38269) -- (  42.00000 ,   24.24871);
\draw [color=blue     ] (  42.00000 ,   24.24871) -- (  42.50000 ,   25.11474);
\draw [color=blue     ] (  42.50000 ,   25.11474) -- (  43.00000 ,   25.98076);
\draw [color=blue     ] (  43.00000 ,   25.98076) -- (  43.50000 ,   26.84679);
\draw [color=blue     ] (  43.50000 ,   26.84679) -- (  44.00000 ,   27.71281);
\node at (  44.50000 ,   28.57884) [diamond    , fill=orange   ] {};
\draw [color=blue     ] (  44.00000 ,   27.71281) -- (  44.50000 ,   28.57884);
\draw [color=blue     ] (  44.50000 ,   28.57884) -- (  45.00000 ,   29.44486);
\draw [color=blue     ] (  45.00000 ,   29.44486) -- (  45.50000 ,   30.31089);
\draw [color=blue     ] (  45.50000 ,   30.31089) -- (  46.00000 ,   31.17691);
\draw [color=blue     ] (  46.00000 ,   31.17691) -- (  46.50000 ,   32.04294);
\draw [color=blue     ] (  46.50000 ,   32.04294) -- (  47.00000 ,   32.90897);
\draw [color=blue     ] (  47.00000 ,   32.90897) -- (  47.50000 ,   33.77499);
\draw [color=blue     ] (  47.50000 ,   33.77499) -- (  48.00000 ,   34.64102);
\draw [color=blue     ] (  48.00000 ,   34.64102) -- (  48.50000 ,   35.50704);
\draw [color=blue     ] (  48.50000 ,   35.50704) -- (  49.00000 ,   36.37307);
\draw [color=blue     ] (  49.00000 ,   36.37307) -- (  49.50000 ,   37.23909);
\draw [color=blue     ] (  49.50000 ,   37.23909) -- (  50.00000 ,   38.10512);
\draw [color=blue     ] (  50.00000 ,   38.10512) -- (  50.50000 ,   38.97114);
\draw [color=blue     ] (  50.50000 ,   38.97114) -- (  51.00000 ,   39.83717);
\draw [color=blue     ] (  51.00000 ,   39.83717) -- (  51.50000 ,   40.70319);
\draw [color=blue     ] (  51.50000 ,   40.70319) -- (  52.00000 ,   41.56922);
\draw [color=blue     ] (  52.00000 ,   41.56922) -- (  52.50000 ,   42.43524);
\draw [color=blue     ] (  52.50000 ,   42.43524) -- (  53.00000 ,   43.30127);
\draw [color=blue     ] (  53.00000 ,   43.30127) -- (  53.50000 ,   44.16730);
\draw [color=blue     ] (  53.50000 ,   44.16730) -- (  54.00000 ,   45.03332);
\draw [color=blue     ] (  54.00000 ,   45.03332) -- (  53.50000 ,   45.89935);
\draw [color=blue     ] (  53.50000 ,   45.89935) -- (  53.00000 ,   46.76537);
\draw [color=blue     ] (  53.00000 ,   46.76537) -- (  53.50000 ,   47.63140);
\draw [color=blue     ] (  53.50000 ,   47.63140) -- (  54.00000 ,   48.49742);
\draw [color=blue     ] (  54.00000 ,   48.49742) -- (  54.50000 ,   49.36345);
\draw [color=blue     ] (  54.50000 ,   49.36345) -- (  55.00000 ,   50.22947);
\draw [color=blue     ] (  55.00000 ,   50.22947) -- (  55.50000 ,   51.09550);
\draw [color=blue     ] (  55.50000 ,   51.09550) -- (  56.00000 ,   51.96152);
\draw [color=blue     ] (  56.00000 ,   51.96152) -- (  56.50000 ,   52.82755);
\draw [color=blue     ] (  56.50000 ,   52.82755) -- (  57.00000 ,   53.69358);
\draw [color=blue     ] (  57.00000 ,   53.69358) -- (  57.50000 ,   54.55960);
\draw [color=blue     ] (  57.50000 ,   54.55960) -- (  58.00000 ,   55.42563);
\draw [color=blue     ] (  58.00000 ,   55.42563) -- (  58.50000 ,   56.29165);
\draw [color=blue     ] (  58.50000 ,   56.29165) -- (  59.00000 ,   57.15768);
\draw [color=blue     ] (  59.00000 ,   57.15768) -- (  59.50000 ,   58.02370);
\draw [color=blue     ] (  59.50000 ,   58.02370) -- (  60.00000 ,   58.88973);
\draw [color=blue     ] (  60.00000 ,   58.88973) -- (  60.50000 ,   59.75575);
\draw [color=blue     ] (  60.50000 ,   59.75575) -- (  61.00000 ,   60.62178);
\draw [color=blue     ] (  61.00000 ,   60.62178) -- (  61.50000 ,   61.48780);
\draw [color=blue     ] (  61.50000 ,   61.48780) -- (  62.00000 ,   62.35383);
\draw [color=blue     ] (  62.00000 ,   62.35383) -- (  62.50000 ,   63.21985);
\draw [color=blue     ] (  62.50000 ,   63.21985) -- (  63.00000 ,   64.08588);
\draw [color=blue     ] (  63.00000 ,   64.08588) -- (  63.50000 ,   64.95191);
\draw [color=blue     ] (  63.50000 ,   64.95191) -- (  64.00000 ,   65.81793);
\draw [color=blue     ] (  64.00000 ,   65.81793) -- (  64.50000 ,   66.68396);
\draw [color=blue     ] (  64.50000 ,   66.68396) -- (  65.00000 ,   67.54998);
\draw [color=blue     ] (  65.00000 ,   67.54998) -- (  65.50000 ,   68.41601);
\draw [color=blue     ] (  65.50000 ,   68.41601) -- (  66.00000 ,   69.28203);
\draw [color=blue     ] (  66.00000 ,   69.28203) -- (  66.50000 ,   70.14806);
\draw [color=blue     ] (  66.50000 ,   70.14806) -- (  67.00000 ,   71.01408);
\draw [color=blue     ] (  67.00000 ,   71.01408) -- (  67.50000 ,   71.88011);
\node at (  67.50000 ,   71.88011) [rectangle    , fill=green    ] {};
\end{tikzpicture}
\end{scaletikzpicturetowidth}
\caption{Case 1: Time to target 95 time steps, agent knockout at time step 45}
\label{fig12b}
\end{figure}
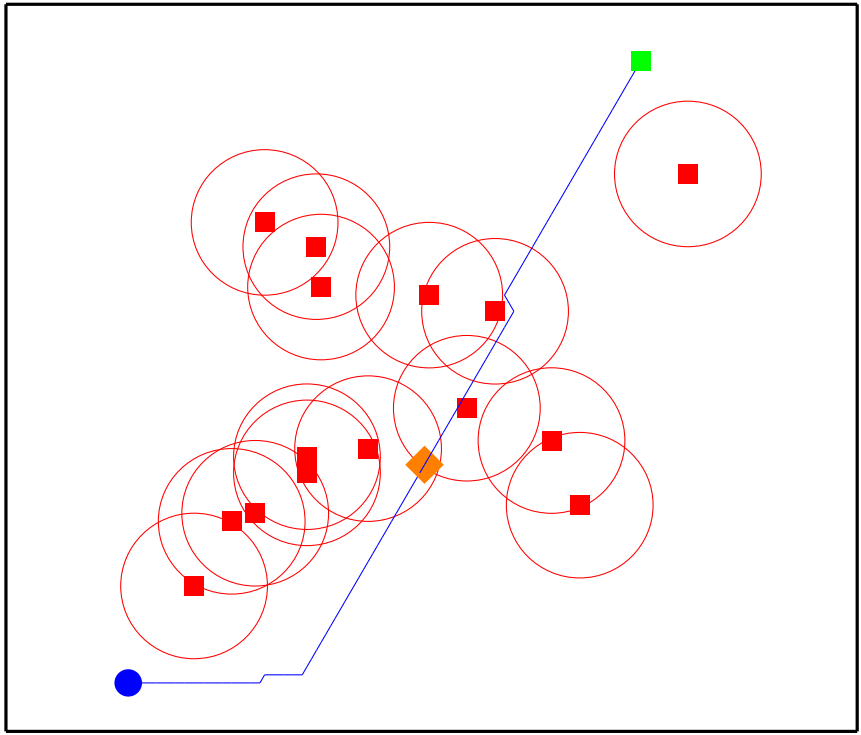

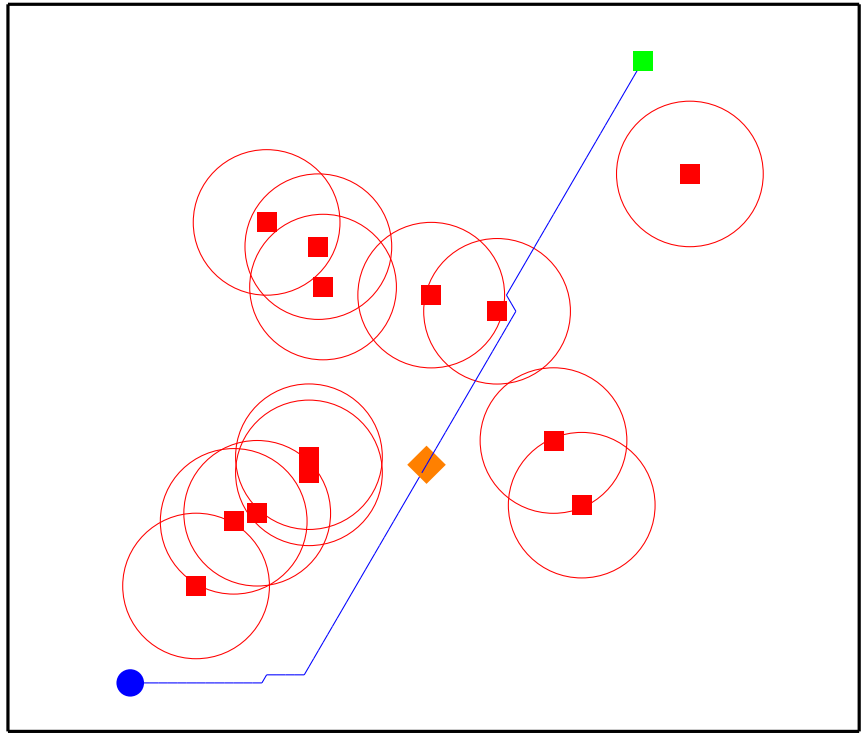
\begin{figure}[!htb]
\centering	  
\begin{scaletikzpicturetowidth}{\textwidth*0.70} 
\begin{tikzpicture}[scale=\tikzscale]	
\draw [red] (  26.50000 ,   23.38269) circle [radius=   7.80000];
\draw [red] (  32.00000 ,   29.44486) circle [radius=   7.80000];
\draw [red] (  72.50000 ,   59.75575) circle [radius=   7.80000];
\draw [red] (  20.00000 ,   15.58846) circle [radius=   7.80000];
\draw [red] (  24.00000 ,   22.51666) circle [radius=   7.80000];
\draw [red] (  33.00000 ,   51.96152) circle [radius=   7.80000];
\draw [red] (  61.00000 ,   24.24871) circle [radius=   7.80000];
\draw [red] (  33.50000 ,   47.63140) circle [radius=   7.80000];
\draw [red] (  27.50000 ,   54.55960) circle [radius=   7.80000];
\draw [red] (  45.00000 ,   46.76537) circle [radius=   7.80000];
\draw [red] (  32.00000 ,   27.71281) circle [radius=   7.80000];
\draw [red] (  52.00000 ,   45.03332) circle [radius=   7.80000];
\draw [red] (  58.00000 ,   31.17691) circle [radius=   7.80000];
\node at (  26.50000 ,   23.38269) [rectangle    , fill=red    ] {};
\node at (  32.00000 ,   29.44486) [rectangle    , fill=red    ] {};
\node at ( 72.50000 ,   59.75575) [rectangle    , fill=red    ] {};
\node at (  20.00000 ,   15.58846) [rectangle    , fill=red    ] {};
\node at ( 24.00000 ,   22.51666) [rectangle    , fill=red    ] {};
\node at (  33.00000 ,   51.96152) [rectangle    , fill=red    ] {};
\node at (  61.00000 ,   24.24871) [rectangle    , fill=red    ] {};
\node at ( 33.50000 ,   47.63140) [rectangle    , fill=red    ] {};
\node at (  27.50000 ,   54.55960) [rectangle    , fill=red    ] {};
\node at (  45.00000 ,   46.76537) [rectangle    , fill=red    ] {};
\node at (  32.00000 ,   27.71281) [rectangle    , fill=red    ] {};
\node at (  52.00000 ,   45.03332) [rectangle    , fill=red    ] {};

\node at (  58.00000 ,   31.17691) [rectangle    , fill=red    ] {};


\draw[very thick] (    .00000 ,     .00000) -- (  90.50000 ,     .00000);
\draw[very thick] (  90.50000 ,     .00000) -- (  90.50000 ,   77.94229);
\draw[very thick] (  90.50000 ,   77.94229) -- (    .00000 ,   77.94229);
\draw[very thick] (    .00000 ,   77.94229) -- (    .00000 ,     .00000);
\node at (  13.00000 ,    5.19615) [circle    , fill=blue     ] {};
\draw [color=blue     ] (  13.00000 ,    5.19615) -- (  14.00000 ,    5.19615);
\draw [color=blue     ] (  14.00000 ,    5.19615) -- (  15.00000 ,    5.19615);
\draw [color=blue     ] (  15.00000 ,    5.19615) -- (  16.00000 ,    5.19615);
\draw [color=blue     ] (  16.00000 ,    5.19615) -- (  17.00000 ,    5.19615);
\draw [color=blue     ] (  17.00000 ,    5.19615) -- (  18.00000 ,    5.19615);
\draw [color=blue     ] (  18.00000 ,    5.19615) -- (  19.00000 ,    5.19615);
\draw [color=blue     ] (  19.00000 ,    5.19615) -- (  20.00000 ,    5.19615);
\draw [color=blue     ] (  20.00000 ,    5.19615) -- (  21.00000 ,    5.19615);
\draw [color=blue     ] (  21.00000 ,    5.19615) -- (  22.00000 ,    5.19615);
\draw [color=blue     ] (  22.00000 ,    5.19615) -- (  23.00000 ,    5.19615);
\draw [color=blue     ] (  23.00000 ,    5.19615) -- (  24.00000 ,    5.19615);
\draw [color=blue     ] (  24.00000 ,    5.19615) -- (  25.00000 ,    5.19615);
\draw [color=blue     ] (  25.00000 ,    5.19615) -- (  26.00000 ,    5.19615);
\draw [color=blue     ] (  26.00000 ,    5.19615) -- (  27.00000 ,    5.19615);
\draw [color=blue     ] (  27.00000 ,    5.19615) -- (  27.50000 ,    6.06218);
\draw [color=blue     ] (  27.50000 ,    6.06218) -- (  28.50000 ,    6.06218);
\draw [color=blue     ] (  28.50000 ,    6.06218) -- (  29.50000 ,    6.06218);
\draw [color=blue     ] (  29.50000 ,    6.06218) -- (  30.50000 ,    6.06218);
\draw [color=blue     ] (  30.50000 ,    6.06218) -- (  31.50000 ,    6.06218);
\draw [color=blue     ] (  31.50000 ,    6.06218) -- (  32.00000 ,    6.92820);
\draw [color=blue     ] (  32.00000 ,    6.92820) -- (  32.50000 ,    7.79423);
\draw [color=blue     ] (  32.50000 ,    7.79423) -- (  33.00000 ,    8.66025);
\draw [color=blue     ] (  33.00000 ,    8.66025) -- (  33.50000 ,    9.52628);
\draw [color=blue     ] (  33.50000 ,    9.52628) -- (  34.00000 ,   10.39230);
\draw [color=blue     ] (  34.00000 ,   10.39230) -- (  34.50000 ,   11.25833);
\draw [color=blue     ] (  34.50000 ,   11.25833) -- (  35.00000 ,   12.12436);
\draw [color=blue     ] (  35.00000 ,   12.12436) -- (  35.50000 ,   12.99038);
\draw [color=blue     ] (  35.50000 ,   12.99038) -- (  36.00000 ,   13.85641);
\draw [color=blue     ] (  36.00000 ,   13.85641) -- (  36.50000 ,   14.72243);
\draw [color=blue     ] (  36.50000 ,   14.72243) -- (  37.00000 ,   15.58846);
\draw [color=blue     ] (  37.00000 ,   15.58846) -- (  37.50000 ,   16.45448);
\draw [color=blue     ] (  37.50000 ,   16.45448) -- (  38.00000 ,   17.32051);
\draw [color=blue     ] (  38.00000 ,   17.32051) -- (  38.50000 ,   18.18653);
\draw [color=blue     ] (  38.50000 ,   18.18653) -- (  39.00000 ,   19.05256);
\draw [color=blue     ] (  39.00000 ,   19.05256) -- (  39.50000 ,   19.91858);
\draw [color=blue     ] (  39.50000 ,   19.91858) -- (  40.00000 ,   20.78461);
\draw [color=blue     ] (  40.00000 ,   20.78461) -- (  40.50000 ,   21.65064);
\draw [color=blue     ] (  40.50000 ,   21.65064) -- (  41.00000 ,   22.51666);
\draw [color=blue     ] (  41.00000 ,   22.51666) -- (  41.50000 ,   23.38269);
\draw [color=blue     ] (  41.50000 ,   23.38269) -- (  42.00000 ,   24.24871);
\draw [color=blue     ] (  42.00000 ,   24.24871) -- (  42.50000 ,   25.11474);
\draw [color=blue     ] (  42.50000 ,   25.11474) -- (  43.00000 ,   25.98076);
\draw [color=blue     ] (  43.00000 ,   25.98076) -- (  43.50000 ,   26.84679);
\draw [color=blue     ] (  43.50000 ,   26.84679) -- (  44.00000 ,   27.71281);
\node at (  44.50000 ,   28.57884) [diamond   , fill=orange   ] {};
\draw [color=blue     ] (  44.00000 ,   27.71281) -- (  44.50000 ,   28.57884);
\draw [color=blue     ] (  44.50000 ,   28.57884) -- (  45.00000 ,   29.44486);
\draw [color=blue     ] (  45.00000 ,   29.44486) -- (  45.50000 ,   30.31089);
\draw [color=blue     ] (  45.50000 ,   30.31089) -- (  46.00000 ,   31.17691);
\draw [color=blue     ] (  46.00000 ,   31.17691) -- (  46.50000 ,   32.04294);
\draw [color=blue     ] (  46.50000 ,   32.04294) -- (  47.00000 ,   32.90897);
\draw [color=blue     ] (  47.00000 ,   32.90897) -- (  47.50000 ,   33.77499);
\draw [color=blue     ] (  47.50000 ,   33.77499) -- (  48.00000 ,   34.64102);
\draw [color=blue     ] (  48.00000 ,   34.64102) -- (  48.50000 ,   35.50704);
\draw [color=blue     ] (  48.50000 ,   35.50704) -- (  49.00000 ,   36.37307);
\draw [color=blue     ] (  49.00000 ,   36.37307) -- (  49.50000 ,   37.23909);
\draw [color=blue     ] (  49.50000 ,   37.23909) -- (  50.00000 ,   38.10512);
\draw [color=blue     ] (  50.00000 ,   38.10512) -- (  50.50000 ,   38.97114);
\draw [color=blue     ] (  50.50000 ,   38.97114) -- (  51.00000 ,   39.83717);
\draw [color=blue     ] (  51.00000 ,   39.83717) -- (  51.50000 ,   40.70319);
\draw [color=blue     ] (  51.50000 ,   40.70319) -- (  52.00000 ,   41.56922);
\draw [color=blue     ] (  52.00000 ,   41.56922) -- (  52.50000 ,   42.43524);
\draw [color=blue     ] (  52.50000 ,   42.43524) -- (  53.00000 ,   43.30127);
\draw [color=blue     ] (  53.00000 ,   43.30127) -- (  53.50000 ,   44.16730);
\draw [color=blue     ] (  53.50000 ,   44.16730) -- (  54.00000 ,   45.03332);
\draw [color=blue     ] (  54.00000 ,   45.03332) -- (  53.50000 ,   45.89935);
\draw [color=blue     ] (  53.50000 ,   45.89935) -- (  53.00000 ,   46.76537);
\draw [color=blue     ] (  53.00000 ,   46.76537) -- (  53.50000 ,   47.63140);
\draw [color=blue     ] (  53.50000 ,   47.63140) -- (  54.00000 ,   48.49742);
\draw [color=blue     ] (  54.00000 ,   48.49742) -- (  54.50000 ,   49.36345);
\draw [color=blue     ] (  54.50000 ,   49.36345) -- (  55.00000 ,   50.22947);
\draw [color=blue     ] (  55.00000 ,   50.22947) -- (  55.50000 ,   51.09550);
\draw [color=blue     ] (  55.50000 ,   51.09550) -- (  56.00000 ,   51.96152);
\draw [color=blue     ] (  56.00000 ,   51.96152) -- (  56.50000 ,   52.82755);
\draw [color=blue     ] (  56.50000 ,   52.82755) -- (  57.00000 ,   53.69358);
\draw [color=blue     ] (  57.00000 ,   53.69358) -- (  57.50000 ,   54.55960);
\draw [color=blue     ] (  57.50000 ,   54.55960) -- (  58.00000 ,   55.42563);
\draw [color=blue     ] (  58.00000 ,   55.42563) -- (  58.50000 ,   56.29165);
\draw [color=blue     ] (  58.50000 ,   56.29165) -- (  59.00000 ,   57.15768);
\draw [color=blue     ] (  59.00000 ,   57.15768) -- (  59.50000 ,   58.02370);
\draw [color=blue     ] (  59.50000 ,   58.02370) -- (  60.00000 ,   58.88973);
\draw [color=blue     ] (  60.00000 ,   58.88973) -- (  60.50000 ,   59.75575);
\draw [color=blue     ] (  60.50000 ,   59.75575) -- (  61.00000 ,   60.62178);
\draw [color=blue     ] (  61.00000 ,   60.62178) -- (  61.50000 ,   61.48780);
\draw [color=blue     ] (  61.50000 ,   61.48780) -- (  62.00000 ,   62.35383);
\draw [color=blue     ] (  62.00000 ,   62.35383) -- (  62.50000 ,   63.21985);
\draw [color=blue     ] (  62.50000 ,   63.21985) -- (  63.00000 ,   64.08588);
\draw [color=blue     ] (  63.00000 ,   64.08588) -- (  63.50000 ,   64.95191);
\draw [color=blue     ] (  63.50000 ,   64.95191) -- (  64.00000 ,   65.81793);
\draw [color=blue     ] (  64.00000 ,   65.81793) -- (  64.50000 ,   66.68396);
\draw [color=blue     ] (  64.50000 ,   66.68396) -- (  65.00000 ,   67.54998);
\draw [color=blue     ] (  65.00000 ,   67.54998) -- (  65.50000 ,   68.41601);
\draw [color=blue     ] (  65.50000 ,   68.41601) -- (  66.00000 ,   69.28203);
\draw [color=blue     ] (  66.00000 ,   69.28203) -- (  66.50000 ,   70.14806);
\draw [color=blue     ] (  66.50000 ,   70.14806) -- (  67.00000 ,   71.01408);
\draw [color=blue     ] (  67.00000 ,   71.01408) -- (  67.50000 ,   71.88011);
\node at (  67.50000 ,   71.88011) [rectangle    , fill=green    ] {};
\end{tikzpicture}
\end{scaletikzpicturetowidth}
\caption{Case 1: Time to target 95 time steps, agent knockout at time step 45, knocked out sensors not shown}
\label{fig12c}
\end{figure}

For case 1, where one knockout is allowed, 
Figure~\ref{fig12b} shows the solution. The minimal time to target when only one knockout is allowed is 95 time steps. Figure~\ref{fig12c} shows the same solution as 
Figure~\ref{fig12b}, but with the sensors that have been knocked out removed for clarity.

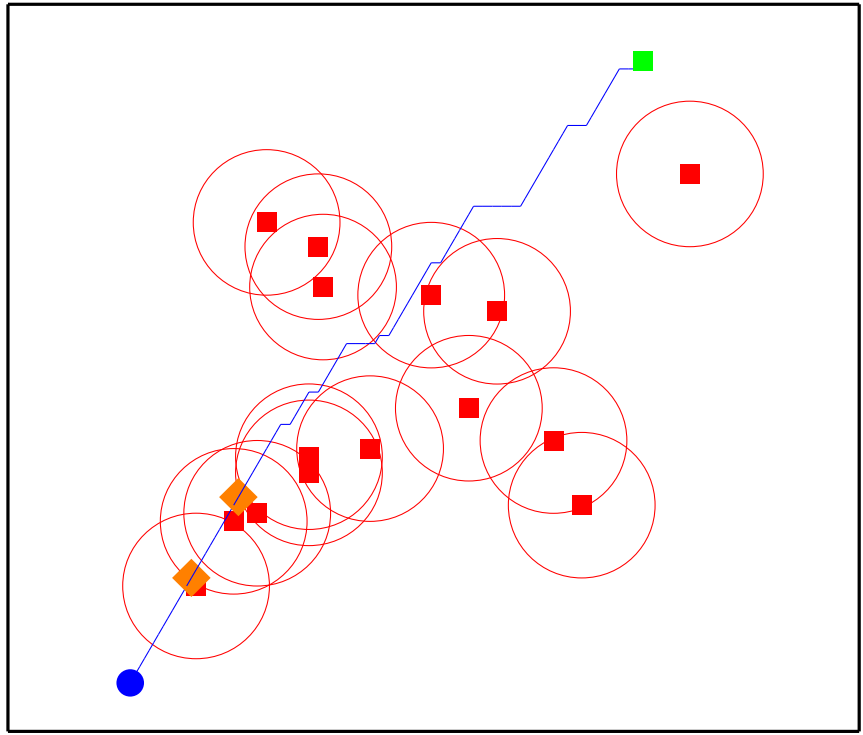
\begin{figure}[!htb]
\centering	  
\begin{scaletikzpicturetowidth}{\textwidth*0.70} 
\begin{tikzpicture}[scale=\tikzscale]	
\draw [red] (  26.50000 ,   23.38269) circle [radius=   7.80000];
\draw [red] (  32.00000 ,   29.44486) circle [radius=   7.80000];
\draw [red] (  72.50000 ,   59.75575) circle [radius=   7.80000];
\draw [red] (  20.00000 ,   15.58846) circle [radius=   7.80000];
\draw [red] (  24.00000 ,   22.51666) circle [radius=   7.80000];
\draw [red] (  33.00000 ,   51.96152) circle [radius=   7.80000];
\draw [red] (  61.00000 ,   24.24871) circle [radius=   7.80000];
\draw [red] (  33.50000 ,   47.63140) circle [radius=   7.80000];
\draw [red] (  27.50000 ,   54.55960) circle [radius=   7.80000];
\draw [red] (  45.00000 ,   46.76537) circle [radius=   7.80000];
\draw [red] (  32.00000 ,   27.71281) circle [radius=   7.80000];
\draw [red] (  52.00000 ,   45.03332) circle [radius=   7.80000];
\draw [red] (  38.50000 ,   30.31089) circle [radius=   7.80000];
\draw [red] (  58.00000 ,   31.17691) circle [radius=   7.80000];
\draw [red] (  49.00000 ,   34.64102) circle [radius=   7.80000];

\node at (  26.50000 ,   23.38269) [rectangle    , fill=red    ] {};
\node at (  32.00000 ,   29.44486) [rectangle    , fill=red    ] {};
\node at ( 72.50000 ,   59.75575) [rectangle    , fill=red    ] {};
\node at (  20.00000 ,   15.58846) [rectangle    , fill=red    ] {};
\node at ( 24.00000 ,   22.51666) [rectangle    , fill=red    ] {};
\node at (  33.00000 ,   51.96152) [rectangle    , fill=red    ] {};
\node at (  61.00000 ,   24.24871) [rectangle    , fill=red    ] {};
\node at ( 33.50000 ,   47.63140) [rectangle    , fill=red    ] {};
\node at (  27.50000 ,   54.55960) [rectangle    , fill=red    ] {};
\node at (  45.00000 ,   46.76537) [rectangle    , fill=red    ] {};
\node at (  32.00000 ,   27.71281) [rectangle    , fill=red    ] {};
\node at (  52.00000 ,   45.03332) [rectangle    , fill=red    ] {};
\node at (  38.50000 ,   30.31089) [rectangle    , fill=red    ] {};
\node at (  58.00000 ,   31.17691) [rectangle    , fill=red    ] {};
\node at (  49.00000 ,   34.64102) [rectangle    , fill=red    ] {};

\draw[very thick] (    .00000 ,     .00000) -- (  90.50000 ,     .00000);
\draw[very thick] (  90.50000 ,     .00000) -- (  90.50000 ,   77.94229);
\draw[very thick] (  90.50000 ,   77.94229) -- (    .00000 ,   77.94229);
\draw[very thick] (    .00000 ,   77.94229) -- (    .00000 ,     .00000);
\node at (  13.00000 ,    5.19615) [circle    , fill=blue     ] {};
\draw [color=blue     ] (  13.00000 ,    5.19615) -- (  13.50000 ,    6.06218);
\draw [color=blue     ] (  13.50000 ,    6.06218) -- (  14.00000 ,    6.92820);
\draw [color=blue     ] (  14.00000 ,    6.92820) -- (  14.50000 ,    7.79423);
\draw [color=blue     ] (  14.50000 ,    7.79423) -- (  15.00000 ,    8.66025);
\draw [color=blue     ] (  15.00000 ,    8.66025) -- (  15.50000 ,    9.52628);
\draw [color=blue     ] (  15.50000 ,    9.52628) -- (  16.00000 ,   10.39230);
\draw [color=blue     ] (  16.00000 ,   10.39230) -- (  16.50000 ,   11.25833);
\draw [color=blue     ] (  16.50000 ,   11.25833) -- (  17.00000 ,   12.12436);
\draw [color=blue     ] (  17.00000 ,   12.12436) -- (  17.50000 ,   12.99038);
\draw [color=blue     ] (  17.50000 ,   12.99038) -- (  18.00000 ,   13.85641);
\draw [color=blue     ] (  18.00000 ,   13.85641) -- (  18.50000 ,   14.72243);
\draw [color=blue     ] (  18.50000 ,   14.72243) -- (  19.00000 ,   15.58846);
\node at (  19.50000 ,   16.45448) [diamond, fill=orange   ] {};
\draw [color=blue     ] (  19.00000 ,   15.58846) -- (  19.50000 ,   16.45448);
\draw [color=blue     ] (  19.50000 ,   16.45448) -- (  20.00000 ,   17.32051);
\draw [color=blue     ] (  20.00000 ,   17.32051) -- (  20.50000 ,   18.18653);
\draw [color=blue     ] (  20.50000 ,   18.18653) -- (  21.00000 ,   19.05256);
\draw [color=blue     ] (  21.00000 ,   19.05256) -- (  21.50000 ,   19.91858);
\draw [color=blue     ] (  21.50000 ,   19.91858) -- (  22.00000 ,   20.78461);
\draw [color=blue     ] (  22.00000 ,   20.78461) -- (  22.50000 ,   21.65064);
\draw [color=blue     ] (  22.50000 ,   21.65064) -- (  23.00000 ,   22.51666);
\draw [color=blue     ] (  23.00000 ,   22.51666) -- (  23.50000 ,   23.38269);
\draw [color=blue     ] (  23.50000 ,   23.38269) -- (  24.00000 ,   24.24871);
\node at (  24.50000 ,   25.11474) [diamond    , fill=orange   ] {};
\draw [color=blue     ] (  24.00000 ,   24.24871) -- (  24.50000 ,   25.11474);
\draw [color=blue     ] (  24.50000 ,   25.11474) -- (  25.00000 ,   25.98076);
\draw [color=blue     ] (  25.00000 ,   25.98076) -- (  25.50000 ,   26.84679);
\draw [color=blue     ] (  25.50000 ,   26.84679) -- (  26.00000 ,   27.71281);
\draw [color=blue     ] (  26.00000 ,   27.71281) -- (  26.50000 ,   28.57884);
\draw [color=blue     ] (  26.50000 ,   28.57884) -- (  27.00000 ,   29.44486);
\draw [color=blue     ] (  27.00000 ,   29.44486) -- (  27.50000 ,   30.31089);
\draw [color=blue     ] (  27.50000 ,   30.31089) -- (  28.00000 ,   31.17691);
\draw [color=blue     ] (  28.00000 ,   31.17691) -- (  28.50000 ,   32.04294);
\draw [color=blue     ] (  28.50000 ,   32.04294) -- (  29.00000 ,   32.90897);
\draw [color=blue     ] (  29.00000 ,   32.90897) -- (  30.00000 ,   32.90897);
\draw [color=blue     ] (  30.00000 ,   32.90897) -- (  30.50000 ,   33.77499);
\draw [color=blue     ] (  30.50000 ,   33.77499) -- (  31.00000 ,   34.64102);
\draw [color=blue     ] (  31.00000 ,   34.64102) -- (  31.50000 ,   35.50704);
\draw [color=blue     ] (  31.50000 ,   35.50704) -- (  32.00000 ,   36.37307);
\draw [color=blue     ] (  32.00000 ,   36.37307) -- (  33.00000 ,   36.37307);
\draw [color=blue     ] (  33.00000 ,   36.37307) -- (  33.50000 ,   37.23909);
\draw [color=blue     ] (  33.50000 ,   37.23909) -- (  34.00000 ,   38.10512);
\draw [color=blue     ] (  34.00000 ,   38.10512) -- (  34.50000 ,   38.97114);
\draw [color=blue     ] (  34.50000 ,   38.97114) -- (  35.00000 ,   39.83717);
\draw [color=blue     ] (  35.00000 ,   39.83717) -- (  35.50000 ,   40.70319);
\draw [color=blue     ] (  35.50000 ,   40.70319) -- (  36.00000 ,   41.56922);
\draw [color=blue     ] (  36.00000 ,   41.56922) -- (  37.00000 ,   41.56922);
\draw [color=blue     ] (  37.00000 ,   41.56922) -- (  38.00000 ,   41.56922);
\draw [color=blue     ] (  38.00000 ,   41.56922) -- (  39.00000 ,   41.56922);
\draw [color=blue     ] (  39.00000 ,   41.56922) -- (  39.50000 ,   42.43524);
\draw [color=blue     ] (  39.50000 ,   42.43524) -- (  40.50000 ,   42.43524);
\draw [color=blue     ] (  40.50000 ,   42.43524) -- (  41.00000 ,   43.30127);
\draw [color=blue     ] (  41.00000 ,   43.30127) -- (  41.50000 ,   44.16730);
\draw [color=blue     ] (  41.50000 ,   44.16730) -- (  42.00000 ,   45.03332);
\draw [color=blue     ] (  42.00000 ,   45.03332) -- (  42.50000 ,   45.89935);
\draw [color=blue     ] (  42.50000 ,   45.89935) -- (  43.00000 ,   46.76537);
\draw [color=blue     ] (  43.00000 ,   46.76537) -- (  43.50000 ,   47.63140);
\draw [color=blue     ] (  43.50000 ,   47.63140) -- (  44.00000 ,   48.49742);
\draw [color=blue     ] (  44.00000 ,   48.49742) -- (  44.50000 ,   49.36345);
\draw [color=blue     ] (  44.50000 ,   49.36345) -- (  45.00000 ,   50.22947);
\draw [color=blue     ] (  45.00000 ,   50.22947) -- (  46.00000 ,   50.22947);
\draw [color=blue     ] (  46.00000 ,   50.22947) -- (  46.50000 ,   51.09550);
\draw [color=blue     ] (  46.50000 ,   51.09550) -- (  47.00000 ,   51.96152);
\draw [color=blue     ] (  47.00000 ,   51.96152) -- (  47.50000 ,   52.82755);
\draw [color=blue     ] (  47.50000 ,   52.82755) -- (  48.00000 ,   53.69358);
\draw [color=blue     ] (  48.00000 ,   53.69358) -- (  48.50000 ,   54.55960);
\draw [color=blue     ] (  48.50000 ,   54.55960) -- (  49.00000 ,   55.42563);
\draw [color=blue     ] (  49.00000 ,   55.42563) -- (  49.50000 ,   56.29165);
\draw [color=blue     ] (  49.50000 ,   56.29165) -- (  50.50000 ,   56.29165);
\draw [color=blue     ] (  50.50000 ,   56.29165) -- (  51.50000 ,   56.29165);
\draw [color=blue     ] (  51.50000 ,   56.29165) -- (  52.50000 ,   56.29165);
\draw [color=blue     ] (  52.50000 ,   56.29165) -- (  53.50000 ,   56.29165);
\draw [color=blue     ] (  53.50000 ,   56.29165) -- (  54.50000 ,   56.29165);
\draw [color=blue     ] (  54.50000 ,   56.29165) -- (  55.00000 ,   57.15768);
\draw [color=blue     ] (  55.00000 ,   57.15768) -- (  55.50000 ,   58.02370);
\draw [color=blue     ] (  55.50000 ,   58.02370) -- (  56.00000 ,   58.88973);
\draw [color=blue     ] (  56.00000 ,   58.88973) -- (  56.50000 ,   59.75575);
\draw [color=blue     ] (  56.50000 ,   59.75575) -- (  57.00000 ,   60.62178);
\draw [color=blue     ] (  57.00000 ,   60.62178) -- (  57.50000 ,   61.48780);
\draw [color=blue     ] (  57.50000 ,   61.48780) -- (  58.00000 ,   62.35383);
\draw [color=blue     ] (  58.00000 ,   62.35383) -- (  58.50000 ,   63.21985);
\draw [color=blue     ] (  58.50000 ,   63.21985) -- (  59.00000 ,   64.08588);
\draw [color=blue     ] (  59.00000 ,   64.08588) -- (  59.50000 ,   64.95191);
\draw [color=blue     ] (  59.50000 ,   64.95191) -- (  60.50000 ,   64.95191);
\draw [color=blue     ] (  60.50000 ,   64.95191) -- (  61.50000 ,   64.95191);
\draw [color=blue     ] (  61.50000 ,   64.95191) -- (  62.00000 ,   65.81793);
\draw [color=blue     ] (  62.00000 ,   65.81793) -- (  62.50000 ,   66.68396);
\draw [color=blue     ] (  62.50000 ,   66.68396) -- (  63.00000 ,   67.54998);
\draw [color=blue     ] (  63.00000 ,   67.54998) -- (  63.50000 ,   68.41601);
\draw [color=blue     ] (  63.50000 ,   68.41601) -- (  64.00000 ,   69.28203);
\draw [color=blue     ] (  64.00000 ,   69.28203) -- (  64.50000 ,   70.14806);
\draw [color=blue     ] (  64.50000 ,   70.14806) -- (  65.00000 ,   71.01408);
\draw [color=blue     ] (  65.00000 ,   71.01408) -- (  66.00000 ,   71.01408);
\draw [color=blue     ] (  66.00000 ,   71.01408) -- (  67.00000 ,   71.01408);
\draw [color=blue     ] (  67.00000 ,   71.01408) -- (  67.50000 ,   71.88011);
\node at (  67.50000 ,   71.88011) [rectangle    , fill=green    ] {};
\end{tikzpicture}
\end{scaletikzpicturetowidth}
\caption{Case 2: Time to target 93 time steps, agent knockouts at time steps 13 and 23}
\label{fig12d}
\end{figure}

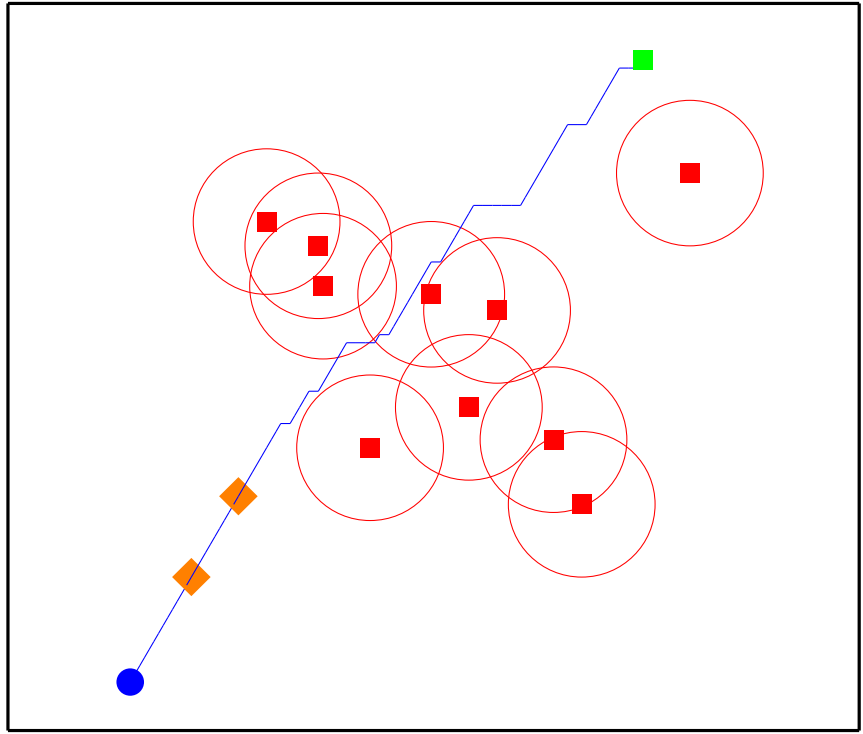
\begin{figure}[!htb]
\centering	  
\begin{scaletikzpicturetowidth}{\textwidth*0.70} 
\begin{tikzpicture}[scale=\tikzscale]	
\draw [red] (  72.50000 ,   59.75575) circle [radius=   7.80000];
\draw [red] (  33.00000 ,   51.96152) circle [radius=   7.80000];
\draw [red] (  61.00000 ,   24.24871) circle [radius=   7.80000];
\draw [red] (  33.50000 ,   47.63140) circle [radius=   7.80000];
\draw [red] (  27.50000 ,   54.55960) circle [radius=   7.80000];
\draw [red] (  45.00000 ,   46.76537) circle [radius=   7.80000];
\draw [red] (  52.00000 ,   45.03332) circle [radius=   7.80000];
\draw [red] (  38.50000 ,   30.31089) circle [radius=   7.80000];
\draw [red] (  58.00000 ,   31.17691) circle [radius=   7.80000];
\draw [red] (  49.00000 ,   34.64102) circle [radius=   7.80000];


\node at ( 72.50000 ,   59.75575) [rectangle    , fill=red    ] {};

\node at (  33.00000 ,   51.96152) [rectangle    , fill=red    ] {};
\node at (  61.00000 ,   24.24871) [rectangle    , fill=red    ] {};
\node at ( 33.50000 ,   47.63140) [rectangle    , fill=red    ] {};
\node at (  27.50000 ,   54.55960) [rectangle    , fill=red    ] {};
\node at (  45.00000 ,   46.76537) [rectangle    , fill=red    ] {};

\node at (  52.00000 ,   45.03332) [rectangle    , fill=red    ] {};
\node at (  38.50000 ,   30.31089) [rectangle    , fill=red    ] {};
\node at (  58.00000 ,   31.17691) [rectangle    , fill=red    ] {};
\node at (  49.00000 ,   34.64102) [rectangle    , fill=red    ] {};

\draw[very thick] (    .00000 ,     .00000) -- (  90.50000 ,     .00000);
\draw[very thick] (  90.50000 ,     .00000) -- (  90.50000 ,   77.94229);
\draw[very thick] (  90.50000 ,   77.94229) -- (    .00000 ,   77.94229);
\draw[very thick] (    .00000 ,   77.94229) -- (    .00000 ,     .00000);
\node at (  13.00000 ,    5.19615) [circle    , fill=blue     ] {};
\draw [color=blue     ] (  13.00000 ,    5.19615) -- (  13.50000 ,    6.06218);
\draw [color=blue     ] (  13.50000 ,    6.06218) -- (  14.00000 ,    6.92820);
\draw [color=blue     ] (  14.00000 ,    6.92820) -- (  14.50000 ,    7.79423);
\draw [color=blue     ] (  14.50000 ,    7.79423) -- (  15.00000 ,    8.66025);
\draw [color=blue     ] (  15.00000 ,    8.66025) -- (  15.50000 ,    9.52628);
\draw [color=blue     ] (  15.50000 ,    9.52628) -- (  16.00000 ,   10.39230);
\draw [color=blue     ] (  16.00000 ,   10.39230) -- (  16.50000 ,   11.25833);
\draw [color=blue     ] (  16.50000 ,   11.25833) -- (  17.00000 ,   12.12436);
\draw [color=blue     ] (  17.00000 ,   12.12436) -- (  17.50000 ,   12.99038);
\draw [color=blue     ] (  17.50000 ,   12.99038) -- (  18.00000 ,   13.85641);
\draw [color=blue     ] (  18.00000 ,   13.85641) -- (  18.50000 ,   14.72243);
\draw [color=blue     ] (  18.50000 ,   14.72243) -- (  19.00000 ,   15.58846);
\node at (  19.50000 ,   16.45448) [diamond    , fill=orange   ] {};
\draw [color=blue     ] (  19.00000 ,   15.58846) -- (  19.50000 ,   16.45448);
\draw [color=blue     ] (  19.50000 ,   16.45448) -- (  20.00000 ,   17.32051);
\draw [color=blue     ] (  20.00000 ,   17.32051) -- (  20.50000 ,   18.18653);
\draw [color=blue     ] (  20.50000 ,   18.18653) -- (  21.00000 ,   19.05256);
\draw [color=blue     ] (  21.00000 ,   19.05256) -- (  21.50000 ,   19.91858);
\draw [color=blue     ] (  21.50000 ,   19.91858) -- (  22.00000 ,   20.78461);
\draw [color=blue     ] (  22.00000 ,   20.78461) -- (  22.50000 ,   21.65064);
\draw [color=blue     ] (  22.50000 ,   21.65064) -- (  23.00000 ,   22.51666);
\draw [color=blue     ] (  23.00000 ,   22.51666) -- (  23.50000 ,   23.38269);
\draw [color=blue     ] (  23.50000 ,   23.38269) -- (  24.00000 ,   24.24871);
\node at (  24.50000 ,   25.11474) [diamond    , fill=orange   ] {};
\draw [color=blue     ] (  24.00000 ,   24.24871) -- (  24.50000 ,   25.11474);
\draw [color=blue     ] (  24.50000 ,   25.11474) -- (  25.00000 ,   25.98076);
\draw [color=blue     ] (  25.00000 ,   25.98076) -- (  25.50000 ,   26.84679);
\draw [color=blue     ] (  25.50000 ,   26.84679) -- (  26.00000 ,   27.71281);
\draw [color=blue     ] (  26.00000 ,   27.71281) -- (  26.50000 ,   28.57884);
\draw [color=blue     ] (  26.50000 ,   28.57884) -- (  27.00000 ,   29.44486);
\draw [color=blue     ] (  27.00000 ,   29.44486) -- (  27.50000 ,   30.31089);
\draw [color=blue     ] (  27.50000 ,   30.31089) -- (  28.00000 ,   31.17691);
\draw [color=blue     ] (  28.00000 ,   31.17691) -- (  28.50000 ,   32.04294);
\draw [color=blue     ] (  28.50000 ,   32.04294) -- (  29.00000 ,   32.90897);
\draw [color=blue     ] (  29.00000 ,   32.90897) -- (  30.00000 ,   32.90897);
\draw [color=blue     ] (  30.00000 ,   32.90897) -- (  30.50000 ,   33.77499);
\draw [color=blue     ] (  30.50000 ,   33.77499) -- (  31.00000 ,   34.64102);
\draw [color=blue     ] (  31.00000 ,   34.64102) -- (  31.50000 ,   35.50704);
\draw [color=blue     ] (  31.50000 ,   35.50704) -- (  32.00000 ,   36.37307);
\draw [color=blue     ] (  32.00000 ,   36.37307) -- (  33.00000 ,   36.37307);
\draw [color=blue     ] (  33.00000 ,   36.37307) -- (  33.50000 ,   37.23909);
\draw [color=blue     ] (  33.50000 ,   37.23909) -- (  34.00000 ,   38.10512);
\draw [color=blue     ] (  34.00000 ,   38.10512) -- (  34.50000 ,   38.97114);
\draw [color=blue     ] (  34.50000 ,   38.97114) -- (  35.00000 ,   39.83717);
\draw [color=blue     ] (  35.00000 ,   39.83717) -- (  35.50000 ,   40.70319);
\draw [color=blue     ] (  35.50000 ,   40.70319) -- (  36.00000 ,   41.56922);
\draw [color=blue     ] (  36.00000 ,   41.56922) -- (  37.00000 ,   41.56922);
\draw [color=blue     ] (  37.00000 ,   41.56922) -- (  38.00000 ,   41.56922);
\draw [color=blue     ] (  38.00000 ,   41.56922) -- (  39.00000 ,   41.56922);
\draw [color=blue     ] (  39.00000 ,   41.56922) -- (  39.50000 ,   42.43524);
\draw [color=blue     ] (  39.50000 ,   42.43524) -- (  40.50000 ,   42.43524);
\draw [color=blue     ] (  40.50000 ,   42.43524) -- (  41.00000 ,   43.30127);
\draw [color=blue     ] (  41.00000 ,   43.30127) -- (  41.50000 ,   44.16730);
\draw [color=blue     ] (  41.50000 ,   44.16730) -- (  42.00000 ,   45.03332);
\draw [color=blue     ] (  42.00000 ,   45.03332) -- (  42.50000 ,   45.89935);
\draw [color=blue     ] (  42.50000 ,   45.89935) -- (  43.00000 ,   46.76537);
\draw [color=blue     ] (  43.00000 ,   46.76537) -- (  43.50000 ,   47.63140);
\draw [color=blue     ] (  43.50000 ,   47.63140) -- (  44.00000 ,   48.49742);
\draw [color=blue     ] (  44.00000 ,   48.49742) -- (  44.50000 ,   49.36345);
\draw [color=blue     ] (  44.50000 ,   49.36345) -- (  45.00000 ,   50.22947);
\draw [color=blue     ] (  45.00000 ,   50.22947) -- (  46.00000 ,   50.22947);
\draw [color=blue     ] (  46.00000 ,   50.22947) -- (  46.50000 ,   51.09550);
\draw [color=blue     ] (  46.50000 ,   51.09550) -- (  47.00000 ,   51.96152);
\draw [color=blue     ] (  47.00000 ,   51.96152) -- (  47.50000 ,   52.82755);
\draw [color=blue     ] (  47.50000 ,   52.82755) -- (  48.00000 ,   53.69358);
\draw [color=blue     ] (  48.00000 ,   53.69358) -- (  48.50000 ,   54.55960);
\draw [color=blue     ] (  48.50000 ,   54.55960) -- (  49.00000 ,   55.42563);
\draw [color=blue     ] (  49.00000 ,   55.42563) -- (  49.50000 ,   56.29165);
\draw [color=blue     ] (  49.50000 ,   56.29165) -- (  50.50000 ,   56.29165);
\draw [color=blue     ] (  50.50000 ,   56.29165) -- (  51.50000 ,   56.29165);
\draw [color=blue     ] (  51.50000 ,   56.29165) -- (  52.50000 ,   56.29165);
\draw [color=blue     ] (  52.50000 ,   56.29165) -- (  53.50000 ,   56.29165);
\draw [color=blue     ] (  53.50000 ,   56.29165) -- (  54.50000 ,   56.29165);
\draw [color=blue     ] (  54.50000 ,   56.29165) -- (  55.00000 ,   57.15768);
\draw [color=blue     ] (  55.00000 ,   57.15768) -- (  55.50000 ,   58.02370);
\draw [color=blue     ] (  55.50000 ,   58.02370) -- (  56.00000 ,   58.88973);
\draw [color=blue     ] (  56.00000 ,   58.88973) -- (  56.50000 ,   59.75575);
\draw [color=blue     ] (  56.50000 ,   59.75575) -- (  57.00000 ,   60.62178);
\draw [color=blue     ] (  57.00000 ,   60.62178) -- (  57.50000 ,   61.48780);
\draw [color=blue     ] (  57.50000 ,   61.48780) -- (  58.00000 ,   62.35383);
\draw [color=blue     ] (  58.00000 ,   62.35383) -- (  58.50000 ,   63.21985);
\draw [color=blue     ] (  58.50000 ,   63.21985) -- (  59.00000 ,   64.08588);
\draw [color=blue     ] (  59.00000 ,   64.08588) -- (  59.50000 ,   64.95191);
\draw [color=blue     ] (  59.50000 ,   64.95191) -- (  60.50000 ,   64.95191);
\draw [color=blue     ] (  60.50000 ,   64.95191) -- (  61.50000 ,   64.95191);
\draw [color=blue     ] (  61.50000 ,   64.95191) -- (  62.00000 ,   65.81793);
\draw [color=blue     ] (  62.00000 ,   65.81793) -- (  62.50000 ,   66.68396);
\draw [color=blue     ] (  62.50000 ,   66.68396) -- (  63.00000 ,   67.54998);
\draw [color=blue     ] (  63.00000 ,   67.54998) -- (  63.50000 ,   68.41601);
\draw [color=blue     ] (  63.50000 ,   68.41601) -- (  64.00000 ,   69.28203);
\draw [color=blue     ] (  64.00000 ,   69.28203) -- (  64.50000 ,   70.14806);
\draw [color=blue     ] (  64.50000 ,   70.14806) -- (  65.00000 ,   71.01408);
\draw [color=blue     ] (  65.00000 ,   71.01408) -- (  66.00000 ,   71.01408);
\draw [color=blue     ] (  66.00000 ,   71.01408) -- (  67.00000 ,   71.01408);
\draw [color=blue     ] (  67.00000 ,   71.01408) -- (  67.50000 ,   71.88011);
\node at (  67.50000 ,   71.88011) [rectangle    , fill=green    ] {};
\end{tikzpicture}
\end{scaletikzpicturetowidth}
\caption{Case 2: Time to target 93 time steps, agent knockouts at time steps 13 and 23, knocked out sensors not shown}
\label{fig12e}
\end{figure}

For case 2, when two knockouts are allowed, 
Figure~\ref{fig12d} shows the solution. The minimal time to target reduces to 93 time steps. Figure~\ref{fig12e} shows the same solution as Figure~\ref{fig12d}, but with the sensors that have been knocked out removed for clarity.

The computation time required to find, and prove to be optimal, the solution seen in 
Figure~\ref{fig12b} 
was over 2 minutes. The computation time required to find, and prove to be optimal,  the solution seen in
Figure~\ref{fig12d} was 7 seconds.

\subsubsection{Cases 3,4,5}

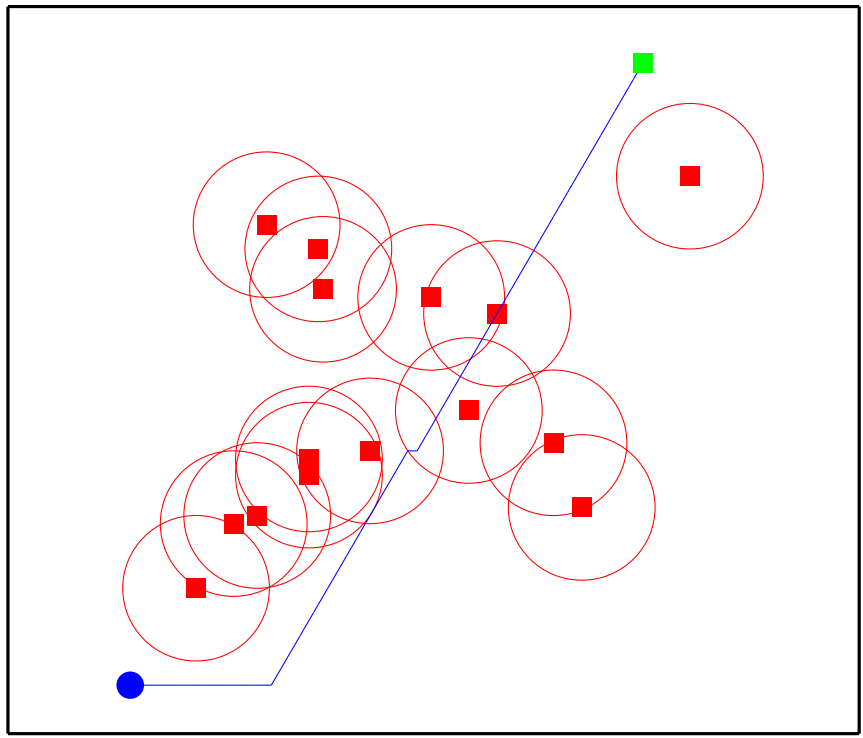
\begin{figure}[!htb]
\centering	  
\begin{scaletikzpicturetowidth}{\textwidth*0.70} 
\begin{tikzpicture}[scale=\tikzscale]	
\draw [red] (  26.50000 ,   23.38269) circle [radius=   7.80000];
\draw [red] (  32.00000 ,   29.44486) circle [radius=   7.80000];
\draw [red] (  72.50000 ,   59.75575) circle [radius=   7.80000];
\draw [red] (  20.00000 ,   15.58846) circle [radius=   7.80000];
\draw [red] (  24.00000 ,   22.51666) circle [radius=   7.80000];
\draw [red] (  33.00000 ,   51.96152) circle [radius=   7.80000];
\draw [red] (  61.00000 ,   24.24871) circle [radius=   7.80000];
\draw [red] (  33.50000 ,   47.63140) circle [radius=   7.80000];
\draw [red] (  27.50000 ,   54.55960) circle [radius=   7.80000];
\draw [red] (  45.00000 ,   46.76537) circle [radius=   7.80000];
\draw [red] (  32.00000 ,   27.71281) circle [radius=   7.80000];
\draw [red] (  52.00000 ,   45.03332) circle [radius=   7.80000];
\draw [red] (  38.50000 ,   30.31089) circle [radius=   7.80000];
\draw [red] (  58.00000 ,   31.17691) circle [radius=   7.80000];
\draw [red] (  49.00000 ,   34.64102) circle [radius=   7.80000];
\node at (  26.50000 ,   23.38269) [rectangle    , fill=red    ] {};
\node at (  32.00000 ,   29.44486) [rectangle    , fill=red    ] {};
\node at ( 72.50000 ,   59.75575) [rectangle    , fill=red    ] {};
\node at (  20.00000 ,   15.58846) [rectangle    , fill=red    ] {};
\node at ( 24.00000 ,   22.51666) [rectangle    , fill=red    ] {};
\node at (  33.00000 ,   51.96152) [rectangle    , fill=red    ] {};
\node at (  61.00000 ,   24.24871) [rectangle    , fill=red    ] {};
\node at ( 33.50000 ,   47.63140) [rectangle    , fill=red    ] {};
\node at (  27.50000 ,   54.55960) [rectangle    , fill=red    ] {};
\node at (  45.00000 ,   46.76537) [rectangle    , fill=red    ] {};
\node at (  32.00000 ,   27.71281) [rectangle    , fill=red    ] {};
\node at (  52.00000 ,   45.03332) [rectangle    , fill=red    ] {};
\node at (  38.50000 ,   30.31089) [rectangle    , fill=red    ] {};
\node at (  58.00000 ,   31.17691) [rectangle    , fill=red    ] {};
\node at (  49.00000 ,   34.64102) [rectangle    , fill=red    ] {};

\draw[very thick] (    .00000 ,     .00000) -- (  90.50000 ,     .00000);
\draw[very thick] (  90.50000 ,     .00000) -- (  90.50000 ,   77.94229);
\draw[very thick] (  90.50000 ,   77.94229) -- (    .00000 ,   77.94229);
\draw[very thick] (    .00000 ,   77.94229) -- (    .00000 ,     .00000);
\node at (  13.00000 ,    5.19615) [circle    , fill=blue     ] {};
\draw [color=blue     ] (  13.00000 ,    5.19615) -- (  14.00000 ,    5.19615);
\draw [color=blue     ] (  14.00000 ,    5.19615) -- (  15.00000 ,    5.19615);
\draw [color=blue     ] (  15.00000 ,    5.19615) -- (  16.00000 ,    5.19615);
\draw [color=blue     ] (  16.00000 ,    5.19615) -- (  17.00000 ,    5.19615);
\draw [color=blue     ] (  17.00000 ,    5.19615) -- (  18.00000 ,    5.19615);
\draw [color=blue     ] (  18.00000 ,    5.19615) -- (  19.00000 ,    5.19615);
\draw [color=blue     ] (  19.00000 ,    5.19615) -- (  20.00000 ,    5.19615);
\draw [color=blue     ] (  20.00000 ,    5.19615) -- (  21.00000 ,    5.19615);
\draw [color=blue     ] (  21.00000 ,    5.19615) -- (  22.00000 ,    5.19615);
\draw [color=blue     ] (  22.00000 ,    5.19615) -- (  23.00000 ,    5.19615);
\draw [color=blue     ] (  23.00000 ,    5.19615) -- (  24.00000 ,    5.19615);
\draw [color=blue     ] (  24.00000 ,    5.19615) -- (  25.00000 ,    5.19615);
\draw [color=blue     ] (  25.00000 ,    5.19615) -- (  26.00000 ,    5.19615);
\draw [color=blue     ] (  26.00000 ,    5.19615) -- (  27.00000 ,    5.19615);
\draw [color=blue     ] (  27.00000 ,    5.19615) -- (  28.00000 ,    5.19615);
\draw [color=blue     ] (  28.00000 ,    5.19615) -- (  28.50000 ,    6.06218);
\draw [color=blue     ] (  28.50000 ,    6.06218) -- (  29.00000 ,    6.92820);
\draw [color=blue     ] (  29.00000 ,    6.92820) -- (  29.50000 ,    7.79423);
\draw [color=blue     ] (  29.50000 ,    7.79423) -- (  30.00000 ,    8.66025);
\draw [color=blue     ] (  30.00000 ,    8.66025) -- (  30.50000 ,    9.52628);
\draw [color=blue     ] (  30.50000 ,    9.52628) -- (  31.00000 ,   10.39230);
\draw [color=blue     ] (  31.00000 ,   10.39230) -- (  31.50000 ,   11.25833);
\draw [color=blue     ] (  31.50000 ,   11.25833) -- (  32.00000 ,   12.12436);
\draw [color=blue     ] (  32.00000 ,   12.12436) -- (  32.50000 ,   12.99038);
\draw [color=blue     ] (  32.50000 ,   12.99038) -- (  33.00000 ,   13.85641);
\draw [color=blue     ] (  33.00000 ,   13.85641) -- (  33.50000 ,   14.72243);
\draw [color=blue     ] (  33.50000 ,   14.72243) -- (  34.00000 ,   15.58846);
\draw [color=blue     ] (  34.00000 ,   15.58846) -- (  34.50000 ,   16.45448);
\draw [color=blue     ] (  34.50000 ,   16.45448) -- (  35.00000 ,   17.32051);
\draw [color=blue     ] (  35.00000 ,   17.32051) -- (  35.50000 ,   18.18653);
\draw [color=blue     ] (  35.50000 ,   18.18653) -- (  36.00000 ,   19.05256);
\draw [color=blue     ] (  36.00000 ,   19.05256) -- (  36.50000 ,   19.91858);
\draw [color=blue     ] (  36.50000 ,   19.91858) -- (  37.00000 ,   20.78461);
\draw [color=blue     ] (  37.00000 ,   20.78461) -- (  37.50000 ,   21.65064);
\draw [color=blue     ] (  37.50000 ,   21.65064) -- (  38.00000 ,   22.51666);
\draw [color=blue     ] (  38.00000 ,   22.51666) -- (  38.50000 ,   23.38269);
\draw [color=blue     ] (  38.50000 ,   23.38269) -- (  39.00000 ,   24.24871);
\draw [color=blue     ] (  39.00000 ,   24.24871) -- (  39.50000 ,   25.11474);
\draw [color=blue     ] (  39.50000 ,   25.11474) -- (  40.00000 ,   25.98076);
\draw [color=blue     ] (  40.00000 ,   25.98076) -- (  40.50000 ,   26.84679);
\draw [color=blue     ] (  40.50000 ,   26.84679) -- (  41.00000 ,   27.71281);
\draw [color=blue     ] (  41.00000 ,   27.71281) -- (  41.50000 ,   28.57884);
\draw [color=blue     ] (  41.50000 ,   28.57884) -- (  42.00000 ,   29.44486);
\draw [color=blue     ] (  42.00000 ,   29.44486) -- (  42.50000 ,   30.31089);
\draw [color=blue     ] (  42.50000 ,   30.31089) -- (  43.50000 ,   30.31089);
\draw [color=blue     ] (  43.50000 ,   30.31089) -- (  44.00000 ,   31.17691);
\draw [color=blue     ] (  44.00000 ,   31.17691) -- (  44.50000 ,   32.04294);
\draw [color=blue     ] (  44.50000 ,   32.04294) -- (  45.00000 ,   32.90897);
\draw [color=blue     ] (  45.00000 ,   32.90897) -- (  45.50000 ,   33.77499);
\draw [color=blue     ] (  45.50000 ,   33.77499) -- (  46.00000 ,   34.64102);
\draw [color=blue     ] (  46.00000 ,   34.64102) -- (  46.50000 ,   35.50704);
\draw [color=blue     ] (  46.50000 ,   35.50704) -- (  47.00000 ,   36.37307);
\draw [color=blue     ] (  47.00000 ,   36.37307) -- (  47.50000 ,   37.23909);
\draw [color=blue     ] (  47.50000 ,   37.23909) -- (  48.00000 ,   38.10512);
\draw [color=blue     ] (  48.00000 ,   38.10512) -- (  48.50000 ,   38.97114);
\draw [color=blue     ] (  48.50000 ,   38.97114) -- (  49.00000 ,   39.83717);
\draw [color=blue     ] (  49.00000 ,   39.83717) -- (  49.50000 ,   40.70319);
\draw [color=blue     ] (  49.50000 ,   40.70319) -- (  50.00000 ,   41.56922);
\draw [color=blue     ] (  50.00000 ,   41.56922) -- (  50.50000 ,   42.43524);
\draw [color=blue     ] (  50.50000 ,   42.43524) -- (  51.00000 ,   43.30127);
\draw [color=blue     ] (  51.00000 ,   43.30127) -- (  51.50000 ,   44.16730);
\draw [color=blue     ] (  51.50000 ,   44.16730) -- (  52.00000 ,   45.03332);
\draw [color=blue     ] (  52.00000 ,   45.03332) -- (  52.50000 ,   45.89935);
\draw [color=blue     ] (  52.50000 ,   45.89935) -- (  53.00000 ,   46.76537);
\draw [color=blue     ] (  53.00000 ,   46.76537) -- (  53.50000 ,   47.63140);
\draw [color=blue     ] (  53.50000 ,   47.63140) -- (  54.00000 ,   48.49742);
\draw [color=blue     ] (  54.00000 ,   48.49742) -- (  54.50000 ,   49.36345);
\draw [color=blue     ] (  54.50000 ,   49.36345) -- (  55.00000 ,   50.22947);
\draw [color=blue     ] (  55.00000 ,   50.22947) -- (  55.50000 ,   51.09550);
\draw [color=blue     ] (  55.50000 ,   51.09550) -- (  56.00000 ,   51.96152);
\draw [color=blue     ] (  56.00000 ,   51.96152) -- (  56.50000 ,   52.82755);
\draw [color=blue     ] (  56.50000 ,   52.82755) -- (  57.00000 ,   53.69358);
\draw [color=blue     ] (  57.00000 ,   53.69358) -- (  57.50000 ,   54.55960);
\draw [color=blue     ] (  57.50000 ,   54.55960) -- (  58.00000 ,   55.42563);
\draw [color=blue     ] (  58.00000 ,   55.42563) -- (  58.50000 ,   56.29165);
\draw [color=blue     ] (  58.50000 ,   56.29165) -- (  59.00000 ,   57.15768);
\draw [color=blue     ] (  59.00000 ,   57.15768) -- (  59.50000 ,   58.02370);
\draw [color=blue     ] (  59.50000 ,   58.02370) -- (  60.00000 ,   58.88973);
\draw [color=blue     ] (  60.00000 ,   58.88973) -- (  60.50000 ,   59.75575);
\draw [color=blue     ] (  60.50000 ,   59.75575) -- (  61.00000 ,   60.62178);
\draw [color=blue     ] (  61.00000 ,   60.62178) -- (  61.50000 ,   61.48780);
\draw [color=blue     ] (  61.50000 ,   61.48780) -- (  62.00000 ,   62.35383);
\draw [color=blue     ] (  62.00000 ,   62.35383) -- (  62.50000 ,   63.21985);
\draw [color=blue     ] (  62.50000 ,   63.21985) -- (  63.00000 ,   64.08588);
\draw [color=blue     ] (  63.00000 ,   64.08588) -- (  63.50000 ,   64.95191);
\draw [color=blue     ] (  63.50000 ,   64.95191) -- (  64.00000 ,   65.81793);
\draw [color=blue     ] (  64.00000 ,   65.81793) -- (  64.50000 ,   66.68396);
\draw [color=blue     ] (  64.50000 ,   66.68396) -- (  65.00000 ,   67.54998);
\draw [color=blue     ] (  65.00000 ,   67.54998) -- (  65.50000 ,   68.41601);
\draw [color=blue     ] (  65.50000 ,   68.41601) -- (  66.00000 ,   69.28203);
\draw [color=blue     ] (  66.00000 ,   69.28203) -- (  66.50000 ,   70.14806);
\draw [color=blue     ] (  66.50000 ,   70.14806) -- (  67.00000 ,   71.01408);
\draw [color=blue     ] (  67.00000 ,   71.01408) -- (  67.50000 ,   71.88011);
\node at (  67.50000 ,   71.88011) [rectangle    , fill=green    ] {};
\end{tikzpicture}
\end{scaletikzpicturetowidth}
\caption{Case 3: Minimal time to target 93 time steps, no confusion, PED much less than 1\%}
\label{fig345b}
\end{figure}

Figure~\ref{fig345b} shows for the example problem considered the solution for case 3. Here the minimal time to target is 93 time steps, but with an associated probability of evading detection (PED) that is much less than 1\%.

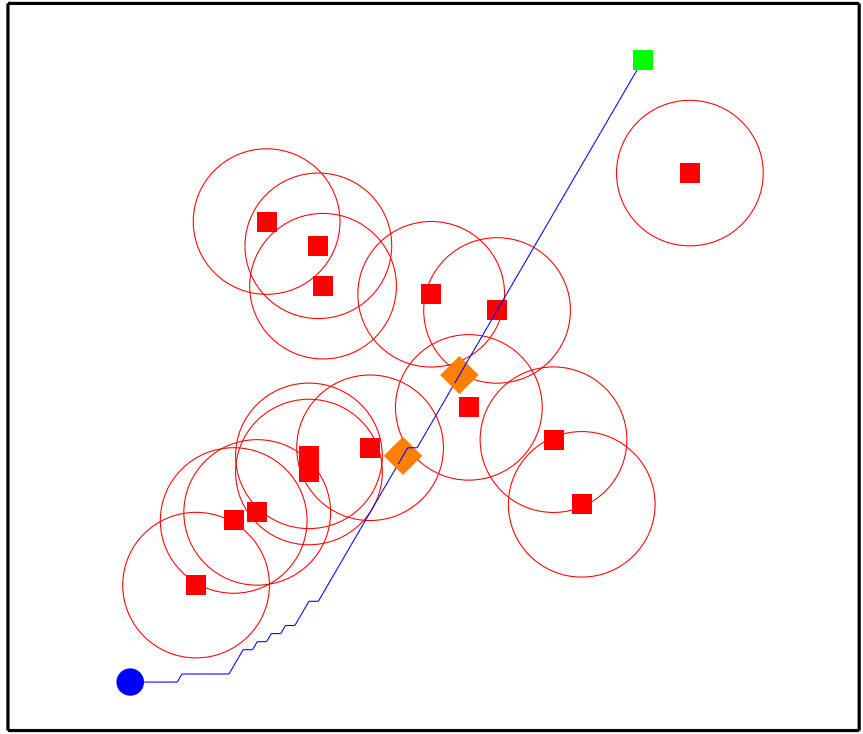
\begin{figure}[!htb]
\centering	  
\begin{scaletikzpicturetowidth}{\textwidth*0.70} 
\begin{tikzpicture}[scale=\tikzscale]	
\draw [red] (  26.50000 ,   23.38269) circle [radius=   7.80000];
\draw [red] (  32.00000 ,   29.44486) circle [radius=   7.80000];
\draw [red] (  72.50000 ,   59.75575) circle [radius=   7.80000];
\draw [red] (  20.00000 ,   15.58846) circle [radius=   7.80000];
\draw [red] (  24.00000 ,   22.51666) circle [radius=   7.80000];
\draw [red] (  33.00000 ,   51.96152) circle [radius=   7.80000];
\draw [red] (  61.00000 ,   24.24871) circle [radius=   7.80000];
\draw [red] (  33.50000 ,   47.63140) circle [radius=   7.80000];
\draw [red] (  27.50000 ,   54.55960) circle [radius=   7.80000];
\draw [red] (  45.00000 ,   46.76537) circle [radius=   7.80000];
\draw [red] (  32.00000 ,   27.71281) circle [radius=   7.80000];
\draw [red] (  52.00000 ,   45.03332) circle [radius=   7.80000];
\draw [red] (  38.50000 ,   30.31089) circle [radius=   7.80000];
\draw [red] (  58.00000 ,   31.17691) circle [radius=   7.80000];
\draw [red] (  49.00000 ,   34.64102) circle [radius=   7.80000];

\node at (  26.50000 ,   23.38269) [rectangle    , fill=red    ] {};
\node at (  32.00000 ,   29.44486) [rectangle    , fill=red    ] {};
\node at ( 72.50000 ,   59.75575) [rectangle    , fill=red    ] {};
\node at (  20.00000 ,   15.58846) [rectangle    , fill=red    ] {};
\node at ( 24.00000 ,   22.51666) [rectangle    , fill=red    ] {};
\node at (  33.00000 ,   51.96152) [rectangle    , fill=red    ] {};
\node at (  61.00000 ,   24.24871) [rectangle    , fill=red    ] {};
\node at ( 33.50000 ,   47.63140) [rectangle    , fill=red    ] {};
\node at (  27.50000 ,   54.55960) [rectangle    , fill=red    ] {};
\node at (  45.00000 ,   46.76537) [rectangle    , fill=red    ] {};
\node at (  32.00000 ,   27.71281) [rectangle    , fill=red    ] {};
\node at (  52.00000 ,   45.03332) [rectangle    , fill=red    ] {};
\node at (  38.50000 ,   30.31089) [rectangle    , fill=red    ] {};
\node at (  58.00000 ,   31.17691) [rectangle    , fill=red    ] {};
\node at (  49.00000 ,   34.64102) [rectangle    , fill=red    ] {};

\draw[very thick] (    .00000 ,     .00000) -- (  90.50000 ,     .00000);
\draw[very thick] (  90.50000 ,     .00000) -- (  90.50000 ,   77.94229);
\draw[very thick] (  90.50000 ,   77.94229) -- (    .00000 ,   77.94229);
\draw[very thick] (    .00000 ,   77.94229) -- (    .00000 ,     .00000);
\node at (  13.00000 ,    5.19615) [circle    , fill=blue     ] {};
\draw [color=blue     ] (  13.00000 ,    5.19615) -- (  14.00000 ,    5.19615);
\draw [color=blue     ] (  14.00000 ,    5.19615) -- (  15.00000 ,    5.19615);
\draw [color=blue     ] (  15.00000 ,    5.19615) -- (  16.00000 ,    5.19615);
\draw [color=blue     ] (  16.00000 ,    5.19615) -- (  17.00000 ,    5.19615);
\draw [color=blue     ] (  17.00000 ,    5.19615) -- (  18.00000 ,    5.19615);
\draw [color=blue     ] (  18.00000 ,    5.19615) -- (  18.50000 ,    6.06218);
\draw [color=blue     ] (  18.50000 ,    6.06218) -- (  19.50000 ,    6.06218);
\draw [color=blue     ] (  19.50000 ,    6.06218) -- (  20.50000 ,    6.06218);
\draw [color=blue     ] (  20.50000 ,    6.06218) -- (  21.50000 ,    6.06218);
\draw [color=blue     ] (  21.50000 ,    6.06218) -- (  22.50000 ,    6.06218);
\draw [color=blue     ] (  22.50000 ,    6.06218) -- (  23.50000 ,    6.06218);
\draw [color=blue     ] (  23.50000 ,    6.06218) -- (  24.00000 ,    6.92820);
\draw [color=blue     ] (  24.00000 ,    6.92820) -- (  24.50000 ,    7.79423);
\draw [color=blue     ] (  24.50000 ,    7.79423) -- (  25.00000 ,    8.66025);
\draw [color=blue     ] (  25.00000 ,    8.66025) -- (  26.00000 ,    8.66025);
\draw [color=blue     ] (  26.00000 ,    8.66025) -- (  26.50000 ,    9.52628);
\draw [color=blue     ] (  26.50000 ,    9.52628) -- (  27.50000 ,    9.52628);
\draw [color=blue     ] (  27.50000 ,    9.52628) -- (  28.00000 ,   10.39230);
\draw [color=blue     ] (  28.00000 ,   10.39230) -- (  29.00000 ,   10.39230);
\draw [color=blue     ] (  29.00000 ,   10.39230) -- (  29.50000 ,   11.25833);
\draw [color=blue     ] (  29.50000 ,   11.25833) -- (  30.50000 ,   11.25833);
\draw [color=blue     ] (  30.50000 ,   11.25833) -- (  31.00000 ,   12.12436);
\draw [color=blue     ] (  31.00000 ,   12.12436) -- (  31.50000 ,   12.99038);
\draw [color=blue     ] (  31.50000 ,   12.99038) -- (  32.00000 ,   13.85641);
\draw [color=blue     ] (  32.00000 ,   13.85641) -- (  33.00000 ,   13.85641);
\draw [color=blue     ] (  33.00000 ,   13.85641) -- (  33.50000 ,   14.72243);
\draw [color=blue     ] (  33.50000 ,   14.72243) -- (  34.00000 ,   15.58846);
\draw [color=blue     ] (  34.00000 ,   15.58846) -- (  34.50000 ,   16.45448);
\draw [color=blue     ] (  34.50000 ,   16.45448) -- (  35.00000 ,   17.32051);
\draw [color=blue     ] (  35.00000 ,   17.32051) -- (  35.50000 ,   18.18653);
\draw [color=blue     ] (  35.50000 ,   18.18653) -- (  36.00000 ,   19.05256);
\draw [color=blue     ] (  36.00000 ,   19.05256) -- (  36.50000 ,   19.91858);
\draw [color=blue     ] (  36.50000 ,   19.91858) -- (  37.00000 ,   20.78461);
\draw [color=blue     ] (  37.00000 ,   20.78461) -- (  37.50000 ,   21.65064);
\draw [color=blue     ] (  37.50000 ,   21.65064) -- (  38.00000 ,   22.51666);
\draw [color=blue     ] (  38.00000 ,   22.51666) -- (  38.50000 ,   23.38269);
\draw [color=blue     ] (  38.50000 ,   23.38269) -- (  39.00000 ,   24.24871);
\draw [color=blue     ] (  39.00000 ,   24.24871) -- (  39.50000 ,   25.11474);
\draw [color=blue     ] (  39.50000 ,   25.11474) -- (  40.00000 ,   25.98076);
\draw [color=blue     ] (  40.00000 ,   25.98076) -- (  40.50000 ,   26.84679);
\draw [color=blue     ] (  40.50000 ,   26.84679) -- (  41.00000 ,   27.71281);
\draw [color=blue     ] (  41.00000 ,   27.71281) -- (  41.50000 ,   28.57884);
\node at (  42.00000 ,   29.44486) [diamond    , fill=orange   ] {};
\draw [color=blue     ] (  41.50000 ,   28.57884) -- (  42.00000 ,   29.44486);
\draw [color=blue     ] (  42.00000 ,   29.44486) -- (  42.50000 ,   30.31089);
\draw [color=blue     ] (  42.50000 ,   30.31089) -- (  43.50000 ,   30.31089);
\draw [color=blue     ] (  43.50000 ,   30.31089) -- (  44.00000 ,   31.17691);
\draw [color=blue     ] (  44.00000 ,   31.17691) -- (  44.50000 ,   32.04294);
\draw [color=blue     ] (  44.50000 ,   32.04294) -- (  45.00000 ,   32.90897);
\draw [color=blue     ] (  45.00000 ,   32.90897) -- (  45.50000 ,   33.77499);
\draw [color=blue     ] (  45.50000 ,   33.77499) -- (  46.00000 ,   34.64102);
\draw [color=blue     ] (  46.00000 ,   34.64102) -- (  46.50000 ,   35.50704);
\draw [color=blue     ] (  46.50000 ,   35.50704) -- (  47.00000 ,   36.37307);
\draw [color=blue     ] (  47.00000 ,   36.37307) -- (  47.50000 ,   37.23909);
\node at (  48.00000 ,   38.10512) [diamond    , fill=orange   ] {};
\draw [color=blue     ] (  47.50000 ,   37.23909) -- (  48.00000 ,   38.10512);
\draw [color=blue     ] (  48.00000 ,   38.10512) -- (  48.50000 ,   38.97114);
\draw [color=blue     ] (  48.50000 ,   38.97114) -- (  49.00000 ,   39.83717);
\draw [color=blue     ] (  49.00000 ,   39.83717) -- (  49.50000 ,   40.70319);
\draw [color=blue     ] (  49.50000 ,   40.70319) -- (  50.00000 ,   41.56922);
\draw [color=blue     ] (  50.00000 ,   41.56922) -- (  50.50000 ,   42.43524);
\draw [color=blue     ] (  50.50000 ,   42.43524) -- (  51.00000 ,   43.30127);
\draw [color=blue     ] (  51.00000 ,   43.30127) -- (  51.50000 ,   44.16730);
\draw [color=blue     ] (  51.50000 ,   44.16730) -- (  52.00000 ,   45.03332);
\draw [color=blue     ] (  52.00000 ,   45.03332) -- (  52.50000 ,   45.89935);
\draw [color=blue     ] (  52.50000 ,   45.89935) -- (  53.00000 ,   46.76537);
\draw [color=blue     ] (  53.00000 ,   46.76537) -- (  53.50000 ,   47.63140);
\draw [color=blue     ] (  53.50000 ,   47.63140) -- (  54.00000 ,   48.49742);
\draw [color=blue     ] (  54.00000 ,   48.49742) -- (  54.50000 ,   49.36345);
\draw [color=blue     ] (  54.50000 ,   49.36345) -- (  55.00000 ,   50.22947);
\draw [color=blue     ] (  55.00000 ,   50.22947) -- (  55.50000 ,   51.09550);
\draw [color=blue     ] (  55.50000 ,   51.09550) -- (  56.00000 ,   51.96152);
\draw [color=blue     ] (  56.00000 ,   51.96152) -- (  56.50000 ,   52.82755);
\draw [color=blue     ] (  56.50000 ,   52.82755) -- (  57.00000 ,   53.69358);
\draw [color=blue     ] (  57.00000 ,   53.69358) -- (  57.50000 ,   54.55960);
\draw [color=blue     ] (  57.50000 ,   54.55960) -- (  58.00000 ,   55.42563);
\draw [color=blue     ] (  58.00000 ,   55.42563) -- (  58.50000 ,   56.29165);
\draw [color=blue     ] (  58.50000 ,   56.29165) -- (  59.00000 ,   57.15768);
\draw [color=blue     ] (  59.00000 ,   57.15768) -- (  59.50000 ,   58.02370);
\draw [color=blue     ] (  59.50000 ,   58.02370) -- (  60.00000 ,   58.88973);
\draw [color=blue     ] (  60.00000 ,   58.88973) -- (  60.50000 ,   59.75575);
\draw [color=blue     ] (  60.50000 ,   59.75575) -- (  61.00000 ,   60.62178);
\draw [color=blue     ] (  61.00000 ,   60.62178) -- (  61.50000 ,   61.48780);
\draw [color=blue     ] (  61.50000 ,   61.48780) -- (  62.00000 ,   62.35383);
\draw [color=blue     ] (  62.00000 ,   62.35383) -- (  62.50000 ,   63.21985);
\draw [color=blue     ] (  62.50000 ,   63.21985) -- (  63.00000 ,   64.08588);
\draw [color=blue     ] (  63.00000 ,   64.08588) -- (  63.50000 ,   64.95191);
\draw [color=blue     ] (  63.50000 ,   64.95191) -- (  64.00000 ,   65.81793);
\draw [color=blue     ] (  64.00000 ,   65.81793) -- (  64.50000 ,   66.68396);
\draw [color=blue     ] (  64.50000 ,   66.68396) -- (  65.00000 ,   67.54998);
\draw [color=blue     ] (  65.00000 ,   67.54998) -- (  65.50000 ,   68.41601);
\draw [color=blue     ] (  65.50000 ,   68.41601) -- (  66.00000 ,   69.28203);
\draw [color=blue     ] (  66.00000 ,   69.28203) -- (  66.50000 ,   70.14806);
\draw [color=blue     ] (  66.50000 ,   70.14806) -- (  67.00000 ,   71.01408);
\draw [color=blue     ] (  67.00000 ,   71.01408) -- (  67.50000 ,   71.88011);
\node at (  67.50000 ,   71.88011) [rectangle    , fill=green    ] {};
\end{tikzpicture}
\end{scaletikzpicturetowidth}
\caption{Case 4: Time to target 93 time steps, confusion actions at time steps 43 and 54, PED $=0.93$}
\label{fig345c}
\end{figure}

For the problem considered and case 4, so allowing two confusion actions, the solution is shown in 
Figure~\ref{fig345c}. The time to target is constrained to be the same as case 3,  but the two confusion actions, shown in orange, increase PED to 93\%.

The total computation time required to find, and prove to be optimal, the solution shown in  Figure~\ref{fig345b} was 8 seconds. The corresponding solution time for  
 Figure~\ref{fig345c} was 11 seconds.

For the problem considered and case 5, so allowing only one confusion action, but requiring a PED of at least 95\%, the solution is shown in Figure~\ref{fig345d}. The time to target is no longer 93, but increases to 97. 
The computation time required to find, and prove to be optimal, the solution shown in  Figure~\ref{fig345d}  was 30 seconds.

\subsection{More challenging problems}

The cases considered in Table~\ref{table1} involve an agent reaching the target in minimal time, subject to certain restrictions. It is therefore of little surprise that, graphically, the figures shown for example solutions are similar in the sense that they involve almost straight paths to target. Note here however that for cases 3 and 4 the time to target is set using  $T=\mbox{min}[D_{\gamma(a)N},~a=1,\ldots,A]$, so that it is equal to the minimal time to target ignoring sensors. Cases 1, 4 and 5 are cases where the time to target is not set, rather it emerges as a result of the optimisation.

Clearly varying mesh size,  as in Table~\ref{table1}, has a significant influence on computation time, both for the optimal and heuristic approaches. To produce more challenging problems we took the  five test problem instances on a   $90 \times 90$ mesh with two agents and 15 sensors as considered in
Table~\ref{table1},
 but  doubled the sensor detection area. Recall here that all our cases involve $\Omega=2$, so paths seek to avoid, as much as possible, vertices covered by two or more sensors.
Doubling the sensor detection area creates more overlaps between sensors and hence makes the problem more challenging. 

The results for these test problems are shown in Table~\ref{table2}. This table has the same format as 
Table~\ref{table1}. For these more challenging problems we imposed a time limit of 300 seconds, 5 minutes. For case 1 in Table~\ref{table2} all five instances with the optimal approach went to time limit and it is clear that for this case the optimal approach is not appropriate. For case 2 the optimal approach solved all instances within the time limit. It required more computation time than the heuristic approach, but that approach missed the optimal solution in two of the five instances.
For cases 3--5 the optimal approach clearly dominates the heuristic approach. The optimal approach has a lower average computation time than the heuristic approach and solves all instances to optimality. 

For illustration Figure~\ref{figsqrt} shows the solution associated with case 5 for the same problem instance as previously illustrated in Figure~\ref{fig345d}, but with doubled senor detection area. It can be seen that the time to target increases from the 97 time steps associated with Figure~\ref{fig345d} to 108 time steps, 
 with the agent that is used also changing.
The computation time required to find, and prove to be optimal, the solution shown in  Figure~\ref{figsqrt}  was 78 seconds.
Over cases 1, 2 and 5 in  Table~\ref{table2} the time to target, as found by the optimal or heuristic approach, increased by 7.4\% on average as compared with Table~\ref{table1}.

\begin{table}[!htb]
\centering
\renewcommand{\tabcolsep}{1mm} 
\renewcommand{\arraystretch}{1.3} 
\begin{tabular}{|c|cc|cc|}
\hline
\multirow{2}*{Case} &   \multicolumn{2}{c|}{Optimal} & \multicolumn{2}{c|} {Heuristic}    \\
& Time & \%  & Time & Not opt \\

\hline
1 &

Time limit	& 	98\%	& 	15.68	&	Not known

\\  

\hline
2 & 
247.11	&	97\%	& 	75.50	& 	2

\\ 

\hline
3  &

7.84	&	12\%	&	13.88	&	1
\\

\hline
4 & 
50.87	&	86\%	&	Time limit	& Not known
\\

\hline
5 & 
63.76	&	89\%	&	175.64	&	2
\\
 
\hline
\multicolumn{1}{|l|}{Percentage of variables excluded by Equations~(\ref{reduction1}) and (\ref{reduction2}) } &
 \multicolumn{2}{c|}{97.35} &  &  
\\
\hline
\end{tabular}
\caption{Computational results: more challenging problems}
\label{table2}
\end{table}

\begin{figure}[!htb]
\centering	  
\begin{scaletikzpicturetowidth}{\textwidth*0.70} 
\begin{tikzpicture}[scale=\tikzscale]	
\draw [red] (  26.50000 ,   23.38269) circle [radius=   7.80000];
\draw [red] (  32.00000 ,   29.44486) circle [radius=   7.80000];
\draw [red] (  72.50000 ,   59.75575) circle [radius=   7.80000];
\draw [red] (  20.00000 ,   15.58846) circle [radius=   7.80000];
\draw [red] (  24.00000 ,   22.51666) circle [radius=   7.80000];
\draw [red] (  33.00000 ,   51.96152) circle [radius=   7.80000];
\draw [red] (  61.00000 ,   24.24871) circle [radius=   7.80000];
\draw [red] (  33.50000 ,   47.63140) circle [radius=   7.80000];
\draw [red] (  27.50000 ,   54.55960) circle [radius=   7.80000];
\draw [red] (  45.00000 ,   46.76537) circle [radius=   7.80000];
\draw [red] (  32.00000 ,   27.71281) circle [radius=   7.80000];
\draw [red] (  52.00000 ,   45.03332) circle [radius=   7.80000];
\draw [red] (  38.50000 ,   30.31089) circle [radius=   7.80000];
\draw [red] (  58.00000 ,   31.17691) circle [radius=   7.80000];
\draw [red] (  49.00000 ,   34.64102) circle [radius=   7.80000];

\node at (  26.50000 ,   23.38269) [rectangle    , fill=red    ] {};
\node at (  32.00000 ,   29.44486) [rectangle    , fill=red    ] {};
\node at ( 72.50000 ,   59.75575) [rectangle    , fill=red    ] {};
\node at (  20.00000 ,   15.58846) [rectangle    , fill=red    ] {};
\node at ( 24.00000 ,   22.51666) [rectangle    , fill=red    ] {};
\node at (  33.00000 ,   51.96152) [rectangle    , fill=red    ] {};
\node at (  61.00000 ,   24.24871) [rectangle    , fill=red    ] {};
\node at ( 33.50000 ,   47.63140) [rectangle    , fill=red    ] {};
\node at (  27.50000 ,   54.55960) [rectangle    , fill=red    ] {};
\node at (  45.00000 ,   46.76537) [rectangle    , fill=red    ] {};
\node at (  32.00000 ,   27.71281) [rectangle    , fill=red    ] {};
\node at (  52.00000 ,   45.03332) [rectangle    , fill=red    ] {};
\node at (  38.50000 ,   30.31089) [rectangle    , fill=red    ] {};
\node at (  58.00000 ,   31.17691) [rectangle    , fill=red    ] {};
\node at (  49.00000 ,   34.64102) [rectangle    , fill=red    ] {};

\draw[very thick] (    .00000 ,     .00000) -- (  90.50000 ,     .00000);
\draw[very thick] (  90.50000 ,     .00000) -- (  90.50000 ,   77.94229);
\draw[very thick] (  90.50000 ,   77.94229) -- (    .00000 ,   77.94229);
\draw[very thick] (    .00000 ,   77.94229) -- (    .00000 ,     .00000);
\node at (  13.00000 ,    5.19615) [circle    , fill=blue     ] {};
\draw [color=blue     ] (  13.00000 ,    5.19615) -- (  13.50000 ,    6.06218);
\draw [color=blue     ] (  13.50000 ,    6.06218) -- (  14.00000 ,    6.92820);
\draw [color=blue     ] (  14.00000 ,    6.92820) -- (  14.50000 ,    7.79423);
\draw [color=blue     ] (  14.50000 ,    7.79423) -- (  15.50000 ,    7.79423);
\draw [color=blue     ] (  15.50000 ,    7.79423) -- (  16.00000 ,    8.66025);
\draw [color=blue     ] (  16.00000 ,    8.66025) -- (  16.50000 ,    9.52628);
\draw [color=blue     ] (  16.50000 ,    9.52628) -- (  17.50000 ,    9.52628);
\draw [color=blue     ] (  17.50000 ,    9.52628) -- (  18.50000 ,    9.52628);
\draw [color=blue     ] (  18.50000 ,    9.52628) -- (  19.50000 ,    9.52628);
\draw [color=blue     ] (  19.50000 ,    9.52628) -- (  20.50000 ,    9.52628);
\draw [color=blue     ] (  20.50000 ,    9.52628) -- (  21.00000 ,   10.39230);
\draw [color=blue     ] (  21.00000 ,   10.39230) -- (  21.50000 ,   11.25833);
\draw [color=blue     ] (  21.50000 ,   11.25833) -- (  22.00000 ,   12.12436);
\draw [color=blue     ] (  22.00000 ,   12.12436) -- (  22.50000 ,   12.99038);
\draw [color=blue     ] (  22.50000 ,   12.99038) -- (  23.50000 ,   12.99038);
\draw [color=blue     ] (  23.50000 ,   12.99038) -- (  24.50000 ,   12.99038);
\draw [color=blue     ] (  24.50000 ,   12.99038) -- (  25.50000 ,   12.99038);
\draw [color=blue     ] (  25.50000 ,   12.99038) -- (  26.50000 ,   12.99038);
\draw [color=blue     ] (  26.50000 ,   12.99038) -- (  27.00000 ,   13.85641);
\draw [color=blue     ] (  27.00000 ,   13.85641) -- (  27.50000 ,   14.72243);
\draw [color=blue     ] (  27.50000 ,   14.72243) -- (  28.00000 ,   15.58846);
\draw [color=blue     ] (  28.00000 ,   15.58846) -- (  29.00000 ,   15.58846);
\draw [color=blue     ] (  29.00000 ,   15.58846) -- (  29.50000 ,   16.45448);
\draw [color=blue     ] (  29.50000 ,   16.45448) -- (  30.00000 ,   17.32051);
\draw [color=blue     ] (  30.00000 ,   17.32051) -- (  31.00000 ,   17.32051);
\draw [color=blue     ] (  31.00000 ,   17.32051) -- (  32.00000 ,   17.32051);
\draw [color=blue     ] (  32.00000 ,   17.32051) -- (  33.00000 ,   17.32051);
\draw [color=blue     ] (  33.00000 ,   17.32051) -- (  34.00000 ,   17.32051);
\draw [color=blue     ] (  34.00000 ,   17.32051) -- (  35.00000 ,   17.32051);
\draw [color=blue     ] (  35.00000 ,   17.32051) -- (  36.00000 ,   17.32051);
\draw [color=blue     ] (  36.00000 ,   17.32051) -- (  36.50000 ,   18.18653);
\draw [color=blue     ] (  36.50000 ,   18.18653) -- (  37.50000 ,   18.18653);
\draw [color=blue     ] (  37.50000 ,   18.18653) -- (  38.50000 ,   18.18653);
\draw [color=blue     ] (  38.50000 ,   18.18653) -- (  39.00000 ,   19.05256);
\draw [color=blue     ] (  39.00000 ,   19.05256) -- (  39.50000 ,   19.91858);
\draw [color=blue     ] (  39.50000 ,   19.91858) -- (  40.50000 ,   19.91858);
\draw [color=blue     ] (  40.50000 ,   19.91858) -- (  41.00000 ,   20.78461);
\draw [color=blue     ] (  41.00000 ,   20.78461) -- (  41.50000 ,   21.65064);
\draw [color=blue     ] (  41.50000 ,   21.65064) -- (  42.50000 ,   21.65064);
\draw [color=blue     ] (  42.50000 ,   21.65064) -- (  43.00000 ,   22.51666);
\draw [color=blue     ] (  43.00000 ,   22.51666) -- (  43.50000 ,   23.38269);
\draw [color=blue     ] (  43.50000 ,   23.38269) -- (  44.00000 ,   24.24871);
\draw [color=blue     ] (  44.00000 ,   24.24871) -- (  44.50000 ,   25.11474);
\draw [color=blue     ] (  44.50000 ,   25.11474) -- (  45.00000 ,   25.98076);
\draw [color=blue     ] (  45.00000 ,   25.98076) -- (  45.50000 ,   26.84679);
\draw [color=blue     ] (  45.50000 ,   26.84679) -- (  46.00000 ,   27.71281);
\draw [color=blue     ] (  46.00000 ,   27.71281) -- (  46.50000 ,   28.57884);
\draw [color=blue     ] (  46.50000 ,   28.57884) -- (  47.00000 ,   29.44486);
\draw [color=blue     ] (  47.00000 ,   29.44486) -- (  47.50000 ,   30.31089);
\draw [color=blue     ] (  47.50000 ,   30.31089) -- (  48.00000 ,   31.17691);
\draw [color=blue     ] (  48.00000 ,   31.17691) -- (  48.50000 ,   32.04294);
\draw [color=blue     ] (  48.50000 ,   32.04294) -- (  49.00000 ,   32.90897);
\draw [color=blue     ] (  49.00000 ,   32.90897) -- (  49.50000 ,   33.77499);
\draw [color=blue     ] (  49.50000 ,   33.77499) -- (  50.00000 ,   34.64102);
\draw [color=blue     ] (  50.00000 ,   34.64102) -- (  50.50000 ,   35.50704);
\draw [color=blue     ] (  50.50000 ,   35.50704) -- (  51.00000 ,   36.37307);
\draw [color=blue     ] (  51.00000 ,   36.37307) -- (  50.50000 ,   37.23909);
\node at (  51.00000 ,   38.10512) [diamond    , fill=orange   ] {};
\draw [color=blue     ] (  50.50000 ,   37.23909) -- (  51.00000 ,   38.10512);
\draw [color=blue     ] (  51.00000 ,   38.10512) -- (  51.50000 ,   38.97114);
\draw [color=blue     ] (  51.50000 ,   38.97114) -- (  51.00000 ,   39.83717);
\draw [color=blue     ] (  51.00000 ,   39.83717) -- (  51.50000 ,   40.70319);
\draw [color=blue     ] (  51.50000 ,   40.70319) -- (  52.00000 ,   41.56922);
\draw [color=blue     ] (  52.00000 ,   41.56922) -- (  51.50000 ,   42.43524);
\draw [color=blue     ] (  51.50000 ,   42.43524) -- (  52.00000 ,   43.30127);
\draw [color=blue     ] (  52.00000 ,   43.30127) -- (  52.50000 ,   44.16730);
\draw [color=blue     ] (  52.50000 ,   44.16730) -- (  53.00000 ,   45.03332);
\draw [color=blue     ] (  53.00000 ,   45.03332) -- (  53.50000 ,   45.89935);
\draw [color=blue     ] (  53.50000 ,   45.89935) -- (  54.00000 ,   46.76537);
\draw [color=blue     ] (  54.00000 ,   46.76537) -- (  53.50000 ,   47.63140);
\draw [color=blue     ] (  53.50000 ,   47.63140) -- (  54.00000 ,   48.49742);
\draw [color=blue     ] (  54.00000 ,   48.49742) -- (  54.50000 ,   49.36345);
\draw [color=blue     ] (  54.50000 ,   49.36345) -- (  55.00000 ,   50.22947);
\draw [color=blue     ] (  55.00000 ,   50.22947) -- (  55.50000 ,   51.09550);
\draw [color=blue     ] (  55.50000 ,   51.09550) -- (  56.00000 ,   51.96152);
\draw [color=blue     ] (  56.00000 ,   51.96152) -- (  56.50000 ,   52.82755);
\draw [color=blue     ] (  56.50000 ,   52.82755) -- (  57.00000 ,   53.69358);
\draw [color=blue     ] (  57.00000 ,   53.69358) -- (  57.50000 ,   54.55960);
\draw [color=blue     ] (  57.50000 ,   54.55960) -- (  58.00000 ,   55.42563);
\draw [color=blue     ] (  58.00000 ,   55.42563) -- (  58.50000 ,   56.29165);
\draw [color=blue     ] (  58.50000 ,   56.29165) -- (  59.00000 ,   57.15768);
\draw [color=blue     ] (  59.00000 ,   57.15768) -- (  59.50000 ,   58.02370);
\draw [color=blue     ] (  59.50000 ,   58.02370) -- (  60.00000 ,   58.88973);
\draw [color=blue     ] (  60.00000 ,   58.88973) -- (  60.50000 ,   59.75575);
\draw [color=blue     ] (  60.50000 ,   59.75575) -- (  61.00000 ,   60.62178);
\draw [color=blue     ] (  61.00000 ,   60.62178) -- (  61.50000 ,   61.48780);
\draw [color=blue     ] (  61.50000 ,   61.48780) -- (  62.00000 ,   62.35383);
\draw [color=blue     ] (  62.00000 ,   62.35383) -- (  62.50000 ,   63.21985);
\draw [color=blue     ] (  62.50000 ,   63.21985) -- (  63.00000 ,   64.08588);
\draw [color=blue     ] (  63.00000 ,   64.08588) -- (  63.50000 ,   64.95191);
\draw [color=blue     ] (  63.50000 ,   64.95191) -- (  64.00000 ,   65.81793);
\draw [color=blue     ] (  64.00000 ,   65.81793) -- (  64.50000 ,   66.68396);
\draw [color=blue     ] (  64.50000 ,   66.68396) -- (  65.00000 ,   67.54998);
\draw [color=blue     ] (  65.00000 ,   67.54998) -- (  65.50000 ,   68.41601);
\draw [color=blue     ] (  65.50000 ,   68.41601) -- (  66.00000 ,   69.28203);
\draw [color=blue     ] (  66.00000 ,   69.28203) -- (  66.50000 ,   70.14806);
\draw [color=blue     ] (  66.50000 ,   70.14806) -- (  67.00000 ,   71.01408);
\draw [color=blue     ] (  67.00000 ,   71.01408) -- (  67.50000 ,   71.88011);
\node at (  67.50000 ,   71.88011) [rectangle    , fill=green    ] {};
\end{tikzpicture}
\end{scaletikzpicturetowidth}
\caption{Case 5: Time to target 97 time steps, confusion action at time step 58, PED $\geq 0.95$}
\label{fig345d}
\end{figure}
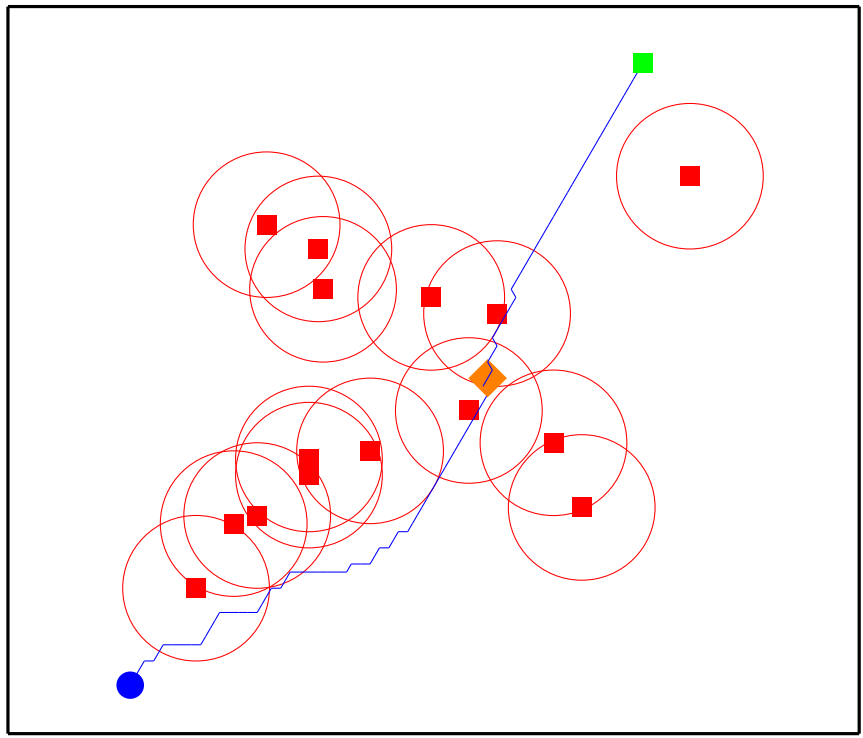

\begin{figure}[!htb]
\centering	  
\begin{scaletikzpicturetowidth}{\textwidth*0.70} 
\begin{tikzpicture}[scale=\tikzscale]	
\draw [red] (  26.50000 ,   23.38269) circle [radius=  11.03087];
\draw [red] (  32.00000 ,   29.44486) circle [radius=  11.03087];
\draw [red] (  72.50000 ,   59.75575) circle [radius=  11.03087];
\draw [red] (  20.00000 ,   15.58846) circle [radius=  11.03087];
\draw [red] (  24.00000 ,   22.51666) circle [radius=  11.03087];
\draw [red] (  33.00000 ,   51.96152) circle [radius=  11.03087];
\draw [red] (  61.00000 ,   24.24871) circle [radius=  11.03087];
\draw [red] (  33.50000 ,   47.63140) circle [radius=  11.03087];
\draw [red] (  27.50000 ,   54.55960) circle [radius=  11.03087];
\draw [red] (  45.00000 ,   46.76537) circle [radius=  11.03087];
\draw [red] (  32.00000 ,   27.71281) circle [radius=  11.03087];
\draw [red] (  52.00000 ,   45.03332) circle [radius=  11.03087];
\draw [red] (  38.50000 ,   30.31089) circle [radius=  11.03087];
\draw [red] (  58.00000 ,   31.17691) circle [radius=  11.03087];
\draw [red] (  49.00000 ,   34.64102) circle [radius=  11.03087];
\draw[very thick] (    .00000 ,     .00000) -- (  90.50000 ,     .00000);
\draw[very thick] (  90.50000 ,     .00000) -- (  90.50000 ,   77.94229);
\draw[very thick] (  90.50000 ,   77.94229) -- (    .00000 ,   77.94229);
\draw[very thick] (    .00000 ,   77.94229) -- (    .00000 ,     .00000);
\node at (   2.50000 ,    7.79423) [circle    , fill=blue     ] {};
\draw [color=blue     ] (   2.50000 ,    7.79423) -- (   3.00000 ,    8.66025);
\draw [color=blue     ] (   3.00000 ,    8.66025) -- (   3.50000 ,    9.52628);
\draw [color=blue     ] (   3.50000 ,    9.52628) -- (   4.00000 ,   10.39230);
\draw [color=blue     ] (   4.00000 ,   10.39230) -- (   3.50000 ,   11.25833);
\draw [color=blue     ] (   3.50000 ,   11.25833) -- (   4.00000 ,   12.12436);
\draw [color=blue     ] (   4.00000 ,   12.12436) -- (   4.50000 ,   12.99038);
\draw [color=blue     ] (   4.50000 ,   12.99038) -- (   5.00000 ,   13.85641);
\draw [color=blue     ] (   5.00000 ,   13.85641) -- (   5.50000 ,   14.72243);
\draw [color=blue     ] (   5.50000 ,   14.72243) -- (   5.00000 ,   15.58846);
\draw [color=blue     ] (   5.00000 ,   15.58846) -- (   5.50000 ,   16.45448);
\draw [color=blue     ] (   5.50000 ,   16.45448) -- (   6.00000 ,   17.32051);
\draw [color=blue     ] (   6.00000 ,   17.32051) -- (   5.50000 ,   18.18653);
\draw [color=blue     ] (   5.50000 ,   18.18653) -- (   6.00000 ,   19.05256);
\draw [color=blue     ] (   6.00000 ,   19.05256) -- (   6.50000 ,   19.91858);
\draw [color=blue     ] (   6.50000 ,   19.91858) -- (   7.00000 ,   20.78461);
\draw [color=blue     ] (   7.00000 ,   20.78461) -- (   7.50000 ,   21.65064);
\draw [color=blue     ] (   7.50000 ,   21.65064) -- (   8.00000 ,   22.51666);
\draw [color=blue     ] (   8.00000 ,   22.51666) -- (   8.50000 ,   23.38269);
\draw [color=blue     ] (   8.50000 ,   23.38269) -- (   8.00000 ,   24.24871);
\draw [color=blue     ] (   8.00000 ,   24.24871) -- (   8.50000 ,   25.11474);
\draw [color=blue     ] (   8.50000 ,   25.11474) -- (   9.00000 ,   25.98076);
\draw [color=blue     ] (   9.00000 ,   25.98076) -- (   9.50000 ,   26.84679);
\draw [color=blue     ] (   9.50000 ,   26.84679) -- (  10.00000 ,   27.71281);
\draw [color=blue     ] (  10.00000 ,   27.71281) -- (  10.50000 ,   28.57884);
\draw [color=blue     ] (  10.50000 ,   28.57884) -- (  11.00000 ,   29.44486);
\draw [color=blue     ] (  11.00000 ,   29.44486) -- (  11.50000 ,   30.31089);
\draw [color=blue     ] (  11.50000 ,   30.31089) -- (  12.00000 ,   31.17691);
\draw [color=blue     ] (  12.00000 ,   31.17691) -- (  12.50000 ,   32.04294);
\draw [color=blue     ] (  12.50000 ,   32.04294) -- (  13.00000 ,   32.90897);
\draw [color=blue     ] (  13.00000 ,   32.90897) -- (  13.50000 ,   33.77499);
\draw [color=blue     ] (  13.50000 ,   33.77499) -- (  13.00000 ,   34.64102);
\draw [color=blue     ] (  13.00000 ,   34.64102) -- (  13.50000 ,   35.50704);
\draw [color=blue     ] (  13.50000 ,   35.50704) -- (  14.00000 ,   36.37307);
\draw [color=blue     ] (  14.00000 ,   36.37307) -- (  14.50000 ,   37.23909);
\draw [color=blue     ] (  14.50000 ,   37.23909) -- (  15.00000 ,   38.10512);
\draw [color=blue     ] (  15.00000 ,   38.10512) -- (  15.50000 ,   38.97114);
\draw [color=blue     ] (  15.50000 ,   38.97114) -- (  16.00000 ,   39.83717);
\draw [color=blue     ] (  16.00000 ,   39.83717) -- (  16.50000 ,   40.70319);
\draw [color=blue     ] (  16.50000 ,   40.70319) -- (  17.00000 ,   41.56922);
\draw [color=blue     ] (  17.00000 ,   41.56922) -- (  16.50000 ,   42.43524);
\draw [color=blue     ] (  16.50000 ,   42.43524) -- (  17.00000 ,   43.30127);
\draw [color=blue     ] (  17.00000 ,   43.30127) -- (  17.50000 ,   44.16730);
\draw [color=blue     ] (  17.50000 ,   44.16730) -- (  18.00000 ,   45.03332);
\draw [color=blue     ] (  18.00000 ,   45.03332) -- (  18.50000 ,   45.89935);
\draw [color=blue     ] (  18.50000 ,   45.89935) -- (  19.00000 ,   46.76537);
\draw [color=blue     ] (  19.00000 ,   46.76537) -- (  19.50000 ,   47.63140);
\draw [color=blue     ] (  19.50000 ,   47.63140) -- (  20.00000 ,   48.49742);
\draw [color=blue     ] (  20.00000 ,   48.49742) -- (  20.50000 ,   49.36345);
\draw [color=blue     ] (  20.50000 ,   49.36345) -- (  21.00000 ,   50.22947);
\draw [color=blue     ] (  21.00000 ,   50.22947) -- (  21.50000 ,   51.09550);
\node at (  22.00000 ,   51.96152) [diamond    , fill=orange   ] {};
\draw [color=blue     ] (  21.50000 ,   51.09550) -- (  22.00000 ,   51.96152);
\draw [color=blue     ] (  22.00000 ,   51.96152) -- (  22.50000 ,   52.82755);
\draw [color=blue     ] (  22.50000 ,   52.82755) -- (  23.00000 ,   53.69358);
\draw [color=blue     ] (  23.00000 ,   53.69358) -- (  23.50000 ,   54.55960);
\draw [color=blue     ] (  23.50000 ,   54.55960) -- (  24.00000 ,   55.42563);
\draw [color=blue     ] (  24.00000 ,   55.42563) -- (  24.50000 ,   56.29165);
\draw [color=blue     ] (  24.50000 ,   56.29165) -- (  25.00000 ,   57.15768);
\draw [color=blue     ] (  25.00000 ,   57.15768) -- (  25.50000 ,   58.02370);
\draw [color=blue     ] (  25.50000 ,   58.02370) -- (  26.00000 ,   58.88973);
\draw [color=blue     ] (  26.00000 ,   58.88973) -- (  26.50000 ,   59.75575);
\draw [color=blue     ] (  26.50000 ,   59.75575) -- (  27.00000 ,   60.62178);
\draw [color=blue     ] (  27.00000 ,   60.62178) -- (  27.50000 ,   61.48780);
\draw [color=blue     ] (  27.50000 ,   61.48780) -- (  28.00000 ,   62.35383);
\draw [color=blue     ] (  28.00000 ,   62.35383) -- (  28.50000 ,   63.21985);
\draw [color=blue     ] (  28.50000 ,   63.21985) -- (  29.00000 ,   64.08588);
\draw [color=blue     ] (  29.00000 ,   64.08588) -- (  30.00000 ,   64.08588);
\draw [color=blue     ] (  30.00000 ,   64.08588) -- (  30.50000 ,   64.95191);
\draw [color=blue     ] (  30.50000 ,   64.95191) -- (  31.50000 ,   64.95191);
\draw [color=blue     ] (  31.50000 ,   64.95191) -- (  32.50000 ,   64.95191);
\draw [color=blue     ] (  32.50000 ,   64.95191) -- (  33.50000 ,   64.95191);
\draw [color=blue     ] (  33.50000 ,   64.95191) -- (  34.00000 ,   65.81793);
\draw [color=blue     ] (  34.00000 ,   65.81793) -- (  35.00000 ,   65.81793);
\draw [color=blue     ] (  35.00000 ,   65.81793) -- (  35.50000 ,   66.68396);
\draw [color=blue     ] (  35.50000 ,   66.68396) -- (  36.50000 ,   66.68396);
\draw [color=blue     ] (  36.50000 ,   66.68396) -- (  37.50000 ,   66.68396);
\draw [color=blue     ] (  37.50000 ,   66.68396) -- (  38.00000 ,   67.54998);
\draw [color=blue     ] (  38.00000 ,   67.54998) -- (  39.00000 ,   67.54998);
\draw [color=blue     ] (  39.00000 ,   67.54998) -- (  40.00000 ,   67.54998);
\draw [color=blue     ] (  40.00000 ,   67.54998) -- (  41.00000 ,   67.54998);
\draw [color=blue     ] (  41.00000 ,   67.54998) -- (  42.00000 ,   67.54998);
\draw [color=blue     ] (  42.00000 ,   67.54998) -- (  43.00000 ,   67.54998);
\draw [color=blue     ] (  43.00000 ,   67.54998) -- (  43.50000 ,   68.41601);
\draw [color=blue     ] (  43.50000 ,   68.41601) -- (  44.50000 ,   68.41601);
\draw [color=blue     ] (  44.50000 ,   68.41601) -- (  45.50000 ,   68.41601);
\draw [color=blue     ] (  45.50000 ,   68.41601) -- (  46.50000 ,   68.41601);
\draw [color=blue     ] (  46.50000 ,   68.41601) -- (  47.50000 ,   68.41601);
\draw [color=blue     ] (  47.50000 ,   68.41601) -- (  48.50000 ,   68.41601);
\draw [color=blue     ] (  48.50000 ,   68.41601) -- (  49.50000 ,   68.41601);
\draw [color=blue     ] (  49.50000 ,   68.41601) -- (  50.50000 ,   68.41601);
\draw [color=blue     ] (  50.50000 ,   68.41601) -- (  51.50000 ,   68.41601);
\draw [color=blue     ] (  51.50000 ,   68.41601) -- (  52.50000 ,   68.41601);
\draw [color=blue     ] (  52.50000 ,   68.41601) -- (  53.50000 ,   68.41601);
\draw [color=blue     ] (  53.50000 ,   68.41601) -- (  54.50000 ,   68.41601);
\draw [color=blue     ] (  54.50000 ,   68.41601) -- (  55.50000 ,   68.41601);
\draw [color=blue     ] (  55.50000 ,   68.41601) -- (  56.00000 ,   69.28203);
\draw [color=blue     ] (  56.00000 ,   69.28203) -- (  56.50000 ,   70.14806);
\draw [color=blue     ] (  56.50000 ,   70.14806) -- (  57.50000 ,   70.14806);
\draw [color=blue     ] (  57.50000 ,   70.14806) -- (  58.00000 ,   71.01408);
\draw [color=blue     ] (  58.00000 ,   71.01408) -- (  58.50000 ,   71.88011);
\draw [color=blue     ] (  58.50000 ,   71.88011) -- (  59.50000 ,   71.88011);
\draw [color=blue     ] (  59.50000 ,   71.88011) -- (  60.50000 ,   71.88011);
\draw [color=blue     ] (  60.50000 ,   71.88011) -- (  61.50000 ,   71.88011);
\draw [color=blue     ] (  61.50000 ,   71.88011) -- (  62.50000 ,   71.88011);
\draw [color=blue     ] (  62.50000 ,   71.88011) -- (  63.50000 ,   71.88011);
\draw [color=blue     ] (  63.50000 ,   71.88011) -- (  64.50000 ,   71.88011);
\draw [color=blue     ] (  64.50000 ,   71.88011) -- (  65.50000 ,   71.88011);
\draw [color=blue     ] (  65.50000 ,   71.88011) -- (  66.50000 ,   71.88011);
\draw [color=blue     ] (  66.50000 ,   71.88011) -- (  67.50000 ,   71.88011);
\node at (  67.50000 ,   71.88011) [rectangle , fill=green    ] {};


\node at (  26.50000 ,   23.38269) [rectangle    , fill=red    ] {};
\node at (  32.00000 ,   29.44486) [rectangle    , fill=red    ] {};
\node at ( 72.50000 ,   59.75575) [rectangle    , fill=red    ] {};
\node at (  20.00000 ,   15.58846) [rectangle    , fill=red    ] {};
\node at ( 24.00000 ,   22.51666) [rectangle    , fill=red    ] {};
\node at (  33.00000 ,   51.96152) [rectangle    , fill=red    ] {};
\node at (  61.00000 ,   24.24871) [rectangle    , fill=red    ] {};
\node at ( 33.50000 ,   47.63140) [rectangle    , fill=red    ] {};
\node at (  27.50000 ,   54.55960) [rectangle    , fill=red    ] {};
\node at (  45.00000 ,   46.76537) [rectangle    , fill=red    ] {};
\node at (  32.00000 ,   27.71281) [rectangle    , fill=red    ] {};
\node at (  52.00000 ,   45.03332) [rectangle    , fill=red    ] {};
\node at (  38.50000 ,   30.31089) [rectangle    , fill=red    ] {};
\node at (  58.00000 ,   31.17691) [rectangle    , fill=red    ] {};
\node at (  49.00000 ,   34.64102) [rectangle    , fill=red    ] {};


\end{tikzpicture}
\end{scaletikzpicturetowidth}
\caption{Case 5, larger sensor detection area: Time to target 108 time steps, confusion action at time step 51, PED $\geq 0.95$}
\label{figsqrt}
\end{figure}
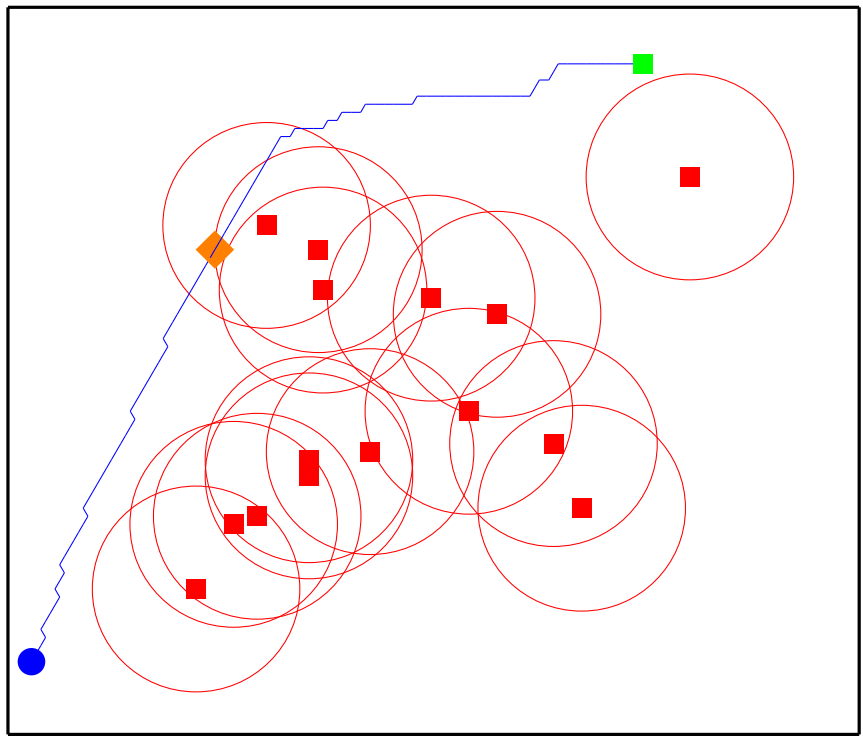

\section{Conclusions} 
\label{sec8}	

In this paper we have considered a problem that arises in military mission planning, namely  
to decide a path for one, or more, agents to reach a target without being detected by enemy sensors.
In this problem agents are not passive, rather they can 
initiate actions which aid evasion, namely knockout, so  completely disable sensors, and confusion, reduce sensor detection probabilities. 

Agent actions were path dependent, so an agent needed to be sufficiently close to a sensor to knock it  out. Agent actions were  also time dependent in  that a time limit was imposed on how long a sensor was knocked out or confused before it reverted back to its original operating state.

We broke the continuous space in which agents move into a discrete space. This enabled the problem to be formulated as a zero-one integer program with linear constraints. A heuristic for the problem based on successive shortest paths was also presented.
Computational results were presented for a number of randomly generated test problems that are made publicly available.

 \clearpage
\newpage
 \pagestyle{empty}
\linespread{1}
\small \normalsize 


\section*{\textbf{Acknowledgments}}
This paper grew out of participation in the workshop: Mathematical Challenges in the Electromagnetic Environment. This workshop was hosted by the Centre for Mathematical Sciences
Cambridge,
 and 
held (virtually) January 8-10 2021.  It was sponsored by the UK Defence Science and Technology Laboratory (Dstl).  The author would accordingly like to acknowledge the support of Dstl, both for this initial workshop and subsequently in developing the work as reported here. 

\vspace{5mm}
\noindent The author also would like to thank the anonymous reviewers for their comments and suggestions to improve the paper.

\bibliography{report}

\end{document}